\documentclass[a4paper,10pt]{amsart}
\usepackage[utf8]{inputenc}

\usepackage[
  margin=30mm,
  marginparwidth=25mm,     
  marginparsep=2mm,       
  bottom=25mm,
  ]{geometry}

\usepackage[bbgreekl]{mathbbol}
\usepackage{amsfonts}
\usepackage{latexsym,amssymb,amsthm,mathrsfs,amsmath,amscd,enumerate,enumitem, amscd,color}
\usepackage{mathrsfs}
\usepackage{tikz}
\usepackage{hyperref}
\usetikzlibrary{plotmarks}
\usepackage{soul}

\DeclareSymbolFontAlphabet{\mathbb}{AMSb}
\DeclareSymbolFontAlphabet{\mathbbl}{bbold}
\DeclareSymbolFontAlphabet{\mathbb}{AMSb}
\DeclareSymbolFontAlphabet{\mathbbl}{bbold}
\DeclareFontEncoding{FMS}{}{}
\DeclareFontSubstitution{FMS}{futm}{m}{n}
\DeclareFontEncoding{FMX}{}{}
\DeclareFontSubstitution{FMX}{futm}{m}{n}
\DeclareSymbolFont{fouriersymbols}{FMS}{futm}{m}{n}
\DeclareSymbolFont{fourierlargesymbols}{FMX}{futm}{m}{n}
\DeclareMathDelimiter{\VERT}{\mathord}{fouriersymbols}{152}{fourierlargesymbols}{147}
\DeclareRobustCommand{\stirling}{\genfrac\{\}{0pt}{}}

\newtheorem{thm}{Theorem}[section]
 \newtheorem{cor}[thm]{Corollary}

\theoremstyle{definition}

 \theoremstyle{remark}

\usepackage{hyperref}

\newcommand{\supp}{\mathop{\mathrm{supp}}}

\numberwithin{equation}{section}
\allowdisplaybreaks

\begin{document}

\footnotetext{Last modification: \today.}

\title[]
 {Littlewood-Paley-Stein theory and Banach spaces in the inverse Gaussian setting}

\author[V. Almeida]{V\'{\i}ctor Almeida}

\author[J.J. Betancor]{Jorge J. Betancor}
\address{V\'{\i}ctor Almeida, Jorge J. Betancor, Juan C. Fari\~na, Lourdes Rodr\'{\i}guez-Mesa\newline
	Departamento de An\'alisis Matem\'atico, Universidad de La Laguna,\newline
	Campus de Anchieta, Avda. Astrof\'isico S\'anchez, s/n,\newline
	38721 La Laguna (Sta. Cruz de Tenerife), Spain}
\email{valmeida@ull.edu.es, jbetanco@ull.es,
jcfarina@ull.es, lrguez@ull.edu.es
}

\author[J.C. Fari\~na]{Juan C. Fari\~na}

\author[L. Rodr\'{\i}guez-Mesa]{Lourdes Rodr\'{\i}guez-Mesa}

\thanks{The authors are partially supported by grant PID2019-106093GB-I00 from the Spanish Government}

\subjclass[2020]{42B20, 42B25, 47B90}

\keywords{Littlewood-Paley functions, inverse Gaussian measure, $q$-uniformly convex, $q$-uniformly smooth and UMD Banach spaces, K\"othe function spaces.}

\date{}


\begin{abstract}
In this paper we consider Littlewood-Paley functions defined by the semigroups associated with the operator $\mathcal{A}=-\frac{1}{2}\Delta-x\nabla$ in the inverse Gaussian setting for Banach valued functions. We characterize the uniformly convex and smooth Banach spaces by using $L^p(\mathbb R^n,\gamma_{-1})$- properties of the $\mathcal{A}$-Littlewood-Paley functions. We also use Littlewood-Paley functions associated with $\mathcal{A}$ to characterize the K\"{o}the function spaces with the UMD property.
\end{abstract}

\maketitle

\section*{}
\begin{flushright}
\textit{....... But there are those who struggle all their lives: These are the indispensable ones.} 

\hfill (Bertolt Brecht)\newline
\newline
\hfill {\bf Dedicated to our friend Jos\'e M\'endez for his retirement}
\end{flushright}

\vspace*{0.5cm} 

\section{Introduction}
Let $\{T_t\}_{t>0}$ be a semigroup of operators on a measure space $(\Omega,\mu)$, $k\in \mathbb N$ and $q \in (1,\infty)$. The Littlewood-Paley function $g_{k,\{T_t\}_{t>0}}^q$ associated with $\{T_t\}_{t>0}$ is defined by
$$
g^q_{k,\{T_t\}_{t>0}}(f)(x)=\left(\int_0^\infty \left| t^k\partial^k_t T_tf(x)\right|^q\frac{dt}{t}\right)^{1/q},\;\;x\in \Omega,
$$
for every $f\in L^p(\Omega,\mu)$, $1\leq p < \infty$.

Littlewood-Paley functions, also called $g$-functions, allow us to obtain equivalent norms in $L^p(\Omega,\mu)$, $1 \leq p <\infty$. This fact makes that $L^p$-boundedness properties of operators connected with the semigroup $\{T_t\}_{t>0}$ can be obtained in a unified way.

Stein developed in \cite{StLP} the Littlewood-Paley theory for symmetric diffusion semigroups. He proved that if $\{T_t\}_{t>0}$ is a symmetric difussion semigroup then, for every $1<p<\infty$ and $k\in \mathbb N$ there exists $C>0$ such that
\begin{equation}\label{(1.1)}
\frac{1}{C}\|f - E_0(f)\|_p \leq \|g_{k,\{T_t\}_{t>0}}^2(f)\|_p \leq C\|f\|_p,\,\, f \in L^p(\Omega,\mu).
\end{equation}
Actually $C$ depends only on $p$ and $k$. Here, $\displaystyle E_0(f)=\lim_{t \rightarrow \infty}T_t(f)$ and $E_0$ is the projection of $L^p(\Omega,\mu)$ on the fixed point space of $\{T_t\}_{t>0}$.

Suppose that $X$ is a Banach space. For every $1\leq p < \infty$ we denote by $L^p_X(\Omega,\mu)$ the $p$-th B\"{o}chner-Lebesgue space. It is wellknown that if $L$ is a bounded operator from $L^p(\Omega,\mu)$ into itself with $1\leq p <\infty$ such that $Lf \geq 0$ when $0 \leq f \in L^p(\Omega,\mu)$, then the tensor extension of $L$ to $L^p(\Omega,\mu) \otimes X$ can be extended to $L^p_X(\Omega,\mu)$ as a bounded operator from $L^p_X(\Omega, \mu)$ into itself that will be also denoted by $L$.

If $\{T_t\}_{t>0}$ is a positive semigroup of operators on $(\Omega,\mu)$, $k \in \mathbb N$ and $q \in (1,\infty)$ the Littlewood-Paley function $g_{k,\{T_t\}_{t>0}}^q$ can be defined on $L^p_X(\Omega,\mu)$ as follows
$$
g_{k,\{T_t\}_{t>0}}^{q,X}(f)(x)=\left(\int_0^\infty\left\|t^k\partial_t^kT_t(f)(x)\right\|_X^q\frac{dt}{t}\right)^{1/q},\;\; x \in \mathbb R^n,
$$
for every $f \in L^p_X(\Omega,\mu)$ with $1\leq p < \infty$. Here $\|\cdot\|_X$ denotes the norm in $X$. The question in this point is whether the inequalities in (\ref{(1.1)}) are still valid in the Banach valued context.

Let $\mathbb T$ be the unit circle. For every $f \in L^1_X(\mathbb T)$ we define
$$
G^X(f)(z)= \left(\int_0^1((1-r)\| \nabla P(f)(rz)\|_X)^2\frac{dr}{1-r}\right)^{1/2},\;\;z\in \mathbb T.
$$
Here $P$ represents the Poisson integral in $\mathbb T$ and $\nabla$ is the gradient. It is known that there exists $C>0$ such that, for every $f\in L^p_X(\mathbb T)$, with $1<p<\infty$,
\begin{equation}\label{1.2}
\frac{1}{C}\|f\|_{L^p_X(\mathbb T)}\leq |\hat{f}(0)| + \|G^X(f)\|_{L^p(\mathbb T)} \leq C\|f\|_{L^p_X(\mathbb T)},
\end{equation}
where $\hat{f}(0) = \int_{\mathbb T} f(w) dw$, if and only if $X$ is isomorphic to a Hilbert space. Xu (\cite{Xu1}) considered $g$-functions defined as follows. Let $ q \in (1,\infty)$, we define
$$
G^{q,X}(f)(z)=\left(\int_0^1((1-r)\|\nabla P(f)(rz)\|_X)^q\frac{dr}{1-r}\right)^{1/q},\;\; z \in \mathbb T,
$$
for every $f\in L^1_X(\mathbb T)$. Xu characterized the Banach spaces $X$ for which one of the inequalities in (\ref{1.2}) holds.

We recall the martingale type and cotype introduced by Pisier (\cite{Pi1}). A Banach space $X$ is said to be of martingale cotype (respectively, martingale type) $q\in(1,\infty)$ when there exists $C>0$ such that for every martingale $(M_n)_{n\in \mathbb N}$ defined on some probability space with values in $X$ the following inequality holds
$$
\sum_{n\in \mathbb N} \mathbb E\|M_n-M_{n-1}\|_X^q \leq C \sup_{n \in \mathbb N} \mathbb E \|M_n\|_X^q
$$

$$(\mbox{respectively,}\;
\sup_{n \in \mathbb N} \mathbb E \|M_n\|_X^q\leq C \sum_{n\in \mathbb N} \mathbb E\|M_n-M_{n-1}\|_X^q) .
$$
Here, $\mathbb E$ denotes the expectation and $M_{-1}=0$. If $X$ has martingale cotype (respectively, type) $q$ then $q\geq 2$ (respectively $q\in (1,2]$).

We define, for every martingale $M=(M_n)_{n \in \mathbb N}$ with values in $X$,
$$
S_{q,X}(M)=\left( \sum_{n\in \mathbb N} \|M_n-M_{n-1}\|_X^q \right)^{1/q},
$$
with $1<q<\infty$. $S_{q,X}$ can be seen as a martingale analogue of $G^{q,X}$. As it is commented by Xu (\cite[p. 208]{Xu1}) $X$ if of martingale cotype (respectively, type) $q$ if and only if for every (respectively, for some) $1<p<\infty$ there exists $C>0$ such that for every $L^p$-martingale $M=(M_n)_{n\in \mathbb N}$ with values in $X$ the following property holds
$$
\mathbb E\left[S_{q,X}(M)\right]^p \leq C \sup_{n\in \mathbb N} \mathbb E\|M_n\|_X^p
$$

$$ (\mbox{respectively},\;
\sup_{n\in \mathbb N} \mathbb E\|M_n\|_X^p \leq C \mathbb E\left[S_{q,X}(M)\right]^p).
$$

The vector valued harmonic analysis is closely connected with the geometry of Banach spaces. The modulus of convexity of $X$ is defined by
$$
\delta_X(\varepsilon)= \inf \left\{1-\left\|\frac{a+b}{2}\right\|_X: a,b \in X,\,\|a\|_X=\|b\|_X=1,\, \|a-b\|_X=\varepsilon\right\},\quad 0<\varepsilon <2,
$$
and the modulus of smoothness of $X$ is defined by
$$
\rho_X(t) = \sup \left\{\frac{\|a+tb\|_X+\|a-tb\|_X}{2}-1: a,b \in X,\,\|a\|_X=\|b\|_X=1,\, \right\},\quad t>0.
$$
$X$ is said to be uniformly convex when $\delta_X(\varepsilon)>0$, for every $0<\varepsilon<2$, and to be uniformly smooth when $\displaystyle \lim_{t\rightarrow 0^+}\frac{\rho_X(t)}{t}=0$. If $q\in (1,\infty)$ we say that $X$ is $q$-uniformly convex (respectively, $q$-uniformly smooth) when there exists $C>0$ such that $\delta_X(\varepsilon) \geq C\varepsilon^q$, $0<\varepsilon < 2$ (respectively, $\rho_X(t) \leq Ct^q,\; t>0$).

The convexity and smoothness  of a Banach space can be characterized by using $G^{q,X}$ and $S_{q,X}$.

\noindent {\bf Theorem A} (\cite[Chapter 10]{Pi2} and \cite[Theorem 3.1 and Corollary 3.2]{Xu1}). Let $X$ be a Banach space and $1<q_1 \leq 2 \leq q_2<\infty$. The following assertions are equivalent.
\begin{enumerate}
\item[(a)] There exists a norm $\VERT \cdot \VERT$ on $X$ that defines the original topology of $X$ and such that $(X;\VERT  \cdot \VERT)$ is $q_2$-uniformly convex (respectively, $q_1$-uniformly smooth).
\item[(b)] $X$ is of martingale cotype $q_2$ (respectively, type $q_1$).
\item[(c)] For every (equivalently, for some) $1<p<\infty$, there exists $C>0$ such that
$$
\|G^{q_2,X}(f)\|_{L^p(\mathbb T)} \leq C\|f\|_{L^p_X(\mathbb T)},\;\; f \in L^p_X(\mathbb T),
$$

$$
(\mbox{respectively},\;\|f\|_{L^p_X(\mathbb T)}\leq C(\|\widehat f(0)\|+ \|G^{q_1,X}(f)\|_{L^p(\mathbb T)} ).
$$
\end{enumerate}
\hfill $\square$

This result also holds when the Poisson integral on $\mathbb T$ is replaced by the Poisson integral in $\mathbb{R}^n \times(0,\infty)$. In this last case the projection $E_0=0$ and the term containing $\widehat f(0)$ does not appear.

Theorem A was extended to symmetric diffusion semigroups by Mart{\'\i}nez, Torrea and Xu (\cite{MTX}) and Xu (\cite{Xu2}). According to Stein (see \cite{StLP}) a uniparametric family $\{T_t\}_{t>0}$ of operators defined on $\displaystyle\bigcup_{1\leq p \leq \infty}L^p(\Omega,\mu)$ is called a symmetric diffusion semigroup when the following properties hold
\begin{itemize}
\item[(i)] $T_{t}$ is a contraction on $L^p(\Omega,\mu)$, for every $1\leq p\leq\infty$ and $t>0$;
\item[(ii)] $T_{t+s}=T_{t}T_{s}$, on $L^p(\Omega,\mu)$, for every $1\leq p\leq\infty$ and $t,s>0$;
\item[(iii)] $\lim\limits_{t\to 0^{}+}T_{t}f=f$, in $L^{2}(\Omega,\mu)$, for every $f\in L^{2}(\Omega,\mu)$;
\item[(iv)] $T_{t}$ is selfadjoint on $L^{2}(\Omega,\mu)$;
\item[(v)] $T_{t}$ is positive preserving, that is, $T_{t}f\geq 0$, $t>0$, provided that $f \geq 0$;
\item[(vi)] $T_{t}$ is Markovian, that is, $T_{t}1=1$, for every $t>0$.
\end{itemize}
 The Poisson subordinated semigroup $\{P^t\}_{t>0}$ of $\{T_t\}_{t>0}$ is defined in the following way
 $$
 P^t(f)=\frac{1}{\sqrt \pi}\int_0^\infty \frac{e^{-s}}{\sqrt s}T_{\frac{t^2}{4s}}(f)\,ds,\quad t>0.
 $$
 $\{P^t\}_{t>0}$ is also a symmetric diffusion semigroup. The Euclidean heat semigroup in $\mathbb R^n$ is the first example of symmetric diffusion semigroup. The Euclidean Poisson semigroup is the Poisson subordinated semigroup of the Euclidean heat semigroup.

 \noindent {\bf Theorem B} (\cite[Theorems 2.1 and 2.2]{MTX}). Let $X$ be a Banach space. Suppose that $\{T_t\}_{t>0}$ is a symmetric diffusion semigroup and that $\{P^t\}_{t>0}$ is the Poisson subordinated semigroup of $\{T_t\}_{t>0}$.
 \begin{enumerate}
 \item[(a)] If $q \geq 2$ (respectively, $q \in (1,2]$) and $X$ is of martingale cotype $q$ (respectively, martingale type $q$) then, for every $1 < p< \infty$, there exists $C>0$ such that
     \begin{equation}\label{1.3}
     \left\| g^{q,X}_{1, \{P^t\}_{t>0}}(f) \right\|_{L^p(\Omega,\mu)} \leq C\|f\|_{L^p_X(\Omega,\mu)},\quad f\in L^p_X(\Omega,\mu),
     \end{equation}

\begin{equation}\label{1.4}
    (\mbox{respectively},\; \|f\|_{L^p_X(\Omega,\mu)}\leq C\left( \|E_0(f)\|_{L^p_X(\Omega,\mu)} + \left\| g^{q,X}_{1, \{P^t\}_{t>0}}(f) \right\|_{L^p(\Omega,\mu)}\right),\;\; f\in L^p_X(\Omega,\mu)).
     \end{equation}
\item[(b)] If $q \geq 2$ (respectively, $q \in (1,2]$) and (\ref{1.3}) (respectively, (\ref{1.4})) holds when $\{P^t\}_{t>0}=\{P_t\}_{t>0}$ is the Euclidean Poisson semigroup for some $1<p<\infty$, then $X$ is of martingale cotype $q$ (respectively, martingale type $q$).
\end{enumerate}

\noindent {\bf Theorem C} (\cite[Theorem 2]{Xu2}) Let $X$ be a Banach space. Suppose that $\{T_t\}_{t>0}$ is a symmetric diffusion semigroup. If $q \geq 2$ (respectively, $q \in (1,2]$) and X is of martingale cotype $q$ (respectively, martingale type $q$), then, for every $1<p<\infty$, there exists $C>0$ such that
$$
\left\|g_{k,\{T_t\}_{t>0}}^{q,X}(f) \right\|_{L^p(\Omega,\mu)}\leq \|f\|_{L^p_X(\Omega,\mu)},\;\; f \in L^p_X(\Omega,\mu),
$$

$$
(\mbox{respectively},\;\|f\|_{L^p_X(\Omega,\mu)} \leq C\left(\|E_0(f)\|_{L^p_X(\Omega,\mu)}+ \left\|g_{k,\{T_t
\}_{t>0}}^{q,X}(f) \right\|_{L^p(\Omega,\mu)}\right),\quad f\in L^p_X(\Omega,\mu)).
$$

The result in Theorem C had been previously established for the Euclidean heat semigroup in $\mathbb R^n$ by Hyt\"onen and Naoir (\cite{HN}). A version of the Theorem B, (b), when the Euclidean Poisson semigroup is replaced by Euclidean heat semigroup was proved in \cite[Theorems 1.5 and 1.7]{BFGM}.

Markovian property for the semigroup is crucial to prove the above results. Characterizations of $q$-uniformly convex and smooth Banach spaces as  the ones established above but using semigroups of operators without Markovian property (associated with Laguerre and Hermite operators) were proved in \cite{BFGM}, \cite{BFMT} and \cite{BFRST}.

Let $\alpha >0$. We choose $m\in \mathbb N$ such that $m-1 \leq \alpha < m$. If $\phi \in C^m(0,\infty)$ the Weyl $\alpha$-derivative $\partial_t^\alpha \phi$ of $\phi$ is defined by
$$
\partial_t^\alpha \phi(t)= \frac{1}{\Gamma(m-\alpha)}\int_t^\infty (\partial^m_u\phi)(u)(u-t)^{m-\alpha - 1}du, \;\; t\in (0,\infty),
$$
provided that the last integral exists.

Assume that $\{T_t\}_{t>0}$ is a symmetric diffusion semigroup on $(\Omega,\mu)$. For every $\alpha >0$, $f\in L^p(\Omega,\mu)$, $1 < p <\infty$ and $t>0$, we have that $\int_t^{\infty}\|\partial_u^m T_u(f)\|_{L^p(\Omega)}(u-t)^{m-\alpha-1}du <\infty$, with $m-1\leq \alpha < m$ (\cite[p. 11]{BFGM}). We define, for each $f\in L^p(\Omega,\mu)$, $1<p<\infty$,
$$
\partial_t^\alpha T_t(f)=\frac{1}{\Gamma(m - \alpha)} \int_t^\infty  \partial_u^mT_u(f)(u-t)^{m-\alpha-1}du, \;\; t>0.
$$

We now consider the Littlewood-Paley function with fractional derivatives. If $\alpha >0$ and $1<q<\infty$, the fractional Littlewood-Paley functions $g^{q,X}_{\alpha,\{T_t\}_{t>0}}$ is defined by
$$
g^{q,X}_{\alpha,\{T_t\}_{t>0}}(f)(x)= \left(\int_0^\infty \|t^\alpha \partial^\alpha_tT_t(f)(x)\|_X^q \frac{dt}{t}\right)^{1/q},\;\; x\in \mathbb R^n,
$$
for every $f\in L^p(\Omega,\mu)$, $1<p < \infty$.

Torrea and Zhang (\cite[Theorems A and B]{TZ}) proved a version of Theorem B involving fractional Littlewood-Paley functions. In \cite[Theorems 1.5 and 1.7]{BFGM} Theorem C is generalized by using fractional $g$-functions.

The  measure $d\gamma(x)=\frac{e^{-|x|^2}}{\pi^{n/2}}dx$, where $dx$ denotes the Lebesgue measure in $\mathbb R^n$, is named Gaussian measure in $\mathbb R^n$. Harmonic analysis associated with the Gaussian-measure was began by Muckenhoupt (\cite{Mu1}). The Ornstein-Uhlenbeck operator $\mathcal{L}$ is defined by
$$
\mathcal{L}\phi=-\frac{\Delta}{2} \phi +x\cdot\nabla \phi, \;\; \phi \in C_c^\infty(\mathbb R^n).
$$
For every $k\in \mathbb N$ we denote by $H_ k$ the $k$-th Hermite polynomial given by
$$
H_k(x) = (-1)^k e^{x^2}\frac{d^k}{dx^k}e^{-x^2},\;\; x \in \mathbb R.
$$
If $k=(k_1,\ldots,k_n) \in \mathbb N^n$ we define
$$
H_k(x) = (-1)^{|k|} \prod_{i=1}^nH_{k_i}(x_i),\;\; x=(x_1,\ldots,x_n) \in \mathbb R^n,
$$
where $|k|=k_1+\ldots +k_n$.

The spectrum of $\mathcal L$ in $L^p(\mathbb R^n,d\gamma)$, $1<p<\infty$, is the discrete set $\mathbb N$ and we have that, for every $k \in \mathbb N^n$,
$$
\mathcal LH_k= |k|H_k.
$$
The Ornstein-Uhlenbeck operator $-\mathcal L$ generates a symmetric diffusion semigroup $\{T^{\mathcal L}_t\}_{t>0}$. Operators defined by $\{T^{\mathcal L}_t\}_{t>0}$ (maximal operators, Riesz transforms, Littlewood-Paley functions, multipliers,\ldots) associated with harmonic analysis in the gaussian setting have been studied by several authors (see \cite{FGS}, \cite{GCMST1}, \cite{GCMST2}, \cite{HTV}, \cite{Pe}, \cite{PS}, \cite{Ur}, and the references therein). In \cite[Theorem 1.12]{HTV} and \cite[Theorems 6.1 and 6.2]{MTX} $q$-uniformly convex and smooth Banach spaces by using Littlewood-Paley functions defined by the Poisson semigroup associated with the Ornstein-Uhlenbeck operator.

We call the inverse gaussian measure on $\mathbb R^n$ to the measure $d\gamma_{-1}(x)= \pi^{n/2}e^{|x|^2} dx$. We will write $\gamma_{-1}$ to refer us to the inverse Gauss measure. $\gamma_{-1}$ is not doubling. The study of the harmonic analysis in the inverse Gaussian setting was began by Salogni (\cite{Sa}).

We consider the differential operator
$$
\mathcal A \phi= -\frac{1}{2} \Delta\phi-x\cdot \nabla\phi,
$$
for every $\phi \in C_c^\infty(\mathbb R^n)$. $\mathcal A$ is essentially selfadjoint with respect to $L^2(\mathbb R^n,\gamma_{-1})$. We continue denoting by $\mathcal A$ to the closure of this operator.

For every $k \in \mathbb N^n$ we define
$$
\widetilde{H}_k(x)= e^{-|x|^2}H_k(x), \;x\in \mathbb R^n.
$$
We have that
$$
\mathcal A\widetilde{H}_k=(|k| +n)\widetilde{H}_k,\; k \in \mathbb N^n,
$$
and the spectrum of $\mathcal A$ in $L^p(\mathbb R^n,\gamma_{-1})$ is the set $\{k+n,\;k \in \mathbb N\}$, for every $1<p<\infty$. The operator $-\mathcal A$ generates the semigroup of operators $\{T^\mathcal A_t\}_{t>0}$ in $L^2(\mathbb R^n,\gamma_{-1})$ defined by
$$
T^\mathcal A_t (f)= \sum_{k \in \mathbb N^n}e^{-(|k|+n)t} c_k(f) \widetilde{H}_k,\;\; f \in L^2(\mathbb R^n,\gamma_{-1}),
$$
where $c_k(f)=\|\widetilde{H}_k\|_{L^2(\mathbb{R}^n,\gamma_{-2})}^{-1}\int_{\mathbb R^n}f(y) \widetilde{H}_k(y)d\gamma_{-1}(y)$, $f \in L^2(\mathbb R^n,\gamma_{-1})$ and $k \in \mathbb N^n$.
By using the Mehler's formula we can write, for every $t>0$,
\begin{equation}\label{1.5}
T_t^\mathcal Af(x) = \int_{\mathbb R^n} T_t^\mathcal A(x,y) f(y)dy, \;\; f \in L^2(\mathbb R^n, \gamma_{-1}),
\end{equation}
being
$$
T_t^\mathcal A(x,y)= \frac{1}{\pi^{n/2}} e^{-nt}\frac{e^{-\frac{|x-e^{-t}y|^2}{1-e^{-2t}}}}{(1-e^{-2t})^{\frac{n}{2}}}, \;\; x,y \in \mathbb R^n \;\;\text{and}\;\;t>0.
$$
Thus, $\{T_t^\mathcal A\}_{t>0}$ is a symmetric diffusion semigroup and, for every $t>0$ and $1\leq p<\infty$, $T^\mathcal A_t$ admits the integral representation (\ref{1.5}) in $L^p(\mathbb R^n,\gamma_{-1})$.

Salogni (\cite[Chapter 3]{Sa}) studied $L^p(\mathbb R^n,\gamma_{-1})$-boundedness properties of the maximal operator
$$
T_*^\mathcal A f= \sup_{t>0} \left|T_t^\mathcal A f\right|
$$
and of a certain spectral multiplier associated with the operator $\mathcal A$.

In \cite{Bru} Bruno proved $L^p(\mathbb R,\gamma_{-1})$-boundedness and $L^p(\mathbb R,\gamma_{-1})$-unboundedness results for the purely imaginary powers and the first order Riesz transforms associated with the translate operator $\mathcal A+ \lambda I$, $\lambda > 0$, from certain Hardy type space adaptaded to $\gamma_{-1}$ to  $L^1(\mathbb R^n,\gamma_{-1}$) and from $L^1(\mathbb R^n,\gamma_{-1})$ into $L^{1,\infty}(\mathbb R^n,\gamma_{-1})$. Recently, Bruno and Sj\"ogren (\cite{BrSj}) have proved that the first and second order Riesz transform defined by $\mathcal A$ are bounded from $L^1(\mathbb R^n,\gamma_{-1})$ into $L^{1,\infty}(\mathbb R^n,\gamma_{-1})$ but Riesz transform with order greater than two are not bounded from $L^1(\mathbb R^n,\gamma_{-1})$ into $L^{1,\infty}(\mathbb R^n,\gamma_{-1})$.

In \cite{BCdL} maximal operators given by
$$
T^\mathcal A_{*,k} (f)= \sup_{t>0} \left|t^k\partial^k_tT_t^\mathcal A(f)\right|, \;\;k \in \mathbb N,
$$
are studied in $L^p(\mathbb R^n,\gamma_{-1})$, $1\leq p<\infty$. Also, the Banach K\"othe function spaces with the Hardy-Littlewood property are characterized by using the maximal operator $T_{*,k}^\mathcal A$, $k \in \mathbb N$. The UMD property for Banach spaces is characterized by using Riesz transforms and imaginary powers of $\mathcal A$ in \cite[Theorems 1.4, 1.5 and 1.6]{BR} where it is established that higher order Riesz transforms in the inverse Gaussian setting are bounded from $L^p(\mathbb R^n, \gamma_{-1})$ into itself, for every $1<p<\infty$.

Our objective in this paper is to study the $L^p(\mathbb R^n,\gamma_{-1})$-boundedness properties of the Littlewood-Paley functions defined by $\{T_t^{\mathcal A}\}_{t>0}$ and also by the Poisson semigroup $\{P_t^{\mathcal A}\}_{t>0}$ subordinated to $\{T_t^{\mathcal A}\}_{t>0}$. We characterize $q$-uniformly convex and smooth Banach spaces by using these Littlewood-Paley functions. Also, the UMD property for K\"othe function spaces is characterized by $g$-functions in the inverse Gaussian setting.

Let $\beta >0$, $q \in (1,\infty)$ and $k=(k_1,\ldots,k_n) \in \mathbb N^n\setminus\{0\}$. Assume that $X$ is a Banach space. We define the Littlewood-Paley g-functions by
$$
g^{q,X}_{\beta,k,\{T_t^\mathcal A\}_{t>0}}(f)(x)=\left(\int_0^\infty \left\|t^{|k|+\beta}\partial _x^k\partial^\beta_tT_t^\mathcal A(f)(x)\right\|_X^q\frac{dt}{t}\right)^{1/q},\;\;x \in \mathbb R^n.
$$
Here and in the sequel $\partial _x^k=\frac{\partial^{|k|}}{\partial x_1^{k_1}\ldots \partial x_n^{k_n}}$, $x=(x_1,...,x_n)\in \mathbb{R}^n$. When $X=\mathbb C$ we write $g^{q}_{\beta,k,\{T_t^\mathcal A\}_{t>0}}$. In a similar way we define $g^{q,X}_{\beta,k,\{P_t^\mathcal A\}_{t>0}}$.
\begin{thm}\label{Th1.1} Let $\beta >0$, $q \in [2,\infty)$ and $k \in \mathbb N^n$.
\begin{enumerate}
\item[(i)] $g^{q}_{\beta,k,\{T_t^\mathcal A\}_{t>0}}$ is bounded from $L^p(\mathbb R^n,\gamma_{-1})$ into itself, for every $1<p<\infty$. When $0<\beta \leq 1$, $g_{\beta, 0,\{T_t^\mathcal{A}\}_{t>0}}^q$ is bounded from $L^1(\mathbb{R}^n,\gamma_{-1})$ into $L^{1,\infty }(\mathbb{R}^n,\gamma_{-1})$.
\item[(ii)]$g^{q}_{\beta,k,\{P_t^\mathcal A\}_{t>0}}$ is bounded from $L^p(\mathbb R^n,\gamma_{-1})$ into itself, for every $1<p<\infty$.
Furthermore, $g^{q}_{\beta,k,\{P_t^\mathcal A\}_{t>0}}$ is bounded from $L^1(\mathbb R^n,\gamma_{-1})$ into $L^{1,\infty}(\mathbb R^n,\gamma_{-1})$, provided that (respectively, if and only if) $|k| \leq 2$ ({respectively, when $\beta \geq 1$}).
\end{enumerate}
\end{thm}
We characterize $q$-uniformly convex and smooth Banach spaces by using $g$-functions associated with $\mathcal A$.

\begin{thm}\label{Th1.2} Let $X$ be a Banach space, $2\leq q<\infty$, $1<p<\infty$ and $\beta >0$. The following assertions are equivalent.
\begin{enumerate}
\item[(i)] There exists a norm $\VERT \cdot \VERT$ on $X$ that defines the original topology of $X$ such that $(X,\VERT \cdot \VERT)$ is $q$-uniformly convex.
\item[(ii)] $g^{q,X}_{\beta,0,\{T_t^\mathcal A\}_{t>0}}$ is bounded from $L^p_X(\mathbb R^n,\gamma_{-1})$ into $L^p(\mathbb R^n,\gamma_{-1})$.
\item[(iii)] $g^{q,X}_{\beta,0,\{P_t^\mathcal A\}_{t>0}}$ is bounded from $L^p_X(\mathbb R^n,\gamma_{-1})$ into $L^p(\mathbb R^n,\gamma_{-1})$.
\item[(iv)] $g^{q,X}_{\beta,0,\{P_t^\mathcal A\}_{t>0}}$ is bounded from $L^1_X(\mathbb R^n,\gamma_{-1})$ into $L^{1,\infty}(\mathbb R^n,\gamma_{-1})$.
\end{enumerate}
\end{thm}

\begin{thm}\label{Th1.3} Let $X$ be a Banach space, $1 \leq q<2$, $1<p<\infty$ and $\beta >0$. The following assertions are equivalent.
\begin{enumerate}
\item[(i)] There exists a norm $\VERT \cdot \VERT$ on $X$ that defines the original topology of $X$ such that $(X,\VERT \cdot \VERT)$ is $q$-uniformly smooth.
\item[(ii)] There exists $C>0$ such that
$$
\|f\|_{L^p_X(\mathbb R^n,\gamma_{-1})} \leq C \left\|g^{q,X}_{\beta,0,\{T_t^\mathcal A\}_{t>0}}(f)\right\|_{L^p(\mathbb R^n,\gamma_{-1})},\;\; f \in L^p_X(\mathbb R^n,\gamma_{-1}).
$$
\item[(iii)] There exists $C>0$ such that
$$
\|f\|_{L^p_X(\mathbb R^n,\gamma_{-1})} \leq C \left\|g^{q,X}_{\beta,0,\{P_t^\mathcal A\}_{t>0}}(f)\right\|_{L^p(\mathbb R^n,\gamma_{-1})},\;\; f \in L^p_X(\mathbb R^n,\gamma_{-1}).
$$
\end{enumerate}
\end{thm}
As a consequence of Theorem \ref{Th1.2} and \ref{Th1.3} we can deduce the following characterization of the Hilbert spaces.

\begin{cor}\label{cor1.4}
Let $X$ be a Banach space, $1<p<\infty$ and $\beta >0$. The following assertions are equivalent.
\begin{enumerate}
\item[(i)] There exists a norm $\VERT \cdot \VERT$ associated with an inner product in $X$ that defines the original topology on $X$ such that $(X,\VERT \cdot \VERT)$ is a Hilbert space.
\item[(ii)] $\frac{1}{C} \|f\|_{L_X^p(\mathbb R^n,\gamma_{-1})}\leq \left\|g^{2,X}_{\beta,0,\{T_t^\mathcal A\}_{t>0}}(f)\right\|_{L^p(\mathbb R^n,\gamma_{-1})}\leq C\|f\|_{L^p_X(\mathbb R^n,\gamma_{-1})},  \;\; f\in L^p_X(\mathbb R^n, \gamma_{-1})$.
\item[(iii)] $\frac{1}{C} \|f\|_{L_X^p(\mathbb R^n,\gamma_{-1})}\leq \left\|g^{2,X}_{\beta,0,\{P_t^\mathcal A\}_{t>0}}(f)\right\|_{L^p(\mathbb R^n,\gamma_{-1})}\leq C\|f\|_{L^p_X(\mathbb R^n,\gamma_{-1})},  \;\; f\in L^p_X(\mathbb R^n, \gamma_{-1})$.
    \end{enumerate}
\end{cor}

Let $(\Omega, \Sigma,\mu)$ be a complete $\sigma$-finite measure space. A Banach space $X$ consisting of equivalence classes modulus equality almost everywhere with respect to $\mu$ of locally integrable real functions on $\Omega$ is said to be a K\"othe function space when the following properties hold:
\begin{enumerate}
\item[(a)] If $|f(w)| \leq |g(w)|$, for almost all $w \in \Omega$, $f$ is measurable and $g \in X$, then $f \in X$ and $\|f\|_X \leq \|g\|_X$.
\item[(b)] For every $E\in \Sigma$ being $\mu(E) < \infty$, the characteristic function $\mathcal{X}_E$ of $E$ is in $X$.
\end{enumerate}

Each K\"othe function space endowed with the natural order is a Banach lattice. This lattice is $\sigma$-order complete. Furthermore, if $X$ is an order continuous Banach lattice having a weak unity, then $X$ is order isometric to a K\"othe function space (\cite[Theorem1.b.14]{LT}). A Banach lattice $E$ is called order continuous provided that every decreasing net in $E$ whose infimum is zero is norm-convergent to zero. We refer \cite{LT} for the main facts about Banach lattices.

We say that a Banach space has the UMD-property when for every (equivalently, for some) $1<p<\infty$ there  exists $C>0$ such that
$$
\left\|\sum^n_{k=1}\varepsilon_kd_k\right\|_{L^p_X(\Omega,\mu)} \leq C \left\|\sum^n_{k=1}d_k\right\|_{L^p_X(\Omega,\mu)},
$$
for every martingale difference sequence $(d_k)_{k=1}^n \in (L_X^p(\Omega,\mu))^n$ and $(\varepsilon_k)^n_{k=1} \in \{-1,1\}^n$. UMD Banach spaces are very important in harmonic analysis. Burkholder \cite{Bu} and Bourgain \cite{Bou} proved that a Banach space $X$ is UMD if and only if the Hilbert transform $\mathrm{H}$ defined by
$$
\mathrm{H}(f)(x) = \lim_{\varepsilon \rightarrow 0^+}\frac{1}{\pi}\int_{|x-y|>\varepsilon}\frac{f(y)}{x-y} dy,\;\;\mbox{ a.e. } x \in \mathbb R,
$$
for every $f \in L^p(\mathbb R)$, $1<p<\infty$, can be extended from $L^p(\mathbb R)\otimes X$ to $L^p_X(\mathbb R)$ as a bounded operator from $L^p_X(\mathbb R)$ into itself, for every (equivalently, for some) $1<p<\infty$. We refer to \cite[Chapters 4 and 5]{HNVW} for information on UMD-spaces.

Xu (\cite[Theorem 4.1]{Xu1}) established the equivalence of the UMD-property of a Banach lattice with a two-sided $L^p$-estimates for $g$-functions involving Poisson integrals on $\mathbb T$. Hyt\"onen (\cite[Theorem 1.6]{Hy}) characterized Banach UMD spaces by Littlewood-Paley function defined by stochastic integrals. He uses again Poisson integrals on the circle. For Banach lattices, the results obtained by Hyt\"onen reduce to the ones given by Xu.

We now state our results about K\"othe function spaces with the UMD-valued functions and $g$-functions in the inverse Gaussian setting.

Let $X$ be a K\"othe function space and $1<p<\infty$. Note that if $f \in L^p_X(\mathbb R^n,\gamma_{-1})$, then $f$ can be identify with a function $f$ defined in $\mathbb R^n \times \Omega$. We consider, for every $k \in \mathbb N$,
$$
g_{\{P_t^\mathcal A\}_{t>0}}^{k,X}(f)(x,\omega)=\left(\int_0^\infty \left|t^k\partial_t^kP_t^\mathcal A(f(\cdot,\omega))(x)\right|^2\frac{dt}{t}\right)^{1/2},\;\;x\in \mathbb R^n\;\text{and}\; \omega \in \Omega,
$$
for every $f\in L^p_X(\mathbb R^n,\gamma_{-1})$.

Note that now $g$-functions are defined in a different way than above. These definitions can be given for a Banach lattice $X$ where the absolute value is defined but they can not be defined in this way for general Banach spaces $X$. 

Next result can be seen as a version of \cite[Theorem 1.6]{Hy} in the inverse Gauss setting.

\noindent{\bf Theorem D} {\it Let $X$ be a K\"othe function space, $1<p<\infty$, $k\in \mathbb N$ and $\{P^t\}_{t>0}$ a subordinated semigroup of a certain symmetric difussion semigroup.
\begin{enumerate}
\item[(i)] If $X$ has the UMD property, then there exists $C>0$ such that
$$
\frac{1}{C}\|f\|_{L^p_X(\mathbb R^n,\gamma_{-1})}\leq \|g^{k,X}_{\{P^t\}_{t>0}}(f)\|_{L^p_X(\mathbb R^n,\gamma_{-1})}\leq C \|f\|_{L^p_X(\mathbb R^n,\gamma_{-1})},
$$
for every $f \in L^p_X(\mathbb R^n,\gamma_{-1})$.
\item[(ii)] If the inequalities in (i) hold for $k=1$ and $\{P^t\}_{t>0}=\{P_t\}_{t>0}$ is the Euclidean Poisson semigroup, then $X$ has the UMD property.
\end{enumerate}}

\begin{thm}\label{Th1.5}
Let $X$ be an order continuous K\"othe function space and $1<p<\infty$. The following assertions are equivalent.
\begin{enumerate}
\item[(i)] $X$ has the UMD-property.
\item[(ii)] There exists $C>0$ such that
$$
\left\| g_{\{P_t^\mathcal A\}_{t>0}}^{1,X}(f)\right\|_{L^p_X(\mathbb R^n,\gamma_{-1})}\leq C\|f\|_{L^p_X(\mathbb R^n,\gamma_{-1})},\;\; f \in L^p_X(\mathbb R^n,\gamma_{-1}),
$$
and
$$
\left\| g_{\{P_t^\mathcal A\}_{t>0}}^{1,X^*}(f)\right\|_{L^{p'}_{X^*}(\mathbb R^n,\gamma_{-1})} \leq C\|f\|_{L^{p'}_{X^*}(\mathbb R^n,\gamma_{-1})},\;\; f \in L^{p'}_{X^*}(\mathbb R^n,\gamma_{-1}),
$$
where $X^*$ denotes the dual of $X$ and $p'=\frac{p}{p-1}$.
\end{enumerate}
\end{thm}

Motivated by \cite[Theorem 4.1]{Xu1} we stated the following result.
\begin{thm}\label{Th1.6} Let $X$ be an order continuous K\"othe function space, $1<p<\infty$ and $n\in \mathbb{N}$, $n\geq 2$. The following assertions are equivalent.
\begin{enumerate}
\item[(i)] $X$ has the UMD property.
\item[(ii)] There exists $C>0$ such that, for every $i=1,\ldots,n$, 
$$
\left\|g^X_{i,\{P_t^\mathcal A\}_{t>0}}(f)\right\|_{L^p_X(\mathbb R^n,\gamma_{-1})} \leq C\|f\|_{L^p_X(\mathbb R^n,\gamma_{-1})},\;\; f \in L^p_X(\mathbb R^n,\gamma_{-1}),
$$
and
$$
\left\|g^{1,X^*}_{\{P_t^{\mathcal A-I}\}_{t>0}}(f)\right\|_{L^{p'}_{X^*}(\mathbb R^n,\gamma_{-1})} \leq C\|f\|_{L^{p'}_{X^*}(\mathbb R^n,\gamma_{-1})},\;\; f \in L^{p'}_{X^*}(\mathbb R^n,\gamma_{-1}).
$$
\end{enumerate}
Here $X^*$ denotes the dual space of $X$ and $p'=\frac{p}{p-1}$.
\end{thm}

\noindent{\bf Remark}.
When $n=1$ we have that property $(ii)$ implies that the space $X$ has the UMD property (see the proof of Theorem \ref{Th1.6}).

In the following sections we prove our results. As it is well-known Littlewood-Paley functions can be seen as vector valued singular integrals. Then, our objective will be to prove $L^p$-boundedness properties of certain singular integral operators. In order to do this we decompose the operator under consideration into two parts: a local part, that is the singular part, and a global one. This method of decomposition of operators in local and global parts was used by Muckenhoupt (\cite{Mu1}) on the Gaussian context and then it has been used successfully in other contexts (see \cite{BFMT},\cite{BFRST},\cite{BS}, \cite{HTV}, \cite{Pe}, and \cite{PS}, for instance). The global operator is controlled by a positive integral operator that must have the $L^p(\mathbb R^n,\gamma_{-1})$-boundedness properties that we want. On the other hand we must prove that the local operator satisfies the same $L^p(\mathbb R^n,\gamma_{-1})$-boundedness properties than the operator appearing when we replace the operator $\mathcal A$ by the Laplacian operator $-\Delta$.

Throughout this paper, $m(x)=\min\{1,|x|^{-2}\}$, $x\in \mathbb{R}^n\setminus\{0\}$, and $m(0)=1$, and for every $\nu>0$ we define the so called local region $N_\nu$ by
$$
N_\nu = \left\{(x,y)\in \mathbb R^n \times \mathbb R^n: |x-y|< \nu n \sqrt{m(x)}\right\}.
$$
We name global regions as $N_\nu^c,\;\nu >0$. Observe that $\sqrt{m(x)}\sim (1+|x|)^{-1}$, $x\in \mathbb{R}^n$. We always use $c$ and $C$ to denote positive constants that can change in each occurrence.

\section{Proof of Theorem \ref{Th1.1}}

\subsection{Proof of Theorem \ref{Th1.1} (i), for 1$<$p$<\infty$}\label{S2.1} Since $\{T_t^\mathcal A\}_{t>0}$ is a symmetric diffusion semigroup, by \cite[Theorem 1.5]{BFGM} we deduce that the $g$-function $g^q_{\beta,0,\{T_t^\mathcal A\}_{t>0}}$ is bounded from $L^p(\mathbb R^n,\gamma_{-1})$ into itself.

Let $f \in L^p(\mathbb R^n,\gamma_{-1})$. We have that
$$
\partial^k_x T_t^\mathcal A(f)(x) = \int_{\mathbb R^n}\partial^k_xT_t^\mathcal A(x,y) f(y) dy,\;\; x \in \mathbb R^n.
$$
As in \cite[Proposition 3.1]{TZ} we obtain
$$
g^q_{\beta_1,k,\{T_t^\mathcal A\}_{t>0}}(f) \leq C g^q_{\beta_2,k,\{T_t^\mathcal A\}_{t>0}}(f), \;\;0<\beta_1 <\beta_2.
$$
Thus, in order to prove (i) it is sufficient to see that $g^q_{m,k;\{T_t^\mathcal A\}_{t>0}}$ is bounded from $L^p(\mathbb R^n,\gamma_{-1})$ into itself, for every $m\in  \mathbb N$.

Let $m\in \mathbb N$. In \cite[Theorem 5.1]{Teu} Teuwen obtained an explicit expression for $\partial^m_tT_t^\mathcal{L}(x,y)$, $x,y \in \mathbb R^n$ and $t>0$, where
$$
T_t^\mathcal L(x,y)= \frac{{e^{-\frac{|y-e^{-t}x|^2}{1-e^{-2t}}}}}{(1-e^{-2t})^{n/2}},\;\; x,y \in \mathbb R^n\;\;\text{and}\;t>0.
$$
We think that some signs in the expression of \cite[Theorem 5.1]{Teu} are not correct. The new corrected equality is the following
\begin{align}\label{2.1}
\partial_t^mT_t^\mathcal L(x,y)&= (-1)^mT^\mathcal L_t(x,y) \sum_{|r|=m} {m\choose {r_1\ldots r_m}} \prod_{i=1}^n\sum_{s_i=0}^{r_i}\sum_{\ell_i=0}^{s_i}\frac{(-1)^{s_i+\ell_i}}{2^{s_i}}\stirling{r_i}{s_i}{s_i \choose \ell_i}\nonumber \\ 
&\quad \times \left(\frac{e^{-t}}{\sqrt{1-e^{-2t}}}\right)^{2s_i-\ell_i}H_{l_i}(x_i)H_{2s_i-\ell_i}\left(\frac{y_i- e^{-t}x_i}{\sqrt{1-e^{-2t}}}\right), 
\end{align}
for each $x=(x_1,\ldots,x_n)$, $y=(y_1,\ldots,y_n) \in \mathbb R^n$ and $t>0$. Here, for every $N,\ell \in \mathbb N$, $N \geq \ell$, the Stirling number of second kind $\stirling{N}{\ell}$ is defined as the number of partitions of an $N$-set into $\ell$ non-empty subset.

Since $T_t^\mathcal{A}(x,y)=\pi^{-n/2}e^{-nt}T_t^\mathcal{L}(y,x)$, $x,y\in \mathbb{R}^n$, $t>0$, by using (\ref{2.1}) we get
\begin{align*}
\partial_t^mT_t^\mathcal A(x,y)&=(-1)^m T_t^\mathcal{A}(x,y)\sum_{j=0}^m{m \choose j}n^{m-j}\sum_{|r|=j} {j\choose r_1\ldots r_n}\prod_{i=1}^n\sum_{s_i=0}^{r_i}\sum_{\ell_i=0}^{s_i}\frac{(-1)^{s_i+\ell_i}}{2^{s_i}}\stirling{r_i}{s_i}{s_i \choose \ell_i}\nonumber\\
&\quad \times \left(\frac{e^{-t}}{\sqrt{1-e^{-2t}}}\right)^{2s_i-\ell_i}H_{\ell_i}(y_i)H_{2s_i-\ell_i}\left(\frac{x_i-e^{-t}y_i}{\sqrt{1-e^{-2t}}}\right),
\end{align*}
for every $x=(x_1,\ldots,x_n)$, $y=(y_1,\ldots,y_n)\in \mathbb R^n$ and $t>0$.

By taking into account that $\frac{d}{dz}\widetilde H_\ell(z)=-\widetilde H_{\ell+1}(z)$, $z\in \mathbb R$ and $\ell \in \mathbb N$, we have that
\begin{align}\label{cuentaderiv}
\partial_x^{k}\partial_t^mT_t^\mathcal{A}(x,y)&\nonumber\\\
&\hspace{-0.7cm}=\frac{(-1)^me^{-nt}}{\pi ^{\frac{n}{2}}(1-e^{-2t})^{\frac{n}{2}}} \sum_{j=0}^m{m \choose j}n^{m-j}\sum_{|r|=j} {j\choose r_1\ldots r_n}\prod_{i=1}^n\sum_{s_i=0}^{r_i}\sum_{\ell_i=0}^{s_i}\frac{(-1)^{s_i+\ell_i+k_i}}{2^{s_i}}\stirling{r_i}{s_i}{s_i \choose \ell_i}\nonumber\\
&\hspace{-0.7cm}\quad \times \left(\frac{e^{-t}}{\sqrt{1-e^{-2t}}}\right)^{2s_i-\ell_i}\frac{1}{(\sqrt{1-e^{-2t}})^{k_i}}H_{\ell_i}(y_i)\widetilde{H}_{2s_i-\ell_i+k_i}\left(\frac{x_i-e^{-t}y_i}{\sqrt{1-e^{-2t}}}\right)\nonumber\\
&\hspace{-0.7cm}=\frac{(-1)^me^{-nt}}{\pi ^{\frac{n}{2}}(1-e^{-2t})^{\frac{n}{2}}} \sum_{j=0}^m{m \choose j}n^{m-j}\sum_{|r|=j} {j\choose r_1\ldots r_n}\sum_{s\in I_r}\sum_{\ell \in I_s}\frac{(-1)^{|s|+|\ell|+|k|}}{2^{|s|}}\frac{1}{(\sqrt{1-e^{-2t}})^{|k|}}\\
&\hspace{-0.7cm}\quad \times \left(\frac{e^{-t}}{\sqrt{1-e^{-2t}}}\right)^{2|s|-|\ell|}H_{\ell}(y)\widetilde{H}_{2s-\ell+k}\left(\frac{x-e^{-t}y}{\sqrt{1-e^{-2t}}}\right)\prod_{i=1}^n\stirling{r_i}{s_i}{s_i \choose \ell_i},\quad x,y\in \mathbb{R}^n,\;t>0.\nonumber
\end{align}
Here, $I_r$, $r=(r_1,...,r_n)\in \mathbb{N}^n$, represents the set of all $s=(s_1,...,s_n)\in \mathbb{N}^n$ such that $0\leq s_j\leq r_j$, $j=1,...,n$. 

Then, it follows that, for every $0<\eta <1$ there exists $C>0$ such that for every $x,y\in \mathbb{R}^n$ and $t>0$,
\begin{equation}\label{acotderiv}
|t^{m+\frac{|k|}{2}}\partial_x^{k}\partial_t^mT_t^\mathcal{A}(x,y)|\leq C\sum_{j=0}^m\sum_{|r|=j}\sum_{s\in I_r}\sum_{\ell \in I_s}\frac{t^{m+\frac{|k|}{2}}e^{-t}(1+|y|)^{|\ell|}(1-e^{-2t})^{\frac{|\ell|}{2}}}{(1-e^{-2t})^{\frac{n}{2}+\frac{|k|}{2}+|s|}}e^{-\eta\frac{|x-e^{-t}y|^2}{1-e^{-2t}}}.
\end{equation}
We now observe that, when $s\in \mathbb{N}^n$, $|s|\leq m$,
$$
\frac{t^{m+\frac{|k|}{2}}e^{-t}}{(1-e^{-2t})^{\frac{|k|}{2}+|s|}}\leq Ce^{-\frac{t}{2}},\quad t>0.
$$
Also, by performance the change of variables $t=\log \frac{1+s}{1-s}$, $t \in (0,\infty)$, and since $|y|\leq \frac{1}{2}(|y+x|+|y-x|)$, $x,y \in \mathbb R^n$, we obtain for every $\ell \in \mathbb{N}^n$,
\begin{align*}
(1+|y|)^\ell (1-e^{2t})^{\frac{|\ell|}{2}}e^{-\eta\frac{|x-e^{-t}y|^2}{1-e^{-2t}}}&\leq C(1+|y|)^\ell s^{\frac{|\ell|}{2}}e^{\eta\frac{|y|^2-|x|^2}{2}}e^{-\frac{\eta}{4} \big(s|x+y|^2+\frac{|x-y|^2}{s}\big)}\\
&\leq C\frac{(1+|y|)^\ell}{(|x+y|+|x-y|)^\ell}e^{\eta\frac{|y|^2-|x|^2}{2}}e^{-\frac{\delta}{4} \big(s|x+y|^2+\frac{|x-y|^2}{s}\big)}\\
&\leq Ce^{(\eta-\delta)\frac{|y|^2-|x|^2}{2}}e^{-\delta\frac{|x-e^{-t}y|^2}{1-e^{-2t}}},\quad x,y\in\mathbb{R}^n,\;t>0,
\end{align*}
where $0<\delta<\eta$. As we will see, we choose, for convenience, $0<\delta <\eta<1$ such that $\frac{1}{n}>q(1-\delta)$ and $\eta -\delta < \frac{2}{p} < \eta + \delta$.

From \eqref{acotderiv} we deduce that
\begin{equation}\label{AcotDeriv}
|t^{m+\frac{|k|}{2}}\partial_x^{k}\partial_t^mT_t^\mathcal{A}(x,y)|\leq C\frac{e^{-\frac{t}{2}}}{(1-e^{-2t})^{\frac{n}{2}}}e^{(\eta-\delta)\frac{|y|^2-|x|^2}{2}}e^{-\delta\frac{|x-e^{-t}y|^2}{1-e^{-2t}}},\quad x,y\in\mathbb{R}^n,\;t>0.
\end{equation}

We can write
$$
g^q_{m,k,\{T_t^\mathcal A\}_{t>0}}(f) \leq g^q_{m,k,\{T_t^\mathcal A\}_{t>0};\,{\rm loc}}(f) + g^q_{m,k,\{T_t^\mathcal A\}_{t>0};\,{\rm glob}}(f),
$$
where
$$
g^q_{m,k,\{T_t^\mathcal A\}_{t>0};\,{\rm loc}}(f)(x)= \left(\int_0^\infty \left|t^{m+\frac{|k|}{2}}\partial^m_t\partial_x^kT_t^\mathcal A(\mathcal{X}_{N_\nu}(x,\cdot)f)(x)\right|^q\frac{dt}{t}\right)^{\frac{1}{q}},\quad x\in \mathbb R^n,
$$
and
$$
g^q_{m,k,\{T_t^\mathcal A\}_{t>0};\,{\rm glob}}(f)(x)=  \left(\int_0^\infty \left|t^{m+\frac{|k|}{2}}\partial^m_t\partial_x^kT_t^\mathcal A(\mathcal{X}_{N_\nu^c}(x,\cdot)f)(x)\right|^q\frac{dt}{t}\right)^{\frac{1}{q}},\quad x\in \mathbb R^n.
$$

Minkowski inequality leads to
$$
g^q_{m,k,\{T_t^\mathcal A\}_{t>0};\,{\rm glob}}(f)(x)\leq \int_{\mathbb R^n} \mathcal{X}_{N_\nu^c}(x,y)\mathcal R_{m,k}^\mathcal A (x,y)f(y)dy,\quad x \in \mathbb R^n,
$$
where
$$
R_{m,k}^\mathcal A (x,y)=\left\|t^{m+\frac{|k|}{2}}\partial^m_t\partial_x^kT_t^\mathcal A(x,y)\right\|_{L^q((0,\infty),\frac{dt}{t})},\quad x,y\in \mathbb R^n.
$$

Our objective is to see that the operator $R_{m,k}^\mathcal A$ defined by
$$
R_{m,k}^\mathcal A (f)(x)= \int_{\mathbb R^n} R_{m,k}^\mathcal A(x,y) \mathcal{X}_{N_1^c}(x,y) f(y) dy,\;\;x\in \mathbb R^n,
$$
is bounded from $L^p(\mathbb R^n,\gamma_{-1})$ into itself.

By \eqref{AcotDeriv} and since $|x-ry|^2=|y-rx|^2+(1-r^2)(|x|^2-|y|^2)$, $x,y\in \mathbb{R}^n$ and $r\in \mathbb{R}$,
\begin{align*}
[R_{m,k}^\mathcal A(x,y)]^q&\leq Ce^{q(\eta-\delta)\frac{|y|^2-|x|^2}{2}}\int_0^\infty \frac{e^{-q\frac{t}{2}}e^{-\delta q\frac{|x-e^{-t}y|^2}{1-e^{-2t}}}}{(1-e^{-2t})^{\frac{nq}{2}}}\frac{dt}{t}\\
&\leq Ce^{q(\eta+\delta)\frac{|y|^2-|x|^2}{2}}\int_0^\infty \frac{e^{-t}e^{-\delta q\frac{|y-e^{-t}x|^2}{1-e^{-2t}}}}{(1-e^{-2t})^{\frac{nq}{2}}}\frac{dt}{t},\quad x,y\in \mathbb{R}^n.
\end{align*}
Let $(x,y)\in N_1^c$. 
As in \cite[p. 861]{PS} we denote $a=|x|^2+|y|^2$ and $b=2\langle x,y\rangle$. Also we make the change of variable $1-s=e^{-2t}$, $t \in (0,\infty)$, and we define $u(s)=\frac{|y-x\sqrt{1-s}|^2}{s}$, $s \in (0,1)$. We have that
\begin{align*}
[R_{m,k}^\mathcal A(x,y)]^q&\leq Ce^{q(\eta+\delta)\frac{|y|^2-|x|^2}{2}}
\int_0^1 \frac{e^{-q\delta u(s)}}{s^{\frac{nq}{2}}(-\log(1-s))\sqrt{1-s}}ds\\
&\leq Ce^{q(\eta+\delta)\frac{|y|^2-|x|^2}{2}}
\int_0^1 \frac{e^{-q\delta u(s)}}{s^{\frac{nq}{2}+1}\sqrt{1-s}}ds.
\end{align*}
Assume now $b\leq 0$. Then, $u(s) \geq \frac{a}{s}-|x|^2$, $s\in (0,1)$, and
\begin{align*}
\int_0^1 \frac{e^{-q\delta u(s)}}{s^{\frac{nq}{2}+1}\sqrt{1-s}}ds&\leq  e^{-q\delta|y|^2}\int_0^1 \frac{e^{-q\delta (\frac{a}{s}-a})}{s^{\frac{nq}{2}+1}\sqrt{1-s}}ds
\leq C e^{-q\delta|y|^2}\int_0^\infty e^{-q\delta r}\frac{(r+a)^{\frac{nq-1}{2}}}{a^{\frac{nq}{2}}}\frac{dr}{\sqrt r}\\
&\leq Ce^{-q\delta|y|^2}\int_0^\infty e^{-q\delta r}(r+1)^{\frac{nq-1}{2}}\frac{dr}{\sqrt r}\leq Ce^{-q\delta|y|^2} ,
\end{align*}
because $a \geq \frac{n}{2}$. We obtain
$$
R_{m,k}^\mathcal A(x,y)\leq C e^{-(\eta+\delta)\frac{|x|^2}{2}+(\eta-\delta)\frac{|y|^2}{2}}.
$$
Suppose that $b>0$. Since $\frac{1}{n}>q(1-\delta)$, we can write (see \cite[p. 862]{PS})
\begin{align*}
\int_0^1 \frac{e^{-q\delta u(s)}}{s^{\frac{nq}{2}+1}\sqrt{1-s}}ds&= \int_0^1 \left(\frac{e^{-u(s)}}{s^{\frac{n}{2}}}\right)^{\frac{nq-1}{n}} \frac{e^{-(q(\delta-1 )+\frac{1}{n})u(s)}}{s^{\frac{3}{2}}\sqrt{1-s}}ds\\
&\leq C\left(\frac{e^{-u_0}}{s_0^{\frac{n}{2}}}\right)^{\frac{nq-1}{n}}\int_0^1\frac{e^{-(q(\delta -1) + \frac{1}{n})u(s)}}{s^{\frac{3}{2}}\sqrt{1-s}}ds,
\end{align*}
where $s_0=2\frac{\sqrt{a^2-b^2}}{a+\sqrt{a^2-b^2}}$ and $u_0=u(s_0)= \frac{|y|^2-|x|^2}{2}+ \frac{|x+y||x-y|}{2}$. 

We now proceed as in the proof of \cite[Lemma 2.3]{PS}. Since $a^2-b^2 \geq 1$, we get
\begin{align*}
\int_0^1 \frac{e^{-(q(\delta -1) + \frac{1}{n})u(s)}}{s^{\frac{3}{2}}\sqrt{1-s}}ds &\leq C\frac{e^{-(q(\delta-1)+\frac{1}{n})u_0}}{\sqrt{s_0}}\int_0^\infty \frac{e^{-(q(\delta -1) + \frac{1}{n})s}}{(a^2-b^2)^{\frac{1}{4}}}\left(\sqrt s + \frac{1}{\sqrt{s}}\right)ds\\
&\leq C\frac{e^{-(q(\delta-1)+\frac{1}{n})u_0}}{\sqrt{s_0}}.
\end{align*}
Then,
$$
\int_0^1 \frac{e^{-q\delta u(s)}}{s^{\frac{nq}{2}+1}\sqrt{1-s}}ds\leq C\frac{e^{-q\delta u_0}}{s_0^{\frac{nq}{2}}}.
$$
According to \cite[p. 863]{PS}, since $|x-y||x+y|\geq n$ and $s_0\sim \frac{\sqrt{a^2-b^2}}{a}$, we obtain
\begin{align*}
R_{m,k}^\mathcal{A} (x,y)&\leq  C \frac{e^{(\eta+ \delta)\frac{|y|^2-|x|^2}{2}-\delta u_0}}{s_0^{n/2}}=C\frac{e^{\frac{\eta}{2}(|y|^2-|x|^2)-\frac{\delta}{2}|x+y||x-y|}}{s_0^{n/2}}\\
&\leq  C\Big(\frac{|x+y|}{|x-y|}\Big)^{n/2}e^{\frac{\eta}{2}(|y|^2-|x|^2)-\frac{\delta}{2}|x+y||x-y|}\leq  C |x+y|^n e^{\frac{\eta}{2}(|y|^2-|x|^2)-\frac{\delta}{2}|x+y||x-y|}.
\end{align*}

By putting together the above estimates we get, when $(x,y) \in N_1^c$,
\begin{equation}\label{A2}
R^\mathcal A_{m,k}(x,y)\leq C \left\{\begin{array}{ll}
e^{-(\eta+\delta)\frac{|x|^2}{2} +(\eta-\delta)\frac{|y|^2}{2}},&\langle x,y \rangle \leq  0, \\[0.2cm]
|x+y|^ne^{\frac{\eta}{2}(|y|^2-|x|^2)-\frac{\delta}{2}|x+y||x-y|},&\langle x,y \rangle > 0. \end{array}\right.
\end{equation}
Since $\left||y|^2-|x|^2\right| \leq |x+y||x-y|$, $x,y\in \mathbb{R}^n$, and that we have chosen $\eta -\delta < \frac{2}{p} < \eta + \delta$, we obtain (see \cite[p. 501]{Pe})
\begin{align*}
\int_{\mathbb R^n} e^{\frac{|x|^2}{p} - \frac{|y|^2}{p}}R^\mathcal A_{m,k}(x,y) \mathcal{X}_{N_1^c}(x,y)dy &\\
&\hspace{-5.5cm}\leq C\left( \int_{\mathbb R^n}
e^{\left(\frac{1}{p}-\frac{\eta +\delta }{2}\right)|x|^2}e^{-|y|^2\left(\frac{1}{p}+\frac{\eta-\delta}{2}\right)}dy+\int_{\mathbb R^n} |x+y|^n e^{-\left(\frac{\delta}{2}-\left|\frac{\eta}{2}-\frac{1}{p}\right|\right)|x-y||x+y|}dy\right)\leq C,\;\;x \in \mathbb R^n.
\end{align*}
In a similar way we can see that
$$
\sup_{y \in \mathbb R^n}\int_{\mathbb R^n}e^{\frac{|x|^2-|y|^2}{p}}R^\mathcal A_{m,k}(x,y) \mathcal{X}_{N_1^c}(x,y)dx < \infty.
$$
It follows that the operator $R_{m,k;p}$ defined by
$$
R^\mathcal A_{m,k;p}(f)(x)= \int_{\mathbb R^n}e^{\frac{|x|^2-|y|^2}{p}}R^\mathcal A_{m,k}(x,y)\mathcal{X}_{N_1^c}(x,y)f(y) dy,\;\;x\in \mathbb R^n,
$$
is bounded from $L^r(\mathbb R^n,dx)$ into itself, for every $1 \leq r < \infty$. Then, the operator $R^\mathcal A_{m,k}$ is bounded from $L^p(\mathbb R^n,\gamma_{-1})$ into itself. Hence, $g_{m,k,\{T_t^\mathcal A\}_{t>0};\,{\rm glob}}^q$ is bounded from $L^p(\mathbb R^n,\gamma_{-1})$ into itself. 

We now deal with $g_{m,k,\{T_t\}_{t>0};\,{\rm loc}}^q$. We define the operator $G_{m,k,\{T^\mathcal A_t\}_{t>0}}$ by
$$
G_{m,k,\{T^\mathcal A_t\}_{t>0}}(f)(x)=t^{m+\frac{|k|}{2}} \partial_t^m\partial_x^kT_t^\mathcal A(f)(x),\;x\in \mathbb R^n\;\;\text{and}\; t>0.
$$
We are going to see that
\begin{enumerate}
\item[(a)]$G_{m,k,\{T_t\}_{t>0}}$ defines a bounded operator from $L^q(\mathbb R^n,\gamma_{-1})$ into $L^q_{L^q((0,\infty)\frac{dt}{t})}(\mathbb R^n,\gamma_{-1})$.
\item[(b)]$R_{m,k}^\mathcal{A}(x,y) \leq \frac{C}{|x-y|^n}$, $(x,y) \in N_2$.
\item[(c)] For every $i=1,\ldots,n$,
$$
\left(\int_0^\infty \left|t^{m+\frac{|k|}{2}}\partial_t^m\partial_{x_i}\partial_x ^kT_t^\mathcal A(x,y)\right|^q\frac{dt}{t}\right)^{\frac{1}{q}} \leq \frac{C}{|x-y|^{n+1}},\;\; (x,y) \in N_2.
$$
\end{enumerate}
This means that $G_{m,k,\{T^\mathcal A_t\}_{t>0}}$ is a local Calder\'on-Zygmund operator with respect to the Banach space $L^q((0,\infty),\frac{dt}{t})$ in the inverse Gaussian setting. Salogni \cite[\S 3.2]{Sa} considered scalar local Calder\'on-Zygmund operators in the inverse Gaussian context. The properties continue being valid in a Banach valued setting.

First we show (b). Assume that $(x,y)\in N_2$. Then, $||x|^2-|y|^2|\leq C$ and according to \eqref{AcotDeriv} we obtain that
$$
R_{m,k}^\mathcal{A}(x,y)^q\leq C\int_0^\infty\frac{e^{-\delta q\frac{|x-e^{-t}y|^2}{1-e^{-2t}}}}{(1-e^{-2t})^{\frac{nq}{2}}}\frac{e^{-\frac{qt}{2}}}{t}dt.
$$
The argument in the proof of \cite[Lemma 3.3.1]{Sa} leads to
$$
\int_0^{m(x)}\frac{e^{-\delta q\frac{|x-e^{-t}y|^2}{1-e^{-2t}}}}{(1-e^{-2t})^{\frac{nq}{2}}}\frac{e^{-\frac{qt}{2}}}{t}dt\leq C\int_0^{m(x)}\frac{e^{-c\frac{|x-y|^2}{t}}}{t^{\frac{nq}{2}+1}}dt.
$$
Also, we have 
$$
\int_{m(x)}^\infty\frac{e^{-\delta q\frac{|x-e^{-t}y|^2}{1-e^{-2t}}}}{(1-e^{-2t})^{\frac{nq}{2}}}\frac{e^{-\frac{qt}{2}}}{t}dt\leq C\int_{m(x)}^\infty\frac{e^{-\frac{qt}{2}}}{(1-e^{-2t})^{\frac{nq}{2}}}\frac{dt}{t}\leq C\int_{m(x)}^\infty \frac{dt}{t^{\frac{nq}{2}+1}}\leq \frac{C}{(m(x))^{\frac{nq}{2}}}\leq \frac{C}{|x-y|^{nq}}.
$$
We conclude that
$$
R_{m,k}^\mathcal{A}(x,y)\leq \frac{C}{|x-y|^n},\;(x,y) \in N_2.
$$
Thus (b) is proved. The estimation in (c) can be seen in a similar way.

Next we handle property (a). Let $f\in L^q(\mathbb R^n,\gamma_{-1})$ and write  $\mathbb{B}_q=L^q((0,\infty),\frac{dt}{t})$. Since $q\geq 2$, by using H\"older's inequality we get
\begin{align*}
\left\| G_{m,k,\{T_t^{\mathcal A}\}_{t>0}}^q(f) \right\|^q_{L^q_{\mathbb{B}_q}(\mathbb{R}^n,\gamma_{-1})} &= \int_{\mathbb R^n}\int_0^\infty\Big|t^{m+\frac{|k|}{2}}\partial_t^m\partial_x^kT_t^{\mathcal A}(f)(x)\Big|^q \frac{dt}{t}d\gamma_{-1}(x) \\
&\hspace{-2cm}=\int_{\mathbb R^n}\int_0^\infty \big|t^{m+\frac{|k|}{2}}\partial_t^m\partial_x^k T_t^{\mathcal A}(f)(x)\big|^2 \sup_{t>0}\big|t^{m+\frac{|k|}{2}}\partial_t^m\partial_x^kT_t^{\mathcal A}(f)(x)\big|^{q-2} \frac{dt}{t}d\gamma_{-1}(x)\\
&\hspace{-2cm} \leq \left(\int_{\mathbb R^n} \big|g^2_{m,k,\{T_t^{\mathcal A}\}_{t>0}}(f)(x)\big|^q d\gamma_{-1}(x)\right)^{\frac{2}{q}}\left(\int_{\mathbb R^n}\big|T^{\mathcal A}_{m,k,*}(f)(x)\big|^q d\gamma_{-1}(x)\right)^{1-\frac{2}{q}}.
\end{align*}
Here $T^\mathcal A_{m,k,*}$ is the maximal operator defined by
$$
T^\mathcal A_{m,k,*}f(x) =  \sup_{t>0} \Big|t^{m+\frac{|k|}{2}} \partial_t^m\partial_x^kT_t^{\mathcal A}(f)(x)\Big|,\quad x\in \mathbb{R}^n.
$$
Thus, to get (a) it is sufficient to see that $g^2_{m,k,\{T_t^\mathcal A\}_{t>0}}$ and $T^\mathcal A_{m,k,*}$ are bounded from $L^q(\mathbb R^n,\gamma_{-1})$ into itself, because these properties allow us to conclude that
$$
\left\| G_{m,k,\{T_t^\mathcal A\}_{t>0}}(f)\right\|^q_{L^q_{\mathbb{B}_q}(\mathbb R^n,\gamma_{-1}))}\leq  C\|f\|^2_{L^q(\mathbb R^n,\gamma_{-1})}\|f\|^{q-2}_{L^q(\mathbb R^n,\gamma_{-1})}
=C\|f\|^q_{L^q(\mathbb R^n,\gamma_{-1})}.
$$

In order to see that $g^2_{m,k,\{T_t^\mathcal A\}_{t>0}}$ is bounded from $L^q(\mathbb R^n,\gamma_{-1})$ into itself  we are going to establish the following properties:
\begin{enumerate}
\item[(a')] $G_{m,k,\{T_t^\mathcal A\}_{t>0}}$ is bounded from $L^2(\mathbb R^n,\gamma_{-1})$ into $L^2_{\mathbb{B}_2}(\mathbb R^n,\gamma_{-1})$.
\item[(b')] $\big\|t^{m+\frac{k}{2}}\partial_t^m\partial_x^kT_t^\mathcal A(x,y)\big\|_{\mathbb{B}_2} \leq \frac{C}{|x-y|^n}$, $(x,y) \in N_2$,
\item[(c')] For every $i=1\ldots,n$,
$$
\Big\|t^{m+\frac{k}{2}}\partial^m_t \partial_{x_i}\partial_x^kT_t^\mathcal A(x,y)\Big\|_{\mathbb{B}_2}\leq \frac{C}{|x-y|^{n+1}},\quad (x,y) \in N_2.
$$
\end{enumerate}

Note that (b') and (c') coincide with (b) and (c), respectively, when $q=2$. So, we only need to prove (a'). Let $f \in L^2(\mathbb R^n,\gamma_{-1})$. We have that
$$
T^\mathcal A_t(f)(x)=\sum_{\ell\in \mathbb N^n}c_\ell (f) e^{-(n+|\ell|)t}\widetilde{H}_\ell(x),\quad x\in \mathbb R^n.
$$

According to \cite[p. 324]{San}, for every $j \in \mathbb N$,
$$
|\widetilde H_j(u)| \leq 2\sqrt{2^jj!}e^{-\frac{u^2}{2}},\quad u\in \mathbb R.
$$
Also,
$$
\|\widetilde H_j\|_{L^2(\mathbb R^n,\gamma_{-1})}= \sqrt{\pi2^jj!},\quad j \in \mathbb N.
$$
Then, for every $\ell=(\ell_1,\ldots,\ell_n)\in \mathbb N^n$,
\begin{align*}
|c_\ell(f)|& \leq \frac{C}{\|\widetilde  H_\ell\|_{L^2(\mathbb R^n,\gamma_{-1})}}\int_{\mathbb R^n}|f(x)||\widetilde H_\ell(x)|e^{|x|^2}dx\leq C\frac{\|f\|_{L^2(\mathbb R^n,\gamma_{-1})}}{\|\widetilde  H_l\|_{L^2(\mathbb R^n,\gamma{-1})}}.
\end{align*}
We get, for every $\ell=(\ell_1,\ldots,\ell_n)\in \mathbb N^n$,
$$
|c_\ell(f)||\widetilde H_\ell(x)| \leq Ce^{-\frac{|x|^2}{2}},\;x\in \mathbb R^n.
$$
Hence, for every $t>0$, the series defining $T_t^\mathcal A(f)$ is absolute and uniformly convergent in $\mathbb R^n$. Also, since $\frac{d}{du}\widetilde H_j(u)=-\widetilde H_{j+1}(u)$, $u\in \mathbb R$ and $j \in \mathbb N$, we can derivate under the sum sign to obtain
$$
\partial^m_t\partial_x^kT_t^\mathcal A(f)(x) = \sum_{\ell \in \mathbb N^n}(-1)^{|k|+n+|\ell|}c_\ell(f)e^{-t(n+|\ell|)}(n+|\ell|)^m\widetilde H_{\ell+k}(x),\quad x \in \mathbb R^n,\;t>0.
$$
By using Plancherel identity for Hermite polynomials it follows that
\begin{align*}
\left\|t^{m+\frac{|k|}{2}}\partial^m_t \partial_x^kT_t^\mathcal A (f)\right\|^2_{L^2_{\mathbb{B}_2}(\mathbb R^n,\gamma_{-1})}&=\int^\infty_0 t^{2m+|k|}\int_{\mathbb R^n}\big| \partial^m_t \partial_x^k T_t^\mathcal A (f)(x)\big|^2d\gamma_{-1}(x)\frac{dt}{t}\\
&\hspace{-3.5cm}= \int^\infty_0 t^{2m+|k|} \sum_{\ell\in \mathbb N^n}|c_\ell(f)|^2 e^{-2t(n+|\ell|)}(n+|\ell|)^{2m}\|\widetilde H_{\ell+k}\|^2_{L^2(\mathbb R^n,\gamma_{-1})}\frac{dt}{t}\\
&\hspace{-3.5cm}=\sum_{\ell\in \mathbb N^n}|c_\ell(f)|^2\|\widetilde H_\ell\|_{L^2(\mathbb R^n,\gamma_{-1})}^2\frac{\|\widetilde H_{\ell+k}\|_{L^2(\mathbb R^n,\gamma_{-1})}^2}{\|\widetilde H_\ell\|_{L^2(\mathbb R^n,\gamma_{-1})}^2}(n+|\ell|)^{2m}\int_0^\infty e^{-2t(n+|\ell|)}t^{2m+|k|}\frac{dt}{t} \\
&\hspace{-3.5cm}\leq C\sum_{\ell \in \mathbb N^n}|c_\ell(f)|^2\|\widetilde H_\ell\|_{L^2(\mathbb R^n,\gamma_{-1})}^2 \frac{\|\widetilde H_{\ell+k}\|_{L^2(\mathbb R^n,\gamma_{-1})}^2}{\|\widetilde H_{\ell}\|_{L^2(\mathbb R^n,\gamma_{-1})}^2}\frac{1}{(n+|\ell|)^{|k|}} \\
&\hspace{-3.5cm}\leq C\sum_{\ell \in \mathbb N^n}\frac{|c_\ell(f)|^2\|\widetilde H_\ell\|_{L^2(\mathbb R^n,\gamma_{-1})}^2}{(n+|\ell|)^{|k|}}\prod_{j=1}^n\frac{(\ell_j+k_j)!}{\ell_j!}\leq C\sum_{\ell \in \mathbb N^n}\frac{|c_\ell(f)|^2\|\widetilde H_\ell\|_{L^2(\mathbb R^n,\gamma_{-1})}^2}{(n+|l|)^{|k|}}\prod_{j=1}^n\ell_j^{k_j}\\
& \hspace{-3.5cm}\leq C\sum_{\ell \in \mathbb N^n}|c_\ell(f)|^2\|\widetilde H_\ell\|_{L^2(\mathbb R^n,\gamma_{-1})}^2=C\|f\|^2_{L^2(\mathbb R^n,\gamma_{-1})}.
\end{align*}
Thus, we obtain that $G_{m,k,\{T_t^\mathcal{A}\}_{t>0}}$ is bounded from $L^2(\mathbb{R}^n,\gamma_{-1})$ into $L_{\mathbb{B}_2}^2(\mathbb{R}^n,\gamma_{-1})$. According to a vector valued $L^2((0,\infty),\frac{dt}{t})$) version of \cite[Proposition 3.2.5]{Sa} we conclude that the operator $G_{m,k,\{T_t^\mathcal{A}\}_{t>0};\,{\rm loc}}$ defined by
$$
G_{m,k,\{T_t^\mathcal{A}\}_{t>0};\,{\rm loc}}(f)(x)=t^{m+\frac{|k|}{2}}\partial_t^m\partial_x^kT_t^\mathcal{A}(\mathcal{X}_{N_2}(x,\cdot)f)(x),\quad x\in \mathbb{R}^n,\;t>0,
$$
is bounded from $L^r(\mathbb{R}^n,\gamma_{-1})$ into $L^r_{\mathbb{B}_2}(\mathbb{R}^n,\gamma_{-1})$, for every $1<r<\infty$. In other words, we have seen that the operator $g_{m,k,\{T_t^\mathcal{A}\}_{t>0};\,{\rm loc}}^2$ is bounded from $L^r(\mathbb{R}^n,\gamma_{-1})$ into itself, for every $1<r<\infty$. As we had proved, $g_{m,k,\{T_t^\mathcal{A}\}_{t>0};\,{\rm glob}}^2$ is also bounded from $L^r(\mathbb{R}^n,\gamma_{-1})$ into itself, for every $1<r<\infty$. Hence, $g_{m,k,\{T_t^\mathcal{A}\}_{t>0}}^2$ is bounded from $L^q(\mathbb{R}^n,\gamma_{-1})$ into itself.

We now consider the maximal operator
$$
T_{m,k,*}^\mathcal{A}(f)(x)=\sup_{t>0}\Big|t^{m+\frac{|k|}{2}}\partial_t^m\partial_x^kT_t^\mathcal{A}(f)(x)\Big|,\quad x\in \mathbb{R}^n.
$$
We are going to see that 
$$
\mathbb{T}_{m,k,*}^\mathcal{A}(f)(x):=\sup_{t>0}\Big|\mathcal{M}(f)(t,x)|,\quad x\in \mathbb{R}^n,
$$
is bounded from $L^q(\mathbb{R}^n,\gamma_{-1})$ into itself, where
$$
\mathcal{M}(f)(t,x)=\int_{\mathbb{R}^n}\mathbb{M}(t,x,y)f(y)dy, \quad x\in \mathbb{R}^n,\;t>0,
$$
and
$$
\mathbb{M}(t,x,y)=\frac{e^{-\frac{t}{2}}}{(1-e^{-2t})^{\frac{n}{2}}}e^{(\eta +\delta)\frac{|y|^2-|x|^2}{2}}e^{-\delta \frac{|y-e^{-t}x|^2}{1-e^{-2t}}},\quad x,y\in \mathbb{R}^n,\;t>0.
$$
Here, $0<\delta<\eta <1$. We observe that, once we prove this fact, according to \eqref{AcotDeriv} and since $|x-ry|^2=|y-rx|^2+(1-r^2)(|x|^2+|y|^2)$, $x,y\in \mathbb{R}^n$, $r>0$, we can deduce that $T_{m,k,*}^\mathcal{A}$ is bounded from $L^q(\mathbb{R}^n,\gamma_{-1})$ into itself. 

Consider $\mathcal{M}_{\rm loc}$ and $\mathcal{M}_{\rm glob}$ defined by 
$$
\mathcal{M}_{\rm loc}(f)(t,x)=\mathcal{M}(\mathcal{X}_{N_\nu}(x,\cdot)f)(t,x)\quad \mbox{ and }\quad \mathcal{M}_{\rm glob}(f)(t,x)=\mathcal{M}(\mathcal{X}_{N_\nu^c}(x,\cdot)f)(t,x),\quad x\in \mathbb{R}^n,\;t>0.
$$
The operators $\mathbb{T}_{m,k,*;\,{\rm loc}}$ and $\mathbb{T}_{m,k,*;\,{\rm glob}}$ are defined in the obvious way. Here $\nu$ will be fixed later.

We first study $\mathbb{T}_{m,k,*;\;{\rm glob}}^\mathcal{A}$. Let $\lambda, \alpha\in (0,\infty )$. Suppose that $(u,v)\in N_\lambda$, that is, $|u-v|< \lambda n\sqrt{m(u)}$. Then,
$$
|\alpha u-\alpha v|<\alpha \lambda n\sqrt{m(u)}\leq \alpha \lambda n\left\{
\begin{array}{ll}
\sqrt{m(\alpha u)},&0<\alpha <1,\\[0.2cm]
\alpha \sqrt{m(\alpha u)},&\alpha \geq 1.
\end{array}
\right.
$$
Hence, $(\alpha u,\alpha v)\in N_{\alpha \lambda}$, when $0<\alpha <1$, and $(\alpha u,\alpha v)\in N_{\alpha ^2\lambda}$, if $\alpha \geq 1$. Since $0<\delta <1$ we deduce that $(\sqrt{\delta}x,\sqrt{\delta}y)\not \in N_1$ provided that $(x,y)\not \in N_{\nu}$, with $\nu =\delta^{-1}$.

Let $\nu =\delta ^{-1}$. By \cite[Proposition 2.1]{MPS} we can write, for every $(x,y)\in N_\nu ^c$, 
\begin{align*}
\sup_{t>0}\mathbb{M}(t,x,y)&\leq Ce^{(\eta +\delta)\frac{|y|^2-|x|^2}{2}}\left\{\begin{array}{ll}
e^{-\delta|y|^2},&\langle x,y \rangle \leq  0, \\[0.2cm]
|x+y|^ne^{-\frac{\delta}{2}(|y|^2-|x|^2+|x+y||x-y|)},&\langle x,y \rangle > 0,
\end{array}\right.\\
&\leq C \left\{\begin{array}{ll}
e^{-(\eta+\delta)\frac{|x|^2}{2} +(\eta-\delta)\frac{|y|^2}{2}},&\langle x,y \rangle \leq  0, \\[0.2cm]
|x+y|^ne^{\frac{\eta}{2}(|y|^2-|x|^2)-\frac{\delta}{2}|x+y||x-y|},&\langle x,y \rangle > 0.
\end{array}
\right.
\end{align*}

By using that $||y|^2-|x|^2|\leq |x+y||x-y|$, we obtain
\begin{align*}
e^{\frac{|x|^2}{p}-\frac{|y|^2}{p}}\sup_{t>0}\mathbb{M}(t,x,y)&\leq C \left\{\begin{array}{ll}
e^{|x|^2(\frac{1}{p}-\frac{\eta+\delta}{2})}e^{|y|^2(\frac{\eta-\delta}{2}-\frac{1}{p})},&\langle x,y \rangle \leq  0, \\[0.2cm]
|x+y|^ne^{-(\frac{\delta}{2}-|\frac{\eta}{2}-\frac{1}{p}|)|x+y||x-y|},&\langle x,y \rangle > 0,
\end{array}
\right.\quad (x,y)\in N_\nu^c.
\end{align*}

We choose $0<\delta<\eta<1$ such that $\eta-\delta<\frac{2}{p}<\eta +\delta$. It follows that (see \cite[p. 501]{Pe})
$$
\sup_{x\in \mathbb{R}^n}\int_{\mathbb{R}^n} e^{\frac{|x|^2}{p}-\frac{|y|^2}{p}}\sup_{t>0}\mathbb{M}(t,x,y)\mathcal{X} _{N_\nu ^c}(x,y)dy<\infty, 
$$
and, in a similar way,
$$
\sup_{y\in \mathbb{R}^n}\int_{\mathbb{R}^n} e^{\frac{|x|^2}{p}-\frac{|y|^2}{p}}\sup_{t>0}\mathbb{M}(t,x,y)\mathcal{X} _{N_\nu ^c}(x,y)dx<\infty. 
$$
Then, the operator  $\mathbb{T}_{m,k,*;\;{\rm glob}}^\mathcal{A}$ is bounded from $L^r(\mathbb{R}^n,\gamma_{-1})$ into itself, for every $1<r<\infty$. In particular,  $\mathbb{T}_{m,k,*;\;{\rm glob}}^\mathcal{A}$ is bounded from $L^q(\mathbb{R}^n,\gamma_{-1})$ into itself.

On the other hand since $||y|^2-|x|^2|\leq |x-y||x+y|\leq C$ when $(x,y)\in N_\nu$, by proceeding as in the proof of \cite[Lemma 3.3.1]{Sa}
\begin{align*}
\mathbb{M}(t,x,y)&\leq C\frac{e^{-\frac{t}{2}}}{(1-e^{-2t})^{\frac{n}{2}}}e^{(\eta -\delta)\frac{|y|^2-|x|^2}{2}}e^{-\delta \frac{|x-e^{-t}y|^2}{1-e^{-2t}}}\leq C\frac{e^{-\frac{t}{2}}e^{-c\frac{|x-y|^2}{1-e^{-2t}}}}{(1-e^{-2t})^{\frac{n}{2}}}\\
&\leq C\frac{e^{-c\frac{|y-x|^2}{t}}}{t^{\frac{n}{2}}},\quad (x,y)\in N_\nu , \;t>0.
\end{align*}
Then
$$
\mathbb{T}_{m,k,*;\,{\rm loc}}^\mathcal{A}(f)(x)\leq \sup_{t>0}\int_{\mathbb{R}^n}\frac{e^{-c\frac{|y-x|^2}{t}}}{t^{\frac{n}{2}}}|f(y)|dy,\quad x\in \mathbb{R}^n.
$$
Hence, $\mathbb{T}_{m,k,*;\,{\rm loc}}^\mathcal{A}$ is bounded from $L^r(\mathbb{R}^n,dx)$ into itself, for every $1<r<\infty$ and by \cite[Proposition 3.2.5]{Sa} we deduce that $\mathbb{T}_{m,k,*;\,{\rm loc}}^\mathcal{A}$ is bounded from $L^q(\mathbb{R}^n,\gamma_{-1})$ into itself. It follows that $T_{m,k,*;\,{\rm loc}}^\mathcal{A}$ is bounded from $L^q(\mathbb{R}^n,\gamma_{-1})$ into itself and property (a) is established.

Since (a), (b) and (c) holds, by \cite[Theorem 3.2.8]{Sa}, the operator $G_{m,k,\{T_t^\mathcal{A}\}_{t>0};\,{\rm loc}}$ given by
$$
G_{m,k,\{T_t^\mathcal{A}\}_{t>0};\,{\rm loc}}(f)(x)=G_{m,k,\{T_t^\mathcal{A}\}_{t>0}}(\mathcal{X}_{N_2}(x,\cdot )f)(x),\quad x\in \mathbb{R}^n,
$$
is bounded from $L^p(\mathbb{R}^n,\gamma_{-1})$ into $L^p_{\mathbb{B}_q}(\mathbb{R}^n,\gamma_{-1})$. In other words, we have shown that $g_{m,k,\{T_t^\mathcal{A}\}_{t>0};\,{\rm loc}}^q$ is bounded from $L^p(\mathbb{R}^n,\gamma_{-1})$ into itself.

Thus, we conclude that $g_{m,k,\{T_t^\mathcal{A}\}_{t>0}}$ is bounded from $L^p(\mathbb{R}^n,\gamma_{-1})$ into itself.

\subsection{Proof of Theorem \ref{Th1.1} (i), for p$=$1}
In order to see that $g_{\beta ,0,\{T_t^\mathcal{A}\}_{t>0}}^q$, with $0<\beta \leq 1$ is bounded from $L^1(\mathbb{R}^n,\gamma_{-1})$ into $L^{1,\infty }(\mathbb{R}^n,\gamma _{-1})$, it is sufficient to prove that $g_{1,0,\{T_t^\mathcal{A}\}_{t>0}}^q$ is bounded from $L^1(\mathbb{R}^n,\gamma_{-1})$ into $L^{1,\infty}(\mathbb{R}^n,\gamma _{-1})$.

We consider the operator
$$
G_{1,0,\{T_t^\mathcal{A}\}_{t>0}}(f)(t,x)=t\partial _tT_t^\mathcal{A}(f)(x),\quad x\in \mathbb{R}^n,\;t>0.
$$
We proved that this operator is a local Calder\'on-Zygmund one in the $L^q((0,\infty ),\frac{dt}{t})$-setting. Then, $G_{1,0,\{T_t^\mathcal{A}\}_{t>0};\,{\rm loc}}$ is bounded from $L^1(\mathbb{R}^n,\gamma_{-1})$ into $L^{1,\infty} _{\mathbb{B}_q}(\mathbb{R}^n,\gamma _{-1})$, where $\mathbb{B}_q=L^q((0,\infty ),\frac{dt}{t})$. In other words, $g_{1,0,\{T_t^\mathcal{A}\}_{t>0};\,{\rm loc}}^q$ is bounded from $L^1(\mathbb{R}^n,\gamma_{-1})$ into $L^{1,\infty}(\mathbb{R}^n,\gamma _{-1})$.

In order to finish the proof of our objective we have to see that $g_{1,0,\{T_t^\mathcal{A}\}_{t>0};\,{\rm glob}}^q$ is bounded from $L^1(\mathbb{R}^n,\gamma_{-1})$ into $L^{1,\infty }(\mathbb{R}^n,\gamma _{-1})$. To establish this result we use some ideas developed in \cite{BrSj}.

We have that, for every $x,y\in  \mathbb{R}^n$ and $t>0$,
\begin{align*}
t\partial _t T_t^\mathcal{A}(x,y)&=-\frac{te^{|y|^2-|x|^2}T_t^\mathcal{A}(y,x)}{1-e^{-2t}}\left(n+2e^{-t}\sum_{i=1}^n(y_i-x_ie^{-t})x_i-2\frac{e^{-2t}|y-e^{-t}x|^2}{1-e^{-2t}}\right).
\end{align*}
Then,
\begin{align}\label{AcotDerivt}
|t\partial _t T_t^\mathcal{A}(x,y)|&\leq Cte^{-nt}\frac{e^{|y|^2-|x|^2}e^{-\frac{|y-e^{-t}x|^2}{1-e^{-2t}}}}{(1-e^{-2t})^{\frac{n}{2}+1}}\Big(1+e^{-t}|x||y-e^{-t}x|+e^{-2t}\frac{|y-e^{-t}x|^2}{1-e^{-2t}}\Big)\nonumber\\
&\leq Cte^{-nt}\frac{e^{|y|^2-|x|^2}e^{-c\frac{|y-e^{-t}x|^2}{1-e^{-2t}}}}{(1-e^{-2t})^{\frac{n}{2}+1}}\Big(1+e^{-t}|x|(1-e^{-2t})^{\frac{1}{2}}\Big),\quad x,y\in \mathbb{R}^n,\;t>0.
\end{align}

Let $f\in L^1(\mathbb{R}^n,\gamma_{-1})$. Minkowski inequality leads to
$$
g_{1,0,\{T_t^\mathcal{A}\}_{t>0};\,{\rm glob}}^q(f)(x)\leq C\int_{\mathbb{R}^n}f(y)\mathcal{X}_{N_1^c}(x,y)\big\|t\partial_tT_t^\mathcal{A}(x,y)\big\|_{L^q((0,\infty),\frac{dt}{t})}dy,\quad x\in \mathbb{R}^n.
$$

Next we analyze the kernel $\big\|t\partial_tT_t^\mathcal{A}(x,y)\big\|_{L^q((0,\infty),\frac{dt}{t})}$, when $(x,y)\in N_1^c$.
From \eqref{AcotDerivt} and by perfoming the change of variables $r= e^{-t}$, $t\in (0,\infty)$ we get
\begin{align*}
    \big\|t\partial_tT_t^\mathcal{A}(x,y)\big\|_{L^q((0,\infty),\frac{dt}{t})}&\leq Ce^{|y|^2-|x|^2}\left(\int_0^\infty \frac{t^{q-1}e^{-nqt}e^{-c\frac{|y-e^{-t}x|^2}{1-e^{-2t}}}}{(1-e^{-2t})^{(\frac{n}{2}+1)q}}\Big(1+|x|^q(1-e^{-2t})^{\frac{q}{2}}\Big)dt\right)^{\frac{1}{q}}\\
    &\hspace{-2cm}\leq Ce^{|y|^2-|x|^2}\left(\int_0^1 \frac{(-\log r)^{q-1}r^{nq-1}e^{-c\frac{|y-rx|^2}{1-r^2}}}{(1-r^2)^{(\frac{n}{2}+1)q}}\Big(1+|x|^q(1-r^2)^{\frac{q}{2}}\Big)dr\right)^{\frac{1}{q}}\\
    &\hspace{-2cm}\leq Ce^{|y|^2-|x|^2}\left(\left(\int_0^{\frac{1}{2}}(-\log r)^{q-1}r^{nq-1}e^{-c|y-rx|^2}(1+|x|^q)dr\right)^{\frac{1}{q}}\right.\\
    &\hspace{-2cm}\quad +\left.\Big(\int_{\frac{1}{2}}^1\frac{e^{-c\frac{|y-rx|^2}{1-r}}}{(1-r)^{\frac{n}{2}q+1}}\Big(1+|x|^q(1-r)^{\frac{q}{2}}\Big)dr\Big)^{\frac{1}{q}}\right)\\
    &\hspace{-2cm}=K_0(x,y)+K_1(x,y),\quad x,y\in \mathbb{R}^n.
\end{align*}
Then,
$$
g_{1,0,\{T_t^\mathcal{A}\}_{t>0};\,{\rm glob}}^q(f)(x)\leq \mathbb{K}_0(|f|)(x)+\mathbb{K}_1(|f|)(x),\quad x\in \mathbb{R}^n,
$$
where 
$$
\mathbb{K}_j(f)(x)=\int_{\mathbb{R}^n}f(y)\mathcal{X}_{N_1^c}(x,y)K_j(x,y)dy,\quad x\in \mathbb{R}^n,\;j=0,1.
$$

We study $\mathbb{K}_j$, $j=0,1$. For every $x,y\in \mathbb{R}^n$ with $x\not=0$ we write $y=y_x+y_\perp$, where $y_x$ is parallel to $x$ and $y_\perp$ is orthogonal to $x$. We define $r_0(x,y)=\frac{|y|}{|x|}\cos \theta (x,y)$, where $\theta (x,y)\in [0,\pi )$ represents the angle between $x$ and $y$.

First we show that
\begin{equation}\label{K0}
K_0(x,y)\leq Ce^{|y|^2-|x|^2}\left(
        |x|^{-n}+|x|^{2-n}+ e^{-c|y_\perp|^2}|x|\Big(\frac{|y|}{|x|}\Big)^{n-1}\mathcal{X}_{|y|\leq 2|x|}(x,y)\right),\quad (x,y)\in N_1^c,
\end{equation}
and thus, by virtue of \cite[Lemmas 4.2 and 4.3]{BrSj}, we can conclude that $\mathbb{K}_0$ is bounded from $L^1(\mathbb{R}^n,\gamma_{-1})$ into $L^{1,\infty }(\mathbb{R}^n,\gamma _{-1})$.

We observe that $|y-rx|^2\geq |x|^2$, when $r\in (0,\frac{1}{2})$ and $|y|\geq 2|x|$. Then,
$$
K_0(x,y)\leq Ce^{|y|^2-|x|^2}(1+|x|)e^{-c|x|^2}\left(\int_0^{\frac{1}{2}}(-\log r)^{q-1}r^{nq-1}dr\right)^{\frac{1}{q}}\leq Ce^{|y|^2-|x|^2}|x|^{-n},
\quad |y|\geq 2|x|.
$$
Consider now $(x,y)\in N_1^c$, with $|y|<2|x|$. Since $|y-rx|^2=(r-r_0)^2|x|^2+|y_\perp|^2$, we have that
\begin{align*}
    K_0(x,y)&\leq Ce^{|y|^2-|x|^2}(1+|x|)e^{-c|y_\perp|^2}\\
    &\quad \times \left(\int_0^{\frac{1}{2}}(-\log r)^{q-1}(|r_0|^{n-1}+|r-r_0|^{n-1})^qe^{-c(r-r_0)^2|x|^2}r^{q-1}dr\right)^{\frac{1}{q}}\\
    &\leq Ce^{|y|^2-|x|^2}(1+|x|)e^{-c|y_\perp|^2}\left(\Big(\frac{|y|}{|x|}\Big)^{n-1}+|x|^{1-n}\right)\left(\int_0^{\frac{1}{2}}(-\log r)^{q-1}r^{q-1}dr\right)^{\frac{1}{q}}\\
    &\leq Ce^{|y|^2-|x|^2}e^{-c|y_\perp|^2}\left(|x|\Big(\frac{|y|}{|x|}\Big)^{n-1}+|x|^{2-n}\right).
\end{align*}
In the last inequality we have taken into account that $|x|\geq C$, when $(x,y)\in N_1^c$ and $|y|<2|x|$. Property \eqref{K0} is then established.

Next we deal with $K_1(x,y)$, $(x,y)\in N_1^c$. We make the following decomposition (see \cite[p.12]{BrSj})
\begin{align*}
(K_1(x,y))^q&=Ce^{(|y|^2-|x|^2)q}\\
&\quad \times \int_{\frac{1}{2}}^1\frac{e^{-c\frac{(r-r_0)^2|x|^2+|y_\perp|^2}{1-r}}}{(1-r)^{\frac{nq}{2}+1}}\Big(1+|x|^q(1-r)^{\frac{q}{2}}\Big)\Big(\mathcal{X}_{\{r_0\not \in (\frac{1}{3},2)\}}(x,y)+\mathcal{X}_{\{r_0\in (\frac{1}{3},2)\}}(x,y)\Big)dr\\
&=(K_{1,1}(x,y))^q+(K_{1,2}(x,y))^q,\quad (x,y)\in N_1^c.
\end{align*}

Let $(x,y)\in N_1^c$. Assume first that $r_0=r_0(x,y)\leq \frac{1}{3}$ or $r_0\geq 2$. In these cases, $|r-r_0|\geq c$, for each $r\in (\frac{1}{2},1)$. Then, it follows that
\begin{align*}
    K_{1,1}(x,y)&\leq Ce^{|y|^2-|x|^2}\left(\int_{\frac{1}{2}}^1\frac{e^{-c\frac{|x|^2}{1-r}}}{(1-r)^{\frac{nq}{2}+1}}(1+|x|^q(1-r)^{\frac{q}{2}})dr\right)^{\frac{1}{q}}\\
    &\leq Ce^{|y|^2-|x|^2}\left(\int_{\frac{1}{2}}^1\frac{e^{-c\frac{|x|^2}{1-r}}}{(1-r)^{\frac{nq}{2}+1}}dr\right)^{\frac{1}{q}}\leq Ce^{|y|^2-|x|^2}|x|^{-n}.
\end{align*}

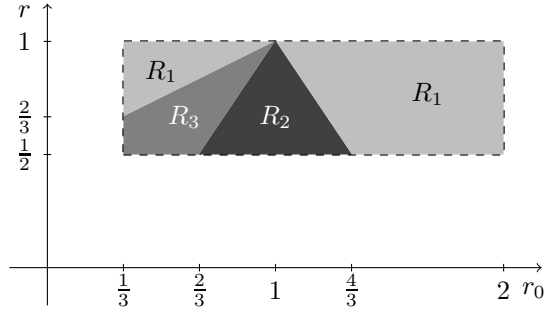
\begin{figure}[htbp]\label{Fig1}
\centering
\begin{minipage}{.5\textwidth}
  \centering
  \begin{tikzpicture}[scale=1]
    \draw[->] (-0.5,0) -- (6.5,0) ;
    \draw[->] (0,-0.5) -- (0,3.5) ;
    \draw (6.4,-0.3) node {$r_0$};
    \draw (-0.3,3.4) node {$r$};
    
    \draw[-] (1,-0.05) -- (1,0.05);
  	\draw (1,-0.3) node {$\frac{1}{3}$};  
    \draw[-] (2,-0.05) -- (2,0.05);
    \draw (2,-0.3) node {$\frac{2}{3}$};  
    \draw[-] (3,-0.05) -- (3,0.05);
    \draw (3,-0.3) node {$1$};
    \draw[-] (4,-0.05) -- (4,0.05);
    \draw (4,-0.3) node {$\frac{4}{3}$};
    \draw[-] (6,-0.05) -- (6,0.05);
    \draw (6,-0.3) node {$2$};
          
    \draw[-] (-0.05, 1.5) -- (0.05,1.5);     
    \draw (-0.3,1.5) node {$\frac12$};
    \draw[-] (-0.05,3) -- (0.05, 3);
    \draw (-0.3,3) node {$1$};
    \draw[-] (-0.05,2) -- (0.05, 2);
    \draw (-0.3,2) node {$\frac{2}{3}$};
     
	\draw[dashed, thick]  (1,1.5) -- (6,1.5) -- (6,3) -- (1,3) -- (1,1.5);  
	\draw[-, thick] (1,2) -- (3,3);
	\draw[-, thick] (2,1.5) -- (3,3);
	\draw[-, thick] (3,3) -- (4,1.5);
	\fill[color=gray]
        (1,1.5) -- (2,1.5) -- (3,3) -- (1,2);
        
    \fill[color=lightgray]
        (1,2) -- (3,3) -- (1,3);
    \fill[color=lightgray]
        (4,1.5) -- (6,1.5) -- (6,3) -- (3,3);
     \fill[color=darkgray]
        (2,1.5) -- (4,1.5) -- (3,3);
        
    \draw (1.5,2.6) node {$R_1$};
    \draw (5,2.3) node {$R_1$};
    \draw (3,2) node {\textcolor{white}{$R_2$}};
    \draw (1.8,2) node {\textcolor{white}{$R_3$}};
 \end{tikzpicture}
  \label{fig:Fig1}
\end{minipage}
\end{figure}

Suppose now that $r_0\in (\frac{1}{3},2)$. Then $|y_x|\sim |x|$ and $|y_\perp|\geq |x|\sin \theta (x,y)$. As in \cite[p.12]{BrSj} we split the interval $(\frac{1}{2},1)$ in the following three parts according to the regions in Figure \ref{Fig1}:
\begin{enumerate}
    \item $I_1=\big\{r\in (\frac{1}{2},1):1-r\leq \max\{\frac{1}{2}(1-r_0),\frac{3}{2}(r_0-1)\}\big\}$;\\
    \item $I_2=\big\{r\in (\frac{1}{2},1): |1-r_0|<\frac{2}{3}(1-r)\big\}$;\\
    \item $I_3=\big\{r\in (\frac{1}{2},1):|r-r_0|\leq \frac{1}{2}(1-r_0),\,r_0\in (\frac{1}{3},1)\big\}$.
\end{enumerate}

When $(r_0,r)\in R_1$, we have that $|r-r_0|\geq c(1-r)$ and $|r-r_0|\sim |1-r_0|=\frac{|x-y_x|}{|x|}$. Then,
\begin{align*}
    \int_{I_1}\frac{e^{-c\frac{(r-r_0)^2|x|^2+|y_\perp|^2}{1-r}}}{(1-r)^{\frac{nq}{2}+1}}\Big(1+|x|^q(1-r)^{\frac{q}{2}}\Big)dr&\leq C\int_{I_1}\frac{e^{-c(1-r)|x|^2}e^{-c\frac{|x-y_x|^2+|y_\perp|^2}{1-r}}}{(1-r)^{\frac{nq}{2}+1}}\Big(1+|x|^q(1-r)^{\frac{q}{2}}\Big)dr\\
    &\hspace{-4.5cm}\leq C\int_0^1\frac{e^{-c\frac{|x-y_x|^2+|y_\perp|^2}{1-r}}}{(1-r)^{\frac{nq}{2}+1}}dr\leq \frac{C}{(|x-y_x|+|y_\perp|)^{nq}}\leq C\min\big\{|x-y|^{-nq},|y_\perp|^{-nq}\big\}\\
    &\hspace{-4.5cm}\leq C\min \big\{(1+|x|)^{nq},|x\sin \theta (x,y)|^{-nq}\big\}.
\end{align*}

If $(r_0,r)\in R_2$, then $|r-r_0|\sim 1-r$ and we can write
\begin{align*}
    \int_{I_2}\frac{e^{-c\frac{(r-r_0)^2|x|^2+|y_\perp|^2}{1-r}}}{(1-r)^{\frac{nq}{2}+1}}\Big(1+|x|^q(1-r)^{\frac{q}{2}}\Big)dr&\leq C\int_{I_2}\frac{e^{-c(1-r)|x|^2}e^{-c\frac{|y_\perp|^2}{1-r}}}{(1-r)^{\frac{nq}{2}+1}}\Big(1+|x|^q(1-r)^{\frac{q}{2}}\Big)dr\\
    &\hspace{-4.5cm}\leq C\int_{I_2}\frac{e^{-c(1-r)|x|^2}e^{-c\frac{|y_\perp|^2}{1-r}}}{(1-r)^{\frac{nq}{2}+1}}dr.
\end{align*}
We have that
$$
\int_{I_2}\frac{e^{-c(1-r)|x|^2}e^{-c\frac{|y_\perp|^2}{1-r}}}{(1-r)^{\frac{nq}{2}+1}}dr\leq C\min\big\{(1+|x|)^{nq},|x\sin \theta (x,y)|^{-nq}\big\}.
$$
Indeed, 
$$
\int_{I_2}\frac{e^{-c(1-r)|x|^2}e^{-c\frac{|y_\perp|^2}{1-r}}}{(1-r)^{\frac{nq}{2}+1}}dr\leq C\int_0^1\frac{e^{-c\frac{|y_\perp|^2}{1-r}}}{(1-r)^{\frac{nq}{2}+1}}dr\leq \frac{C}{|y_\perp|^{nq}}\leq C|x\sin \theta (x,y)|^{-nq}.
$$
On the other hand, by making the change of variables $s=(1-r)|x|^2$, $r\in I_2$, we obtain 
$$
\int_{I_2}\frac{e^{-c(1-r)|x|^2}e^{-c\frac{|y_\perp|^2}{1-r}}}{(1-r)^{\frac{nq}{2}+1}}dr=C|x|^{nq}\int_{s>\frac{3}{2}|x||x-y_x|}\frac{e^{-cs}e^{-c\frac{|y_\perp|^2|x|^2}{s}}}{s^{\frac{nq}{2}+1}}ds.
$$
Since $(x,y)\in N_1^c$, it follows that $|x||x-y|\geq c$ and then
\begin{equation}\label{conseqN1c}
|x||x-y_x|\geq c,\mbox{ when }|y_\perp|\leq |x-y_x|,\quad \mbox{ and }\quad |y_\perp||x|\geq c,\mbox{ when }|x-y_x|\leq |y_\perp|.
\end{equation}
We get
\begin{align*}
\int_{I_2}\frac{e^{-c(1-r)|x|^2}e^{-c\frac{|y_\perp|^2}{1-r}}}{(1-r)^{\frac{nq}{2}+1}}dr&\leq C|x|^{nq}\times\left\{
    \begin{array}{ll}
        \displaystyle \int_0^\infty e^{-cs}ds, & \mbox{ if }|y_\perp|\leq |x-y_x| \\[0.3cm]
         \displaystyle \int_0^\infty \frac{e^{-\frac{c}{s}}}{s^{\frac{nq}{2}+1}}ds,&\mbox{ if }|x-y_x|\leq |y_\perp| 
    \end{array}\right.\leq C(1+|x|)^{nq}.
\end{align*}

Finally, we analyze the integral when $(r_0,r)\in R_3$. In this case, $1-r\sim 1-r_0=\frac{|x-y_x|}{|x|}$. So, by using \cite[(5.3)]{BrSj} we get
\begin{align*}
 \int_{I_3}\frac{e^{-c\frac{(r-r_0)^2|x|^2+|y_\perp|^2}{1-r}}}{(1-r)^{\frac{nq}{2}+1}}\big(1+|x|^q(1-r)^{\frac{q}{2}}\big)dr&\leq C\frac{1+|x|^q(1-r_0)^{\frac{q}{2}}}{(1-r_0)^{\frac{nq}{2}+1}}e^{-c\frac{|y_\perp|^2}{1-r_0}}\int_{I_3}e^{-c\frac{(r-r_0)^2|x|^2}{1-r_0}}dr\\
 &\hspace{-4.5cm}\leq C\frac{1+|x|^q(1-r_0)^{\frac{q}{2}}}{(1-r_0)^{\frac{nq}{2}+1}}e^{-c\frac{|y_\perp|^2}{1-r_0}}\min\Big\{1-r_0,\frac{\sqrt{1-r_0}}{|x|}\Big\}\leq C\frac{1+|x|^q(1-r_0)^{\frac{q}{2}}}{(1-r_0)^{\frac{nq}{2}}}e^{-c\frac{|y_\perp|^2}{1-r_0}}\\
 &\hspace{-4.5cm} =C\frac{|x|^{\frac{nq}{2}}(1+(|x||x-y_x|)^{\frac{q}{2}})}{|x-y_x|^{\frac{nq}{2}}}e^{-c\frac{|y_\perp|^2|x|}{|x-y_x|}}.
\end{align*}
When $|x-y_x|\leq |y_\perp|$, we have that $|x-y|\leq 2|y_\perp|$, and since $(x,y)\in N_1^c$ we get 
\begin{align*}
\frac{|x|^{\frac{nq}{2}}(1+(|x||x-y_x|)^{\frac{q}{2}})}{|x-y_x|^{\frac{nq}{2}}}e^{-c\frac{|y_\perp|^2|x|}{|x-y_x|}}&\leq C\frac{|x|^{\frac{nq}{2}}}{|x-y_x|^{\frac{nq}{2}}}e^{-c\frac{|y_\perp|^2|x|}{|x-y_x|}}\leq \frac{C}{|y_\perp|^{nq}}\\
&\hspace{-4cm}\leq \min\Big\{\frac{1}{|x-y|^{nq}},\frac{1}{|x\sin \theta (x,y)|^{nq}}\Big\}\leq C\min\Big\{(1+|x|)^{nq}, |x\sin \theta(x,y)|^{-nq}\Big\}.
\end{align*}
Assume now that $|y_\perp|\leq |x-y_x|$. From \eqref{conseqN1c}, we obtain that
$$
\frac{|x|^{\frac{nq}{2}}(1+(|x||x-y_x|)^{\frac{q}{2}})}{|x-y_x|^{\frac{nq}{2}}}e^{-c\frac{|y_\perp|^2|x|}{|x-y_x|}}\leq C\left(\frac{|x|^{\frac{nq}{2}}}{|x-y_x|^{\frac{nq}{2}}}+\frac{|x|^{\frac{(n+1)q}{2}}}{|x-y_x|^{\frac{(n-1)q}{2}}}\right)\leq C|x|^{nq}.
$$
Also, we can write
\begin{align*}
    \frac{|x|^{\frac{nq}{2}}(1+(|x||x-y_x|)^{\frac{q}{2}})}{|x-y_x|^{\frac{nq}{2}}}e^{-c\frac{|y_\perp|^2|x|}{|x-y_x|}}&\leq 
    C\left(\frac{1}{|y_\perp|^{nq}}+\frac{|x|^{\frac{(n+1)q}{2}}}{|x-y_x|^{\frac{(n-1)q}{2}}}e^{-c\frac{|y_\perp|^2|x|}{|x-y_x|}}\right)\\
    &\hspace{-2cm}\leq C\left(\frac{1}{|x\sin \theta (x,y)|^{nq}}+\frac{|x|^{\frac{(n+1)q}{2}}}{|x-y_x|^{\frac{(n-1)q}{2}}}e^{-c\frac{|y_\perp|^2|x|}{|x-y_x|}}\mathcal{X}_{\{|x||x-y_x|\geq 1\}}(x,y)\right).
\end{align*}
In the last inequality we have taken into account that if $|x||x-y_x|\leq 1$, then 
$$
\frac{|x|^{\frac{(n+1)q}{2}}}{|x-y_x|^{\frac{(n-1)q}{2}}}e^{-c\frac{|y_\perp|^2|x|}{|x-y_x|}}\leq C\frac{(|x||x-y_x|)^{\frac{q}{2}}}{|y_\perp|^{nq}}\leq \frac{C}{|y_\perp|^{nq}}.
$$
By collecting all the estimates above whe have obtained that
\begin{align*}
\int_{I_3}\frac{e^{-c\frac{(r-r_0)^2|x|^2+|y_\perp|^2}{1-r}}}{(1-r)^{\frac{nq}{2}+1}}\big(1+|x|^q(1-r)^{\frac{q}{2}}\big)dr&\\
&\hspace{-5.5cm}\leq C\left(\min\Big\{(1+|x|)^{nq}, |x\sin \theta(x,y)|^{-nq}\Big\}+\frac{|x|^{\frac{(n+1)q}{2}}}{|x-y_x|^{\frac{(n-1)q}{2}}}e^{-c\frac{|y_\perp|^2|x|}{|x-y_x|}}\mathcal{X}_{\{|x||x-y_x|\geq 1\}}(x,y)\right),
\end{align*}
and we deduce that
\begin{align*}
    K_1(x,y)&\leq Ce^{|y|^2-|x|^2}\left(|x|^{-n}+\min\Big\{(1+|x|)^{n}, |x\sin \theta(x,y)|^{-n}\Big\}\right.\\
    &\left.\quad +\frac{|x|^{\frac{n+1}{2}}}{|x-y_x|^{\frac{n-1}{2}}}e^{-c\frac{|y_\perp|^2|x|}{|x-y_x|}}\mathcal{X}_{\{|x||x-y_x|\geq 1\}}(x,y)\right),\quad (x,y)\in N_1^c.
\end{align*}
By taking into account \cite[Lemma 3.3.4]{Sa} and \cite[Lemmas 4.2 and 4.4]{BrSj} we conclude that the operator $\mathbb{K}_1$ is bounded from $L^1(\mathbb{R}^n,\gamma_{-1})$ into $L^{1,\infty }(\mathbb{R}^n,\gamma_{-1})$.

Hence, $g_{1,0,\{T_t^\mathcal{A}\}_{t>0};\,{\rm glob}}$ defines a bounded operator from $L^1(\mathbb{R}^n,\gamma_{-1})$ into $L^{1,\infty }(\mathbb{R}^n,\gamma_{-1})$. We have proved that $g_{1,0,\{T_t^\mathcal{A}\}_{t>0}}$ is bounded from $L^1(\mathbb{R}^n,\gamma_{-1})$ into $L^{1,\infty }(\mathbb{R}^n,\gamma_{-1})$.

\subsection{Proof of Theorem \ref{Th1.1} (ii), for 1$<$p$<\infty$}
We recall that
$$
P_t^\mathcal{A}(f)=\frac{1}{\sqrt{\pi}}\int_0^\infty\frac{e^{-s}}{\sqrt{s}}T_{t^2/(4s)}^\mathcal{A}(f)ds,\quad t>0.
$$
As in ($i$) it is sufficient to prove ($ii$) with $\beta \in \mathbb{N}$. Let $f\in L^p(\mathbb{R}^n,\gamma_{-1})$ and $m\in \mathbb{N}$. We have that
$$
\partial _t^m\partial _x^kP_t^\mathcal{A}(f)(x)=\frac{1}{\sqrt{\pi}}\int_0^\infty \frac{e^{-s}}{\sqrt{s}}\partial ^k_x\partial _t^mT_{t^2/(4s)}^\mathcal{A}(f)(x)ds,\quad x\in \mathbb{R}^n,\;t>0.
$$
By using Fa\`a di Bruno's formula we get
\begin{align*}
\partial _t^mT_{t^2/(4s)}^\mathcal{A}(f)(x)&=\frac{1}{2^ms^{\frac{m}{2}}}\partial _v^mT_{v^2}^\mathcal{A}(f)(x)_{\big|v=\frac{t}{2\sqrt{s}}}\\
&=\frac{1}{s^{\frac{m}{2}}}\sum_{\ell \in \mathbb{N},\;0\leq \ell \leq \frac{m}{2}}c_{\ell ,m}\Big(\frac{t}{\sqrt{s}}\Big)^{m-2\ell}\partial _v^{m-\ell}T_v^\mathcal{A}(f)(x)_{\big|v=\frac{t^2}{4s}},\quad x\in \mathbb{R}^n,\;t,s>0.
\end{align*}
Here $c_{\ell, m}\in \mathbb{R}$, $\ell\in \mathbb{N}$, $0\leq \ell \leq \frac{m}{2}$.

We can write
$$
\partial _t^m\partial _x^kP_t^\mathcal{A}(f)(x)=\sum_{\ell \in \mathbb{N},\;0\leq \ell \leq \frac{m}{2}}c_{\ell ,m}\int_0^\infty \frac{e^{-s}}{s^{\frac{m+1}{2}}}\Big(\frac{t}{\sqrt{s}}\Big)^{m-2\ell}\partial_x^k\partial _v^{m-\ell}T_v^\mathcal{A}(f)(x)_{\big|v=\frac{t^2}{4s}}ds,\quad x\in \mathbb{R}^n,\;t>0.
$$

By using Minkowski inequality we obtain
\begin{align}\label{gPoissongHeat}
    g_{m,k,\{P_t^\mathcal{A}\}_{t>0}}^q(f)(x)&=\left\|t^{m+|k|}\partial _t^m\partial _x^kP_t^\mathcal{A}(f)(x)\right\|_{L^q((0,\infty ),\frac{dt}{t})}\nonumber\\
    &\leq C\sum_{\ell \in \mathbb{N},\;0\leq \ell \leq \frac{m}{2}}\int_0^\infty \frac{e^{-s}}{s^{m-\ell+\frac{1}{2}}}\left\|t^{2(m-\ell)+|k|}\partial_x^k\partial _v^{m-\ell}T_v^\mathcal{A}(f)(x)_{\big|v=\frac{t^2}{4s}}\right\|_{L^q((0,\infty ),\frac{dt}{t})}ds\nonumber\\
    &\leq C\sum_{\ell \in \mathbb{N},\;0\leq \ell \leq \frac{m}{2}}\int_0^\infty e^{-s}s^{\frac{|k|-1}{2}}\big\|v^{m-\ell +\frac{|k|}{2}}\partial_x^k\partial _v^{m-\ell}T_v^\mathcal{A}(f)(x)\big\|_{L^q((0,\infty),\frac{dt}{t})} ds\nonumber\\
    &\leq C\sum_{\ell \in \mathbb{N},\;0\leq \ell \leq \frac{m}{2}}g_{m-\ell, k,\{T_t^\mathcal{A}\}_{t>0}}^q(f)(x),\quad x\in \mathbb{R}^n.
\end{align}

By using the property established in Section \ref{S2.1} we conclude that $g_{m,k,\{P_t^\mathcal{A}\}_{t>0}}^q$ is bounded from $L^p(\mathbb{R}^n,\gamma_{-1})$ into itself, for every $1<p<\infty$.

\subsection{Proof of Theorem \ref{Th1.1} (ii), for p$=$1}
Let $m\in \mathbb{N}$ and $k\in \mathbb{N}^n$. Our next objective is to see that $g_{m,k,\{P_t^\mathcal{A}\}_{t>0}}^q$ is bounded from $L^1(\mathbb{R}^n,\gamma_{-1})$ into $L^{1,\infty }(\mathbb{R}^n,\gamma_{-1})$ when $|k|\leq 2$.

In Section \ref{S2.1} we proved that $G_{m,k,\{T_t^\mathcal{A}\}_{t>0}}$ satisfies the properties (a), (b) and (c). Then, from the $L^q((0,\infty),\frac{dt}{t})$-version of \cite[Theorem 3.2.8]{Sa} we deduce that $g_{m,k,\{T_t^\mathcal{A}\}_{t>0};\,{\rm loc}}^q$ is bounded from $L^1(\mathbb{R}^n,\gamma_{-1})$ into $L^{1,\infty }(\mathbb{R}^n,\gamma_{-1})$.

By \eqref{gPoissongHeat} it follows that $g_{m,k,\{P_t^\mathcal{A}\}_{t>0};\,{\rm loc}}^q$ is bounded from $L^1(\mathbb{R}^n,\gamma_{-1})$ into $L^{1,\infty }(\mathbb{R}^n,\gamma_{-1})$.

We now study the operator $g_{m,k,\{P_t^\mathcal{A}\}_{t>0};\,{\rm glob}}^q$. We can write
$$
P_t^\mathcal{A}(f)=\frac{t}{2\sqrt{\pi}}\int_0^\infty \frac{e^{-\frac{t^2}{4s}}}{s^{\frac{3}{2}}}T_s^{\mathcal A}(f)ds=-\frac{1}{\sqrt{\pi}}\int_0^\infty \partial_t(e^{-\frac{t^2}{4s}})T_s^\mathcal{A}f\frac{ds}{\sqrt{s}},\quad t>0.
$$
Let $f\in L^1(\mathbb{R}^n,\gamma_{-1})$. By using Minkowski inequality we get
\begin{align*}
    \big\|t^{m+|k|}\partial _t^m\partial_x^kP_t^\mathcal{A}(f)(x)\big\|_{L^q(0,\infty),\frac{dt}{t})}\leq C\int_0^\infty\big\| t^{m+|k|}\partial_t^{m+1}(e^{-\frac{t^2}{4s}})\big\|_{L^q(0,\infty),\frac{dt}{t})}|\partial_x^kT_s^\mathcal{A}f(x)|\frac{ds}{\sqrt{s}}.
\end{align*}
Formula de Fa\`a di Bruno leads to 
$$
\partial_t^{m+1}(e^{-\frac{t^2}{4s}})=\frac{e^{-\frac{t^2}{4s}}}{s^{\frac{m+1}{2}}}\sum_{\ell \in \mathbb{N},\,0\leq \ell \leq\frac{m+1}{2}}c_{m,\ell}\Big(\frac{t^2}{s}\Big)^{\frac{m+1}{2}-\ell},\quad t,s>0,
$$
for certain $c_{m,\ell}\in \mathbb{R}$, $0\leq \ell \leq (m+1)/2$. Then,
$$
\big\| t^{m+|k|}\partial_t^{m+1}(e^{-\frac{t^2}{4s}})\big\|_{L^q((0,\infty),\frac{dt}{t})}^q\leq C\sum_{\ell \in \mathbb{N},\,0\leq \ell \leq\frac{m+1}{2}}\int_0^\infty \frac{e^{-\frac{qt^2}{4s}}t^{(2m+|k|+1-2\ell)q-1}}{s^{(m+1-\ell)q}} dt\leq Cs^{\frac{(|k|-1)q}{2}},\quad s>0,
$$
and we obtain 
\begin{align*}
    \big\|t^{m+|k|}\partial _t^m\partial_x^kP_t^\mathcal{A}(f)(x)\big\|_{L^q(0,\infty),\frac{dt}{t})}&\leq C\int_0^\infty s^{\frac{|k|}{2}-1}|\partial_x^kT_s^\mathcal{A}f(x)|ds\\
    &\leq C\int_{\mathbb{R}^n}|f(y)|\int_0^\infty s^{\frac{|k|}{2}-1}|\partial_x^kT_s^\mathcal{A}(x,y)|dsdy,\quad x\in \mathbb{R}^n.
\end{align*}
According to \cite[Proposition 5.1]{BrSj} we deduce that the operator
\begin{equation}\label{operatorL}
\mathcal{L}(f)(x)=\int_{\mathbb{R}^n}f(y)\mathcal{X}_{N_1^c}(x,y)\int_0^\infty s^{\frac{|k|}{2}-1}|\partial_x^kT_s^\mathcal{A}(x,y)|dsdy,\quad x\in \mathbb{R}^n,
\end{equation}
is bounded from $L^1(\mathbb{R}^n,\gamma_{-1})$ into $L^{1,\infty}(\mathbb{R}^n,\gamma_{-1})$ when $|k|\leq 2$. Note that the proof of \cite[Proposition 5.1]{BrSj} also works for $|k|=0$.

Hence, the operator $g_{m,k,\{P_t^\mathcal{A}\}_{t>0};\,{\rm glob}}^q$ is bounded from $L^1(\mathbb{R}^n,\gamma_{-1})$ into $L^{1,\infty}(\mathbb{R}^n,\gamma_{-1})$ provided that $|k|\leq 2$. We conclude that $g_{m,k,\{P_t^\mathcal{A}\}_{t>0}}^q$ is bounded from $L^1(\mathbb{R}^n,\gamma_{-1})$ into $L^{1,\infty}(\mathbb{R}^n,\gamma_{-1})$ when $|k|\leq 2$.

We are going to see that $g_{m,k,\{P_t^\mathcal{A}\}_{t>0}}^q$ is not bounded from $L^1(\mathbb{R}^n,\gamma_{-1})$ into $L^{1,\infty}(\mathbb{R}^n,\gamma_{-1})$ when $|k|\geq 3$. We use some ideas developed in \cite[\S 6]{BrSj}.

Let $\eta>0$. We define $z=(\eta,...,\eta)$. We will choose $\eta$ large. For every $x\in \mathbb{R}^n$ we write $x=x_z+x_{\perp}$, where $x_z$ is parallel to $z$ and $x_\perp$ is orthogonal to $z$. Let us consider the set
$$
J(z)=\Big\{x\in \mathbb{R}^n: |x_\perp|<1, \,\frac{4}{3}|z|<|x_z|<\frac{3}{2}|z|\Big\}.
$$
We have that (see \cite[\S 6]{BrSj}) there exists $\eta _0>0$ such that when $\eta >\eta_0$ we have that, for every $r\in (0,1)$, $y\in B(z,1)$ and $x\in J(z)$,
$$
\frac{x_i-ry_i}{\sqrt{1-r^2}}\geq C|z|,\quad i=1,...,n,
$$
and hence
$$
H_k\left(\frac{x-ry}{\sqrt{1-r^2}}\right)\geq C|z|^k.
$$

Also we get
$$
e^{-\frac{|y-rx|^2}{1-r^2}}\geq Ce^{-c|r|x_z|-|z||^2},\quad x\in J(z),\,y\in B(z,1)\mbox{ and }r\in \Big(\frac{1}{4},\frac{3}{4}\Big).
$$

We choose $0<t_0<t_1<\infty$ such that for a certain $C>0$
$$
\Big|H_{m+1}\Big(\frac{t}{2(-\log r)^{\frac{1}{2}}} \Big)\Big|\geq C,\quad t\in (t_0,t_1)\mbox{ and }r\in \Big(\frac{1}{4},\frac{3}{4}\Big),
$$
and the sign of $H_{m+1}\Big(\frac{t}{2(-\log r)^{\frac{1}{2}}} \Big)$ is constant when $t\in (t_0,t_1)$ and $r\in (\frac{1}{4},\frac{3}{4})$. We write 
$$
a_0=-{\rm sign}\;H_{m+1}\Big(\frac{t}{2(-\log r)^{\frac{1}{2}}}\Big),\quad t\in (t_0,t_1)\mbox{ and }r\in \Big(\frac{1}{4},\frac{3}{4}\Big).
$$
We now take a function $0\leq f\in L^1(\mathbb{R}^n,\gamma_{-1})$ such that $\mbox{ supp}\,f\subset B(z,1)$ and $\|f\|_{L^1(\mathbb{R}^n,\gamma_{-1})}=1$.

We have that
$$
g_{m,k,\{P_t^\mathcal{A}\}_{t>0}}^q(f)(x)=\Big\|t^{m+|k|}\partial _t^m\partial_x^kP_t^\mathcal{A}(f)(x)\Big\|_{L^q((0,\infty ),\frac{dt}{t})}=\Big\|\int_{\mathbb{R}^n}f(y)H_{k,m}(t,x,y)dy\Big\|_{L^q((0,\infty ),\frac{dt}{t})},
$$
where, for each $x,y\in \mathbb{R}^n$ and $t>0$,
\begin{align*}
H_{k,m}(t,x,y)&=-\frac{t^{m+|k|}}{\sqrt{\pi}}\int_0^\infty \partial _t^{m+1}(e^{-\frac{t^2}{4s}})\partial _x^kT_s^\mathcal{A}(x,y)\frac{ds}{s}\\
&=\frac{(-t)^{m+|k|}}{2^{m+1}\sqrt{\pi}}\int_0^\infty H_{m+1}\Big(\frac{t}{2\sqrt{s}}\Big)\frac{e^{-\frac{t^2}{4s}}}{s^{m+\frac{3}{2}}}e^{-ns}H_k\Big(\frac{x-e^{-s}y}{\sqrt{1-e^{-2s}}}\Big)\frac{e^{-\frac{|x-e^{-s}y|^2}{1-e^{-2s}}}}{(1-e^{-2s})^{\frac{n+|k|}{2}}}ds.
\end{align*}
By making $r=e^{-s}$, $s\in (0,\infty)$, it follows that, when $t\in (t_0,t_1)$, $y\in B(z,1)$, $\eta>\eta_0$ and $x\in J(z)$,
\begin{align*}
a_0H_{k,m}(t,x,y)&=\frac{(-t)^{m+|k|}}{2^{m+1}\sqrt{\pi}}\int_0^1H_{m+1}\Big(\frac{t}{2(-\log r)^{\frac{1}{2}}}\Big)\frac{e^{-\frac{t^2}{4(-\log r)}}}{(-\log r)^{m+\frac{3}{2}}}r^{n-1}H_k\Big(\frac{x-ry}{\sqrt{1-r^2}}\Big)\frac{e^{-\frac{|x-ry|^2}{1-r^2}}}{(1-r^2)^{\frac{n+|k|}{2}}}dr\\
&\geq Ct^{m+|k|}\eta ^{|k|}e^{-ct^2}e^{|y|^2-|x|^2}\int_{1/4}^{3/4}e^{-c|r|x_z|-|z||^2}dr.
\end{align*}
Then,
\begin{align*}
    g_{m,k,\{P_t^\mathcal{A}\}_{t>0}}^q(f)(x)&\geq C\left(\int_{t_0}^{t_1}t^{(m+|k|)q-1}e^{-ct^2}dt\right)^{\frac{1}{q}}\eta ^{|k|}e^{-|x|^2}\int_{1/4}^{3/4}e^{-c|r|x_z|-|z||^2}dr\\
    &\geq C\eta ^{|k|-1}e^{-|x|^2}\geq C|z|^{k-1}e^{-(\frac{3}{2}|z|)^2},\quad x\in J(z),\;\eta >\eta_0.
\end{align*}
Moreover, $\gamma_{-1}(J(z))\geq e^{(\frac{3}{2}|z|)^2}|z|^{-1}$. We get
$$
\sup_{s>0}s\gamma_{-1}(\{x\in \mathbb{R}^n: g_{m,k,\{P_t^\mathcal{A}\}_{t>0}}^q(f)(x)>s\})\geq Ce^{-(\frac{3}{2}|z|)^2}|z|^{|k|-1}\gamma_{-1}(J(z))\geq C|z|^{|k|-2},\quad \eta >\eta_0.
$$
We conclude that $g_{m,k,\{P_t^\mathcal{A}\}_{t>0}}^q$ is not bounded from $L^1(\mathbb{R}^n,\gamma_{-1})$ into $L^{1,\infty }(\mathbb{R}^n,\gamma_{-1})$ when $|k|>2$.

\noindent {\bf Remark}. As it was mentioned during the last proof we have seen that the operator $g_{m,k,\{T_t^\mathcal{A}\}_{t>0}}^q$ is bounded from $L^1(\mathbb{R}^n,\gamma_{-1})$ into $L^{1,\infty }(\mathbb{R}^n,\gamma_{-1})$ except when $m=1$ and $k=0$. As far as we know this question is also open even when $q=2$, $k=0$ and $m=1$ in the Gaussian setting, that is, it is not known whether the Littlewood-Paley function $g_{1,0,\{T_t^\mathcal{L}\}_{t>0}}^2$ is bounded from  $L^1(\mathbb{R}^n,\gamma_{-1})$ into $L^{1,\infty }(\mathbb{R}^n,\gamma_{-1})$. We recall that $\mathcal{L}$ denotes the Ornstein-Uhlenbeck operator.


\section{Proof of Theorems \ref{Th1.2} and \ref{Th1.3} and Corollary \ref{cor1.4}}

\subsection{Proof of Theorem \ref{Th1.2}}
According to \cite[Theorem 2]{Xu2} and \cite[Theorem 1.5]{BFGM}, since $\{T_t^{\mathcal A}\}_{t>0}$ is a symmetric diffusion semigroup, (i) $\Rightarrow$ (ii) holds. By using \cite[Theorem A]{TZ} (see also \cite[Theorem 2.1]{MTX}) we can see that (i) $\Rightarrow$ (iii).

Suppose that (ii) holds. We are going to prove (i). According to \cite[Theorem 1.5 (b)]{BFGM}, there exists a norm $\VERT \cdot\VERT$ on $X$ that defines the original topology of $X$ such that $(X,\VERT\cdot\VERT)$ is $q$-uniformly convex, provided that $g_{\alpha ,0,\{W_t\}_{t>0}}^{q,X}$ is bounded from $L^p_X(\mathbb R^n,dx)$ into $L^p(\mathbb R^n,dx)$, for some $\alpha >0$. Here $\{W_t\}_{t>0}$ denotes the classical heat semigroup being, for every $t>0$, 

$$W_t(f)(x)=\frac{1}{\pi^{n/2}}\int_{\mathbb R^n}\frac{e^{-|x-y|^2/2t}}{(2t)^{n/2}}f(y)dy,\quad x\in\mathbb R^n.$$

Since $g_{\beta_1 ,0,\{W_t\}_{t>0}}^{q,X}(f)\leq g_{\beta_2 ,0,\{W_t\}_{t>0}}^{q,X}(f)$ (\cite[Proposition 3.1]{TZ}), $0<\beta_1<\beta_2$, we can assume that $\beta\in (0,1)$.

Our objective is to see that $g_{\beta ,0,\{W_t\}_{t>0}}^{q,X}$ is bounded from $L^p_X(\mathbb R^n,dx)$ into $L^p(\mathbb R^n,dx)$

We define the local and global operators in the usual way associated to $N_1$. Minkowski inequality leads to

\begin{align*}
g_{\beta ,0,\{W_t\}_{t>0};\,{\rm loc}}^{q,X}(f)(x) &\leq \left|g_{\beta ,0,\{W_t\}_{t>0};\,{\rm loc}}^{q,X}(f)(x)-g_{\beta ,0,\{T_t^{\mathcal A}\}_{t>0};\,{\rm loc}}^{q,X}(f)(x)\right|+g_{\beta ,0,\{T_t^{\mathcal A}\}_{t>0};\,{\rm loc}}^{q,X}(f)(x) \\
& \leq \mathcal H_1^\beta(f)(x)+g_{\beta ,0,\{T_t^{\mathcal A}\}_{t>0};\,{\rm \,{\rm loc}}}^{q,X}(f)(x),\;\;\;\;x\in\mathbb R^n.
\end{align*}
where 
$$
\mathcal H_1^\beta(f)(x)=\int_{\mathbb R^n}\|f(y)\|_X\left(\int_0^\infty |t^\beta\partial_t^\beta(W_t(x-y)-T_t^{\mathcal A}(x,y))|^q\frac{dt}{t}\right)^{1/q}\mathcal{X}_{N_1}(x,y)dy,\;\;x\in\mathbb R^n.
$$
Since $\beta\in (0,1)$, we can write

$$\partial_t^\beta(W_t(x-y)-T_t^{\mathcal A}(x,y))=\frac{1}{\Gamma(1-\beta)}\int_t^\infty \partial_u(W_u(x-y)-T_u^{\mathcal A}(x,y))(u-t)^{-\beta}du,\;\;\;x,y\in\mathbb R^n.$$
By proceeding as in the proof \cite[Proposition 3.1]{TZ} we get

$$\left(\int_0^\infty |t^\beta\partial_t^\beta(W_t(x-y)-T_t^{\mathcal A}(x,y))|^q\frac{dt}{t}\right)^{1/q}\leq C\left(\int_0^\infty |t\partial_t(W_t(x-y)-T_t^{\mathcal A}(x,y))|^q\frac{dt}{t}\right)^{1/q},\;\;\;x,y\in\mathbb R^n.$$
Then,

$$g_{\beta ,0,\{W_t\}_{t>0};\,{\rm loc}}^{q,X}(f)\leq C\left(\mathcal H_1^1(f)+g_{\beta ,0,\{T_t^{\mathcal A}\}_{t>0};\,{\rm loc}}^{q,X}(f)\right).$$
We consider the operator

$$G_{\beta ,\{T_t^{\mathcal A}\}_{t>0}}^{X}(f)(t,x)=t^\beta\partial_t^\beta T_t^{\mathcal A}(f)(x),\;\;\;\;x\in\mathbb R^n\;\mbox{and}\;t>0.$$
Since (ii) holds, $G_{\beta ,\{T_t^{\mathcal A}\}_{t>0}}^{X}$ is bounded from $L^p_X(\mathbb R^n,\gamma_{-1})$ into $L^p_{L^q_X((0,\infty),\frac{dt}{t})}(\mathbb R^n,\gamma_{-1})$. By the property (b) in the proof of Theorem \ref{Th1.1} we have that

$$\left(\int_0^\infty |t^\beta\partial_t^\beta T_t^{\mathcal A}(x,y)|^q\frac{dt}{t}\right)^{1/q}\leq\frac{C}{|x-y|^n},\quad (x,y)\in N_2.$$
A vector valued version of \cite[Proposition 3.2.7]{Sa} (see \cite[Proposition 2.3]{HTV}) allow us to obtain that $G_{\beta ,\{T_t^{\mathcal A}\}_{t>0};\,{\rm loc}}^{X}$ is bounded from $L^p_X(\mathbb R^n,dx)$ into $L^p_{L^q_X((0,\infty),\frac{dt}{t})}(\mathbb R^n,dx)$. or, in other words $g_{\beta ,0,\{T_t^{\mathcal A}\}_{t>0};\,{\rm loc}}^{q,X}$ is bounded from $L^p_X(\mathbb R^n,dx)$ into $L^p(\mathbb R^n,dx)$.

We now study the operator $\mathcal H_1^1$. We have that

\begin{align*}
  t\partial_ t(W_t(x-y) -T_t^{\mathcal A}(x,y))&=-nt\left(\frac{W_t(x-y)}{2t}-\frac{T_t^{\mathcal A}(x,y)}{1-e^{-2t}}\right) \\
 & \quad +2t\left( \frac{|x-y|^2}{(2t)^2}W_t(x-y)-\frac{e^{-2t}|x-e^{-t}y|^2}{(1-e^{-2t})^2}T_t^\mathcal{A}(x,y)\right) \\
 &\quad+2te^{-t}\frac{\langle y,x-e^{-t}y\rangle }{1-e^{-2t}}T_t^\mathcal{A}(x,y)=:\sum_{j=1}^3K_j(t,x,y), \;\;\;\;x,y\in\mathbb R^n,\;t>0.
\end{align*}

We are going to estimate $\|K_j(\cdot,x,y)\|_{L^q((0,\infty),\frac{dt}{t})}$, $j=1,2,3$, for every $(x,y)\in N_1$. By performing the change of variable $t=\log\frac{1+s}{1-s}$, $s\in (0,\infty)$, we obtain

\begin{align*}
 \|K_3(\cdot,x,y)\|_{L^q((0,\infty),\frac{dt}{t})}^q& \leq C|y|^q\int_0^\infty t^qe^{-q(n+1)t}\frac{e^{-c\frac{|x-e^{-t}y|^2}{1-e^{-2t}}}}{(1-e^{-2t})^{\frac{n+1}{2}q}}\frac{dt}{t}\\
 &\hspace{-3.2cm}\leq C|y|^qe^{-c(|x|^2-|y|^2)}\int_0^\infty e^{-qnt}\frac{e^{-c\frac{|y-e^{-t}x|^2}{1-e^{-2t}}}}{(1-e^{-2t})^{\frac{n-1}{2}q}}\frac{dt}{t} \\
 &\hspace{-3.2cm}\leq C|y|^qe^{-c(|x|^2-|y|^2)}\int_0^1\left(\frac{1-s}{1+s}\right)^{qn}\frac{e^{-c(s|x+y|^2+\frac{1}{s}|x-y|^2)}}{s^{\frac{n-1}{2}q}}\frac{ds}{(1-s)\log\frac{1+s}{1-s}}\\
 &\hspace{-3.2cm}\leq C|y|^qe^{-c(|x|^2-|y|^2)}\int_0^1\frac{e^{-c\frac{|x-y|^2}{s}}}{s^{\frac{n-1}{2}q}}e^{-cs(|x+y|^2+|x-y|^2)}\frac{ds}{-\log (1-s)},\quad x,y\in\mathbb R^n.
\end{align*}

Since $2|y|\leq (|x+y|+|x-y|)$, $x,y\in \mathbb{R}^n$, $-\log (1-s)\geq cs$, $s\in (0,\infty )$ and $||y|^2-|x|^2|\leq C$, $(x,y)\in N_1$, it follows that

$$
 \|K_3(\cdot,x,y)\|_{L^q((0,\infty),\frac{dt}{t})}^q\leq C|y|^{q/2} \int_0^1\frac{e^{-c\frac{|x-y|^2}{s}}}{s^{(\frac{n-1}{2}+\frac{1}{4})q+1}} ds\; 
  \leq C\frac{|y|^{q/2}}{|x-y|^{(n-1/2)q}},\quad (x,y)\in N_1.
$$
Then
$$
\|K_3(\cdot,x,y)\|_{L^q((0,\infty),\frac{dt}{t})}\leq C\frac{\sqrt{1+|x|}}{|x-y|^{n-1/2}},\quad (x,y)\in N_1.
$$

Next we deal with $\|K_j(\cdot,x,y)\|_{L^q((0,\infty),\frac{dt}{t})}$, $j=1,2$. Let us first observe that
\begin{align*}
|K_1(t,x,y)|+|K_2(t,x,y)|&\leq C\left(\Big(1+\frac{|x-y|^2}{t}\Big)W_t(x-y)+\frac{t}{1-e^{-2t}}\Big(1+\frac{|x-e^{-t}y|^2}{1-e^{-2t}}\Big)T_t^\mathcal{A}(x,y)\right)\\
& \leq C\left(\frac{1}{t^{\frac{n}{2}}}+\frac{te^{-nt}}{(1-e^{-2t})^{\frac{n}{2}+1}}\right)\leq C\frac{1}{t^{\frac{n}{2}}},\quad x,y\in \mathbb{R}^n,\;t>0.
\end{align*}
Then, by using that $|x-y|\leq n\sqrt{m(x)}$ when $(x,y)\in N_1$, and that $\sqrt{m(x)}\sim (1+|x|)^{-1}$, $x\in \mathbb{R}^n$, we get that
\begin{align}\label{K12infty}
\int_{m(x)}^\infty (|K_1(t,x,y)|^q+|K_2(t,x,y)|^q)\frac{dt}{t}&\leq C\int_{m(x)}^\infty\frac{dt}{t^{\frac{n}{2}q+1}}= \frac{C}{m(x)^{\frac{n}{2}q}}\nonumber\\
&\leq\frac{C}{m(x)^{\frac{q}{4}}|x-y|^{(n-\frac{1}{2})q}}\leq C\Big(\frac{\sqrt{1+|x|}}{|x-y|^{n-\frac{1}{2}}}\Big)^{q},\quad (x,y)\in N_1.
\end{align}
 
On the other hand, we can write
\begin{align*}
    K_1(t,x,y)&=-nt\left(\frac{1}{2t}\big(W_t(x-y)-T_t^\mathcal{A}(x,y)\big)+\Big(\frac{1}{2t}-\frac{1}{1-e^{-t}}\Big)T_t^\mathcal{A}(x,y)\right),\quad x,y\in \mathbb{R}^n,\;t>0,
\end{align*}
and
\begin{align*}
    K_2(t,x,y)&=\frac{|x-y|^2}{2t}\big(W_t(x-y)-T_t^\mathcal{A}(x,y)\big)+2t\Big(\frac{|x-y|^2}{(2t)^2}-\frac{|x-e^{-t}y|^2}{(1-e^{-2t})^2}\Big)T_t^\mathcal{A}(x,y)\\
    &\quad +2t\frac{|x-e^{-t}y|^2}{1-e^{-2t}}T_t^\mathcal{A}(x,y)\quad x,y\in \mathbb{R}^n,\;t>0.
\end{align*}

We have that, for each $x,y\in \mathbb{R}^n$ and $t\in (0,1)$,
$$
\big||x-y|^2-|x-e^{-t}y|^2\big|=|-(1-e^{-t})^2|y|^2+2(1-e^{-t})\langle x-y,y\rangle|\leq t^2|y|^2+2t|x-y||y|,
$$
and
\begin{align*}
\left|\frac{|x-y|^2}{(2t)^2}-\frac{|x-e^{-t}y|^2}{(1-e^{-2t})^2}\right|&= \left||x-y|^2\left(\frac{1}{(2t)^2}-\frac{1}{(1-e^{-2t})^2}\right)-|y|^2+2\frac{\langle x-y,y\rangle}{1-e^{-2t}}\right|\\
&\leq C\left(\frac{|x-y|^2}{t}+|y|^2+\frac{|x-y||y|}{t}\right).
\end{align*}
We also get that, for every $x,y\in \mathbb{R}^n$ and $t\in (0,1)$,
\begin{align*}
 \left|e^{-\frac{|x-y|^2}{2t}}-e^{-\frac{|x-e^{-t}y|^2}{1-e^{-2t}}}\right|& \leq C\left|\frac{|x-y|^2}{2t}-\frac{|x-e^{-t}y|^2}{1-e^{-2t}}\right|\exp\Big(-\min\Big\{\frac{|x-y|^2}{2t},\frac{|x-e^{-t}y|^2}{1-e^{-2t}}\Big\}\Big)\\
 & \hspace{-3cm}\leq C(|x-y|^2+t|y|^2+|x-y||y|)\exp\Big(-\min\Big\{\frac{|x-y|^2}{2t},\frac{|x-e^{-t}y|^2}{1-e^{-2t}}\Big\}\Big),
 \end{align*}
 and since 
 $$
 e^{-\frac{|x-e^{-t}y|^2}{1-e^{-2t}}}\leq Ce^{-c\frac{|x-y|^2}{t}}, \quad (x,y)\in N_1, \;t\in (0,1),
 $$
 we obtain that
 \begin{equation}\label{dife}
\left|e^{-\frac{|x-y|^2}{2t}}-e^{-\frac{|x-e^{-t}y|^2}{1-e^{-2t}}}\right| \leq C\Big(t(1+|y|^2)+|x-y||y|\Big)e^{-c\frac{|x-y|^2}{t}},\quad (x,y)\in N_1, \;t\in (0,1).
 \end{equation}
 Hence,
 \begin{align}\label{3.2.2}
     \big|(W_t(x-y)-T_t^\mathcal{A}(x,y)\big|&\leq \frac{1}{\pi ^{\frac{n}{2}}}\left(\frac{1}{(2t)^{\frac{n}{2}}}\left|e^{-\frac{|x-y|^2}{2t}}-e^{-\frac{|x-e^{-t}y|^2}{1-e^{-2t}}}\right|\right.\nonumber \\
    &\left.\quad +\left|\frac{1}{(2t)^{\frac{n}{2}}}-\frac{1}{(1-e^{-2t})^{\frac{n}{2}}}\right|e^{-\frac{|x-e^{-t}y|^2}{1-e^{-2t}}}+(1-e^{-nt})\frac{e^{-\frac{|x-e^{-t}y|^2}{1-e^{-2t}}}}{(1-e^{-2t})^{\frac{n}{2}}}\right) \nonumber \\
    &\leq C\frac{e^{-c\frac{|x-y|^2}{t}}}{t^{\frac{n}{2}-1}}\left(1+|y|^2+\frac{|x-y||y|}{t}\right),\quad (x,y)\in N_1, \;t\in (0,1). 
 \end{align}
 We have used that $|(2t)^r-(1-e^{-2t})^r|\leq Ct^{1+r}$, when $t\in (0,1)$ and $r\in \mathbb{R}$.
 
By considering all these estimations we obtain
\begin{align*}
    |K_1(t,x,y)|+|K_2(t,x,y)|&\leq C\left[\frac{e^{-c\frac{|x-y|^2}{t}}}{t^{\frac{n}{2}-1}}\left(1+|y|^2+\frac{|x-y||y|}{t}\right)\right.\\
    &\left.\quad +tT_t^\mathcal{A}(x,y)\left(1+\frac{|x-y|^2}{t}+|y|^2+\frac{|x-y||y|}{t}+\frac{|x-e^{-t}y|^2}{1-e^{-2t}}\right)\right]\\
    &\leq C\frac{e^{-c\frac{|x-y|^2}{t}}}{t^{\frac{n}{2}-1}}\left(1+|y|^2+\frac{|x-y||y|}{t}\right)\\
    &\leq C\frac{e^{-c\frac{|x-y|^2}{t}}}{t^{\frac{n}{2}-1}}\left(1+|y|^2+\frac{|y|}{\sqrt{t}}\right),\quad (x,y)\in N_1,\;t\in (0,1).
\end{align*}
By taking into account that $m(x)|y|^2\leq C$, $(x,y)\in N_1$, we can write 

\begin{align}\label{3.2.3}
t|y|^2+t^{\frac{1}{2}}|y| \leq t^{\frac{1}{4}}|y|^{\frac{1}{2}}\big(m(x)^{\frac{3}{4}}|y|^{\frac{3}{2}}+m(x)^{\frac{1}{4}}|y|^{\frac{1}{2}}\big)\leq Ct^{\frac{1}{4}}|y|^{\frac{1}{2}},\quad (x,y)\in N_1,\;t\in (0,m(x)).
\end{align}
Then, we get, for $j=1,2$,
\begin{align}\label{3.2.4}
\int_0^{m(x)}|K_j(t,x,y)|^q\frac{dt}{t}&\leq C\int_0^{m(x)}\frac{e^{-c\frac{|x-y|^2}{t}}}{t^{(\frac{n}{2}-1)q+1}}\left(1+|y|^2+\frac{|y|}{\sqrt{t}}\right)^qdt \nonumber\\
&\leq C\left(\int_0^1\frac{e^{-c\frac{|x-y|^2}{t}}}{t^{(\frac{n}{2}-1)q+1}}dt+|y|^{\frac{q}{2}}\int_0^{m(x)}\frac{e^{-c\frac{|x-y|^2}{t}}}{t^{(\frac{n}{2}-\frac{1}{4})q+1}}dt\right) \nonumber \\
&\leq C\left(\frac{1}{|x-y|^{(n-\frac{1}{2})q}}\int_0^1 t^{\frac{3}{4}q-1}dt+|y|^{\frac{q}{2}}\int_0^\infty\frac{e^{-c\frac{|x-y|^2}{t}}}{t^{(\frac{n}{2}-\frac{1}{4})q+1}}dt\right) \nonumber \\
&\leq C\frac{1+|y|^{\frac{q}{2}}}{|x-y|^{(n-\frac{1}{2})q}}\leq C\left(\frac{\sqrt{1+|x|}}{|x-y|^{n-\frac{1}{2}}}\right)^q,\quad (x,y)\in N_1.
\end{align}

We obtain
 $$\|K_1(\cdot,x,y)\|_{L^q((0,\infty),\frac{dt}{t})}+\|K_2(\cdot,x,y)\|_{L^q((0,\infty),\frac{dt}{t})}\leq C\frac{\sqrt{1+|x|}}{|x-y|^{n-1/2}},\quad (x,y)\in N_1. $$
 By putting together all the above estimates we conclude that
 \begin{equation}\label{diferencia}
  \Big\|t\partial_t(W_t(x-y)-T_t^{\mathcal A}(x,y))\Big\|_{L^q((0,\infty),\frac{dt}{t})}\leq C\frac{\sqrt{1+|x|}}{|x-y|^{n-\frac{1}{2}}},\;\;\;\;(x,y)\in N_1. 
\end{equation}
 We have that
$$
 \int_{\mathbb R^n} \frac{\sqrt{1+|x|}}{|x-y|^{n-\frac{1}{2}}}\mathcal{X}_{N_1}(x,y)dy\leq C\sqrt{1+|x|}\int_0^{\sqrt{m(x)}}\frac{d\rho}{\sqrt{\rho}}\leq Cm(x)^{\frac{1}{4}}\sqrt{1+|x|}\leq C,\quad x\in \mathbb R^n.
$$
 Since $1+|x|\sim 1+|y|$ when $(x,y)\in N_1$ we can also see that
  $$
  \sup_{y\in\mathbb R^n}\int_{\mathbb R^n} \frac{\sqrt{1+|x|}}{|x-y|^{n-\frac{1}{2}}}\mathcal{X}_{N_1}(x,y)dx<\infty.
  $$
 Then, the operator $\mathcal H_1^1$ is bounded from $L^q(\mathbb R^n,dx)$ into itself, for every $1\leq q<\infty$.
 
We conclude that $g_{\beta ,0,\{W_t\}_{t>0};\,{\rm loc}}^{q,X}$ is bounded  from $L^p_X(\mathbb R^n,dx)$ into $L^p(\mathbb R^n,dx)$.
 
The properties of invariance of the operator $g_{\beta ,0,\{W_t\}_{t>0}}^{q,X}$ allow us to prove (see \cite[p. 21]{HTV}) that  $g_{\beta ,0,\{W_t\}_{t>0}}^{q,X}$ can be extended to $L^p_X(\mathbb R^n,dx)$ as a bounded operator from $L^p_X(\mathbb R^n,dx)$ into $L^p(\mathbb R^n,dx)$. Our objective is established and (i) is proved.

 We are going to see that (iii) $\Rightarrow$ (i) and (iv) $\Rightarrow$ (i). Assume that (iii) holds. As in the proof of (ii) $\Rightarrow$ (i) we can assume that $0<\beta <1$. We can also see that
 
 $$g_{\beta ,0,\{P_t\}_{t>0};\,{\rm loc}}^{q,X}(f)\leq\left(\mathcal H(f)+g_{\beta ,0,\{P_t^{\mathcal A}\}_{t>0};\,{\rm loc}}^{q,X}(f)\right),$$
 where
  $$
  \mathcal H(f)(x)=\int_{\mathbb R^n}\|f(y)\|\left\|t\partial_t(P_t(x-y)-P_t^{\mathcal A}(x,y))\right\|_{L^q((0,\infty),\frac{dt}{t})}dy,\quad x\in\mathbb {R}^n.
  $$
 Here, $\{P_t\}_{t>0}$ denotes the Euclidean semigroup in $\mathbb R^n$.
 We consider the operator
  $$G_{\beta ,\{P_t^{\mathcal A}\}_{t>0}}^{X}(f)(t,x)=t^\beta\partial_t^\beta P_t^{\mathcal A}(f)(x),\;\;\;\;x\in\mathbb R^n\;\mbox{and}\;t>0.$$
 Since (iii) holds, $G_{\beta ,\{P_t^{\mathcal A}\}_{t>0}}^{X}$ is bounded from $L^p_X(\mathbb R^n,\gamma_{-1})$ into $L^p_{L^q_X((0,\infty),\frac{dt}{t})}(\mathbb R^n,\gamma_{-1})$.
 
 According to \cite[Lemma 4]{BCCFR} we obtain
  \begin{align*}
\| t^\beta\partial_t^\beta P_t^{\mathcal A}(x,y)\|_{L^q((0,\infty),\frac{dt}{t})} &\leq C\int_0^\infty\left(\int_0^\infty |t^\beta\partial_t^\beta (te^{-t^2/4u})|^q\frac{dt}{t}\right)^{1/q}T_u^{\mathcal A}(x,y)\frac{du}{u^{3/2}} \\
 & \leq C \int_0^\infty\left(\int_0^\infty \left(\frac{t^\beta e^{-t^2/8u}}{u^{\frac{\beta -1}{2}}}\right)^q\frac{dt}{t}\right)^{1/q}T_u^{\mathcal A}(x,y)\frac{du}{u^{3/2}} \\
 & \leq C\int_0^\infty\frac{T_u^{\mathcal A}(x,y)}{u}du,\;\;\;\;x,y\in\mathbb R^n.
 \end{align*}
Since $e^{-\frac{|x-e^{-t}y |^2}{1-e^{-2t}}}\leq Ce^{-\frac{|x-y |^2}{1-e^{-2t}}}$, $(x,y)\in N_1$, $t>0$, we get
  \begin{align*}
 \| t^\beta\partial_t^\beta P_t^{\mathcal A}(x,y)\|_{L^q((0,\infty),\frac{dt}{t})}&\leq C\left(\int_0^1\frac{e^{-\frac{|x-y|^2}{u}}}{u^{\frac{n}{2}+1}}du+\int_{m(x)}^\infty \frac{1}{u^{\frac{n}{2}+1}}du\right) \\
 & \leq C\left(\frac{1}{|x-y|^n}+\frac{1}{m(x)^{\frac{n}{2}}}\right)\leq\frac{C}{|x-y|^n} ,\quad (x,y)\in N_1.
 \end{align*}
In an analogous way we can see that, for every $i=1,...,n$,
 $$
 \| \partial_{x_i}t^\beta\partial_t^\beta P_t^{\mathcal A}(x,y)\|_{L^q((0,\infty),\frac{dt}{t})}\leq\frac{C}{|x-y|^{n+1}} ,\quad (x,y)\in N_1.
 $$
 
 From a $L^q((0,\infty),\frac{dt}{t})$-version of \cite[Proposition 3.2.7]{Sa} (see \cite[Proposition 2.3]{HTV}) we deduce that $G_{\beta ,\{P_t^{\mathcal A}\}_{t>0};\,{\rm loc}}^{X}$ is bounded from $L^p_X(\mathbb R^n,\gamma_{-1})$ into $L^p_{L^q_X((0,\infty),\frac{dt}{t})}(\mathbb R^n,\gamma_{-1})$ and from $L^p_X(\mathbb R^n,dx)$ into $L^p_{L^q_X((0,\infty),\frac{dt}{t})}(\mathbb R^n,dx)$. Hence, $g_{\beta ,0,\{P_t^\mathcal{A}\}_{t>0};\,{\rm loc}}^{q,X}$ is bounded from $L^p_X(\mathbb R^n,dx)$ into $L^p(\mathbb R^n,dx)$.
 
On the other hand, we can write
 \begin{align}\label{difPoisson}
     \big\|t\partial_t(P_t(x-y)-P_t^{\mathcal A}(x,y))\big\|_{L^q((0,\infty),\frac{dt}{t})}&\nonumber\\
     &\hspace{-5cm}=\frac{1}{\sqrt{\pi}}\left\|\int_0^\infty e^{-u} t\partial _t\big(W_{t^2/4u}(x-y)-T_{t^2/4u}^\mathcal{A}(x,y)\big)\frac{du}{\sqrt{u}}\right\|_{L^q((0,\infty),\frac{dt}{t})}\nonumber\\
     &\hspace{-5cm}\leq C\int_0^\infty e^{-u}\left\|t\partial _t\big(W_{t^2/4u}(x-y)-T_{t^2/4u}^\mathcal{A}(x,y)\big)\right\|_{L^q((0,\infty),\frac{dt}{t})}\frac{du}{\sqrt{u}}\nonumber\\
     &\hspace{-5cm}\leq C\left\|s\partial _s\big(W_s(x-y)-T_s^\mathcal{A}(x,y)\big)\right\|_{L^q((0,\infty),\frac{ds}{s})}\int_0^\infty e^{-u}\frac{du}{\sqrt{u}}\leq C\frac{\sqrt{1+|x|}}{|x-y|^{n-\frac{1}{2}}},\quad (x,y)\in N_1.
 \end{align}
In the last inequality we have taken into account the estimation \eqref{diferencia}. 

\noindent We deduce that the operator $\mathcal H$ is bounded from $L^q_X(\mathbb R^n,dx)$ into $L^q(\mathbb R^n,dx)$, for every $1\leq q\leq \infty$.
 
 It is concluded that $g_{\beta ,0,\{P_t\}_{t>0};\,{\rm loc}}^{q,X}$ is bounded  from $L^p_X(\mathbb R^n,dx)$ into $L^p(\mathbb R^n,dx)$. As above, the properties of invariance of the operator $g_{\beta ,0,\{P_t\}_{t>0}}^{q,X}$ allow us to deduce that  $g_{\beta ,0,\{P_t\}_{t>0}}^{q,X}$ defines a bounded operator from $L^p_X(\mathbb R^n,dx)$ into $L^p(\mathbb R^n,dx)$. By \cite[Theorem C]{TZ} (i) holds.
 
If (iv) holds the same arguments lead to that (i) is true (see \cite[Theorem 5.2]{TZ}).

Suppose now that (iii) holds. We are going to prove that (iv) is true.

\noindent As we have just proved we deduce that $g_{\beta ,0,\{P_t\}_{t>0}}^{q,X}$ is bounded  from $L^p_X(\mathbb R^n,dx)$ into $L^p(\mathbb R^n,dx)$. According to \cite[Theorem 5.2]{TZ}, $g_{\beta ,0,\{P_t\}_{t>0}}^{q,X}$ is bounded  from $L^1_X(\mathbb R^n,dx)$ into $L^{1,\infty}(\mathbb R^n,dx)$.

We can write
$$g_{\beta ,0,\{P_t^{\mathcal A}\}_{t>0}}^{q,X}(f)\leq \mathcal H(f)+g_{\beta ,0,\{P_t\}_{t>0};\,{\rm loc}}^{q,X}(f)+g_{\beta ,0,\{P_t^{\mathcal A}\}_{t>0};\,{\rm glob}}^{q,X}(f).$$
According to \cite[Propositions 3.2.5 and 3.2.7]{Sa} in the $L^q((0,\infty),\frac{dt}{t})$-setting, $g_{\beta ,0,\{P_t\}_{t>0};\,{\rm loc}}^{q,X}$ is bounded  from $L^1_X(\mathbb R^n,\gamma_{-1})$ into $L^{1,\infty}(\mathbb R^n,\gamma_{-1})$. We proved that $\mathcal H$ is a bounded operator from $L^1_X(\mathbb R^n,dx)$ into $L^1(\mathbb R^n,dx)$. Since $\mathcal H$ is local, \cite[Proposition 3.2.5]{Sa} says us that $\mathcal H$ is bounded from $L^1_X(\mathbb R^n,\gamma_{-1})$ into $L^1(\mathbb R^n,\gamma_{-1})$. By taking into account (\ref{operatorL}) we have that
$$
g_{1 ,0,\{P_t^{\mathcal A}\}_{t>0};\,{\rm glob}}^{q,X}(f)(x)\leq \int_{\mathbb R^n}\|f(y)\|L(x,y)dy,\;\;\;\;f\in L_X^1(\mathbb R^n,\gamma_{-1}),
$$
where $L\geq 0$ and being the operator
$$
\mathcal L(f)(x)=\int_{\mathbb R^n}\|f(y)\|L(x,y)dy,\;\;\;\;f\in L_X^1(\mathbb R^n,\gamma_{-1}),
$$
bounded from $L^1_X(\mathbb R^n,\gamma_{-1})$ into $L^{1,\infty}(\mathbb R^n,\gamma_{-1})$. Since $0<\beta<1$, $g_{\beta ,0,\{P_t^{\mathcal A}\}_{t>0};\,{\rm glob}}^{q,X}$ is bounded from $L^1_X(\mathbb R^n,\gamma_{-1})$ into $L^{1,\infty}(\mathbb R^n,\gamma_{-1})$.

We conclude that $g_{\beta ,0,\{P_t^{\mathcal A}\}_{t>0}}^{q,X}$ is bounded from $L^1_X(\mathbb R^n,\gamma_{-1})$ into $L^{1,\infty}(\mathbb R^n,\gamma_{-1})$ and (iv) holds.

The proof is complete.

\subsection{Proof of Theorem \ref{Th1.3}}
 $\{T_t^{\mathcal A}\}_{t>0}$ is a symmetric diffusion semigroup. According to \cite[Theorem 2, (ii)]{Xu2} and \cite[Theorem 1.7]{BFGM} we deduce that the property (i) $\Rightarrow$ (ii) holds. From \cite[Theorem B]{TZ} (see also \cite[Theorem 2.2]{MTX}) it deduces that the property (i) $\Rightarrow$ (iii) is true.
 
 We are going to see that (ii) $\Rightarrow$ (i). Assume that (ii) is true. According to \cite[Proposition 3.1]{TZ} we have that, for every $\delta\geq\beta$,
 
 $$\|f\|_{L^p_X(\mathbb R^n,\gamma_{-1})}\leq C\|g_{\delta ,0,\{T_t^{\mathcal A}\}_{t>0}}^{q,X}(f)\|_{L^p(\mathbb R^n,\gamma_{-1})},\;\;\;\;f\in L^p_X(\mathbb R^n,\gamma_{-1}).$$
 We can assume that $\beta\in\mathbb N\setminus\{0\}$.
 
 Since $X$ is a Banach space, $X$ is a closed subspace of $(X^*)^*$. Then according to \cite[Proposition 1.3.1]{HNVW} $L^p_X(\mathbb R^n,\gamma_{-1})$ is norming for $L^{p'}_{X^*}(\mathbb R^n,\gamma_{-1})$, where $p'=\frac{p}{p-1}$. This means that, for every $g\in L^{p'}_{X^*}(\mathbb R^n,\gamma_{-1})$
 
 $$\|g\|_{L^{p'}_{X^*}(\mathbb R^n,\gamma_{-1})}=\sup_{\begin{array}{c} f\in L^{p}_{X}(\mathbb R^n,\gamma_{-1}) \\ \|f\|_{L^{p}_{X}(\mathbb R^n,\gamma_{-1})} \leq 1\end{array}}\left|\int_{\mathbb R^n}\langle f(x),g(x)\rangle d\gamma_{-1}(x)\right|.$$
 We now define the operator
 
 $$Q_\beta^{\mathcal A}(h)(x)=\int_0^\infty s^{\beta -1}(\partial_s^\beta T_s^{\mathcal A})(h(s,\cdot))(x)ds,\;\;\;\;x\in\mathbb R^n.$$
 Our objective is to prove that the operator $\mathcal L_{q,\beta}^{\mathcal A}$ defined by
 
 $$\mathcal L_{q,\beta}^{\mathcal A}(h)=g_{\beta ,0,\{T_t^{\mathcal A}\}_{t>0}}^{q,X}(Q_\beta^{\mathcal A}(h))$$
 is bounded from $L^r_{L^q_X((0,\infty),\frac{dt}{t})}(\mathbb R^n,\gamma_{-1})$ into $L^r(\mathbb R^n,\gamma_{-1})$, $1<r<\infty$.
 
When this objective is established we can see that (i) holds following the procedure developed in the proof of \cite[Theorem 3.1]{MTX}. Indeed, let $f\in C_c^\infty(\mathbb R^n)\otimes X^*$. According to Theorem \ref{Th1.1}, $g_{\beta ,0,\{T_t^{\mathcal A}\}_{t>0}}^{q',X^*}(f)\in L^{p}(\mathbb R^n,\gamma_{-1})$, or in other words, $t^\beta\partial_t^\beta T_t^{\mathcal A}(f)\in L^{p}_{L^{q'}_{X^*}((0,\infty),\frac{dt}{t})}(\mathbb R^n,\gamma_{-1})\subseteq L^{p}_{(L^q_X((0,\infty),\frac{dt}{t}))^*}(\mathbb R^n,\gamma_{-1})$. Then, without loss of generality we can suppose that there exists $h\in L^{p'}_{L^q_X((0,\infty),\frac{dt}{t})}(\mathbb R^n,\gamma_{-1})$, that can be chosen smooth enough, such that $\|h\|_{L^{p'}_{L^q_X((0,\infty),\frac{dt}{t})}(\mathbb R^n,\gamma_{-1})}\leq 1$ and
 
\begin{align*}
 \|g_{\beta ,0,\{T_t^{\mathcal A}\}_{t>0}}^{q',X^*}(f)\|_{L^{p}(\mathbb R^n,\gamma_{-1})}&=\|t^\beta\partial_t^\beta T_t^{\mathcal A}(f)\|_{L^{p}_{L^{q'}_{X^*}((0,\infty),\frac{dt}{t})}(\mathbb R^n,\gamma_{-1})} \\
 & = \int_{\mathbb R^n}\int_0^\infty\langle t^\beta\partial_t^\beta T_t^{\mathcal A}(f)(x),h(t,x)\rangle_{X^*,X}\frac{dt}{t}d\gamma_{-1}(x).
 \end{align*}
 By interchanging the order of integration that is justified by the smoothness of $f$ and $h$ we obtain
 
 \begin{align*}
 \|g_{\beta ,0,\{T_t^{\mathcal A}\}_{t>0}}^{q',X^*}(f)\|_{L^{p}(\mathbb R^n,\gamma_{-1})}&= \int_{\mathbb R^n}\int_0^\infty\langle t^\beta\partial_t^\beta T_t^{\mathcal A}(f)(x),h(t,x)\rangle_{X^*,X}\frac{dt}{t}d\gamma_{-1}(x) \\
& =\int_{\mathbb R^n} \langle f(x),Q_\beta^{\mathcal A}(h)(x)\rangle_{X^*,X}d\gamma_{-1}(x).
 \end{align*}
 By using H\"older's inequality, and our objective, that is assumed proved, we get
   \begin{align*}
 \|g_{\beta ,0,\{T_t^{\mathcal A}\}_{t>0}}^{q',X^*}(f)\|_{L^{p}(\mathbb R^n,\gamma_{-1})}
 &\leq \|f\|_{L^{p}_{X^*}(\mathbb R^n,\gamma_{-1})} \|Q_\beta^{\mathcal A}(h)\|_{L^{p'}_{X}(\mathbb R^n,\gamma_{-1})}\\
 &\hspace{-3cm}\leq C \|f\|_{L^{p}_{X^*}(\mathbb R^n,\gamma_{-1})} \|\mathcal L_{q,\beta}^{\mathcal A}(h)\|_{L^{p'}(\mathbb R^n,\gamma_{-1})}\leq C \|f\|_{L^{p}_{X^*}(\mathbb R^n,\gamma_{-1})} \|h\|_{L^{p'}_{L^q_X((0,\infty),\frac{dt}{t})}(\mathbb R^n,\gamma_{-1})} \\
& \hspace{-3cm}\leq C \|f\|_{L^{p}_{X^*}(\mathbb R^n,\gamma_{-1})}.
 \end{align*}
 By Theorem \ref{Th1.2}, there exists a norm $\VERT\cdot\VERT_{X^*}$ in $X^*$ defining the original topology of $X^*$ and such that $(X^*,\VERT\cdot\VERT_{X^*})$ is $q'$-uniformly convex. This is equivalent to that $X^*$ has $q'$ martingale cotype. Then, $X$ has $q$-martingale type (see \cite[Theorem 3.1]{MTX}). Hence, there exists a norm $\VERT\cdot\VERT_{X}$ in $X$ that defines the topology of $X$ being $(X,\VERT\cdot\VERT_{X})$ $q$-uniformly smooth.
 
 We are going to prove our objective. Suppose that $h\in C_c^\infty(\mathbb R^n)\otimes C_c(0,\infty)\otimes X$. Note that $C_c^\infty(\mathbb R^n)\otimes C_c(0,\infty)\otimes X$ is a dense subspace of $L^r_{L^q_X((0,\infty),\frac{dt}{t})}(\mathbb R^n,\gamma_{-1})$. We can write
 
 \begin{align*}
 \partial_t^\beta T_t^{\mathcal A}(Q_\beta^{\mathcal A}(h))(x)&=\int_{\mathbb R^n}\partial_t^\beta T_t^{\mathcal A}(x,y)Q_\beta^{\mathcal A}(h)(y)dy \\
& = \int_{\mathbb R^n}\partial_t^\beta T_t^{\mathcal A}(x,y)\int_0^\infty s^{\beta -1}\int_{\mathbb R^n} \partial_s^\beta T_s^{\mathcal A}(y,z)h(s,z)dzdsdy \\
& = \int_0^\infty\int_{\mathbb R^n} h(s,z)\int_{\mathbb R^n} s^{\beta -1} \partial_s^\beta T_s^{\mathcal A}(y,z)\partial_t^\beta T_t^{\mathcal A}(x,y) dydzds,\quad x\in\mathbb R^n\;\mbox{and}\;t>0.
 \end{align*}
 By using the semigroup property we get
  \begin{align*}
 \int_{\mathbb R^n}  \partial_s^\beta T_s^{\mathcal A}(y,z)\partial_t^\beta T_t^{\mathcal A}(x,y) dy&=\partial_t^\beta \partial_s^\beta\int_{\mathbb R^n}  T_s^{\mathcal A}(y,z) T_t^{\mathcal A}(x,y) dy \\
& =\partial_t^\beta \partial_s^\beta T_{t+s}^{\mathcal A}(x,z)=\partial_u^{2\beta}T_u^{\mathcal A}(x,z)_{\big|u=t+s},\quad x,z\in\mathbb R^n\;\mbox{and}\;t,s\in (0,\infty).
 \end{align*}
 We consider the operator $L_\beta^{\mathcal A}$ defined by
  $$
  L_\beta^{\mathcal A}(h)(t,x)=\int_0^\infty\int_{\mathbb R^n}h(s,z)K_\beta^{\mathcal A}(t,x;s,z)dz\frac{ds}{s},\;\;\;x\in\mathbb R^n\;\mbox{and}\;t>0,
  $$
 where
  $$
  K_\beta^{\mathcal A}(t,x;s,z)=(st)^\beta\partial_u^{2\beta} T_u^{\mathcal A}(x,z)_{\big|u=t+s},\;\;\;x,z\in\mathbb R^n\;\mbox{and}\;t,s>0.$$
 Note that
  $$
  \mathcal L_{q,\beta}^{\mathcal A}(h)(x)=\|L_\beta^{\mathcal A}(h)(\cdot,x)\|_{L^q_X((0,\infty),\frac{dt}{t})},\;\;\;\;x\in\mathbb R^n.$$
 We are going to see that $ L_{\beta}^{\mathcal A}$ is bounded from $L^r_{L^q_X((0,\infty),\frac{dt}{t})}(\mathbb R^n,\gamma_{-1})$ into itself, for every $1<r<\infty$. We consider $L_{\beta ,{\rm loc}}^{\mathcal A}$ and $L_{\beta, {\rm glob}}^{\mathcal A}$ defined in the usual way. 
 
We first observe that, since $q'/q\geq 1$ and $\displaystyle\int_0^\infty\frac{(st)^{\beta q}}{(s+t)^{2\beta q}}\frac{dt}{t}=\int_0^\infty\frac{z^{\beta q-1}}{(1+z)^{2\beta q}}dz=B(\beta q,\beta q)$, $s>0$, by using Jensen's inequality we obtain that, if $F$ is a function defined on $(0,\infty)$ such that $u^{2\beta}F\in L^{q'}((0,\infty ),\frac{du}{u})$, then
 \begin{align}\label{Jensen}
     \Big\|(st)^\beta F(t+s)\Big\|_{L^{q'}_{L^q((0,\infty),\frac{dt}{t})}((0,\infty),\frac{ds}{s})}&= C\left(\int_0^\infty \left(\int_0^\infty (st)^{\beta q}|F(t+s)|^q\frac{dt}{t}\right)^{q'/q}\frac{ds}{s}\right)^{1/q'}\nonumber\\
     &\hspace{-3cm}\leq C\left(\int_0^\infty \int_0^\infty\frac{(st)^{\beta q-1}}{(s+t)^{2\beta q}}|(t+s)^{2\beta} F(t+s)|^{q'}dtds\right)^{1/q'}\nonumber\\
     &\hspace{-3cm} =C\left(\int_0^\infty u^{2\beta (q'-q)}|F(u)|^{q'}\int_0^u((u-t)t)^{\beta q-1}dtdu\right)^{1/q'}=C\|u^{2\beta} F\|_{L^{q'}((0,\infty ),\frac{du}{u})}.
 \end{align}
 
 We now study the operator $L_{\beta, {\rm glob}}^{\mathcal A}$.
 
According to \eqref{AcotDeriv} we have that, for $0<\delta<\eta<1$,
 \begin{equation}\label{acotacionderivbeta}
 |u^{2\beta}\partial_u^{2\beta}T_u^\mathcal A(x,z)|\leq C\frac{e^{-\frac{u}{2}}}{(1-e^{-2u})^{\frac{n}{2}}}e^{(\eta-\delta)\frac{|z|^2-|x|^2}{2}}e^{-\delta\frac{ |z-e^{-u}x|^2}{1-e^{-2u}}},\quad x,z\in\mathbb R^n\;\mbox{and}\;u>0. 
 \end{equation}

By using Minkowski and H\"older's inequalities we get 
\begin{align*}
  \|L_{\beta, {\rm glob}}^{\mathcal A}(h)(\cdot,x)\|_{L^q_X((0,\infty),\frac{dt}{t})}&\leq\int_{\mathbb R^n}\int_0^\infty \|K_{\beta}^{\mathcal A}(\cdot,x;s,z)\|_{L^q((0,\infty),\frac{dt}{t})}\|h(s,z)\|_X\mathcal{X}_{N_1^c}(x,z)\frac{ds}{s}dz \\
  & \hspace{-4cm}\leq \int_{\mathbb R^n}\Big\|K_{\beta}^{\mathcal A}(\cdot,x;\cdot,z)\Big\|_{L^{q'}_{L^q((0,\infty),\frac{dt}{t})}((0,\infty),\frac{ds}{s})}\|h(\cdot,z)\|_{L^q_X((0,\infty),\frac{ds}{s})}\mathcal{X}_{N_1^c}(x,z)dz, \quad x\in \mathbb{R}^n.
\end{align*}
By considering \eqref{Jensen} and \eqref{acotacionderivbeta} it follows that
\begin{align*}
\Big\|K_{\beta}^{\mathcal A}(\cdot,x;\cdot,z)\Big\|_{L^{q'}_{L^q((0,\infty),\frac{dt}{t})}((0,\infty),\frac{ds}{s})}&= \Big\|(st)^{\beta}\partial_u^{2\beta}T_u^{\mathcal A}(x,z)_{|u=t+s}|\Big\|_{L^{q'}_{L^q((0,\infty),\frac{dt}{t})}((0,\infty),\frac{ds}{s})}\\
&\hspace{-4cm}\leq C\|u^{2\beta} \partial_u^{2\beta}T_u^{\mathcal A}(x,z)\|_{L^{q'}((0,\infty ),\frac{du}{u})}
    \leq Ce^{(\eta-\delta)\frac{|z|^2-|x|^2}{2}}\Big\|\frac{e^{-\frac{u}{2}}e^{-\delta \frac{|z-e^{-u}x|^2}{1-e^{-2u}}}}{(1-e^{-2u})^{\frac{n}{2}}}\Big\|_{L^{q'}((0,\infty),\frac{du}{u})}\\
    &\hspace{-4cm}\leq Ce^{(\eta-\delta)\frac{|z|^2-|x|^2}{2}}\Big\|\frac{e^{-\frac{u}{2}}e^{-\delta \frac{|z-e^{-u}x|^2}{1-e^{-2u}}}}{(1-e^{-2u})^{\frac{n}{2}}}\Big\|_{L^{q'}((0,\infty),\frac{du}{u})} ,\quad x,z\in \mathbb{R}^n.
\end{align*}
Then,
$$
  \|L_{\beta, {\rm glob}}^{\mathcal A}(h)(\cdot,x)\|_{L^q_X((0,\infty),\frac{dt}{t})}\leq C\int_{\mathbb{R}^n}H(x,z)\|h(\cdot,z)\|_{L^q_X((0,\infty),\frac{ds}{s})}\mathcal{X}_{N_1^c}(x,z)dz, \quad x\in \mathbb{R}^n,
$$
where 
$$
H(x,z)=e^{(\eta-\delta)\frac{|z|^2-|x|^2}{2}}\Big\|\frac{e^{-\frac{u}{2}}e^{-\delta \frac{|z-e^{-u}x|^2}{1-e^{-2u}}}}{(1-e^{-2u})^{\frac{n}{2}}}\Big\|_{L^{q'}((0,\infty),\frac{du}{u})},\quad x,z\in \mathbb{R}^n.
$$
By proceeding as in the estimation \eqref{A2} we get that,
$$
H(x,z)\leq C \left\{\begin{array}{ll}
e^{-(\eta+\delta)\frac{|x|^2}{2} +(\eta-\delta)\frac{|z|^2}{2}},&\langle x,z \rangle \leq  0, \\[0.2cm]
|x+z|^ne^{\frac{\eta}{2}(|z|^2-|x|^2)-\frac{\delta}{2}|x+z||x-z|},&\langle x,z \rangle > 0. \end{array}\right., \quad (x,z)\in N_1^c.
$$

\noindent Let $1<r<\infty$. By choosing $0<\eta<\delta <1$ such that $\eta -\delta<\frac{2}{r}<\eta +\delta$, we deduce that
$$
\sup_{x\in\mathbb R^n}\int_{\mathbb R^n}\mathcal{X}_{N_1^c}(x,z)H(x,z)e^{\frac{|x|^2-|z|^2}{r}}dz<\infty,
$$
and
$$
\sup_{z\in\mathbb R^n}\int_{\mathbb R^n}\mathcal{X}_{N_1^c}(x,z)H(x,z)e^{\frac{|x|^2-|z|^2}{r}}dx<\infty.
$$
We conclude that the operator defined by
$$
\mathbb H(g)(x)=\int_{\mathbb R^n}H(x,z)g(z)dz,\quad x\in \mathbb{R}^n,
$$
is bounded from $L^r(\mathbb R^n,\gamma_{-1})$ into itself. Hence, $L_{\beta, {\rm glob}}^{\mathcal A}$ is bounded from $L^r_{L^q_X((0,\infty),\frac{dt}{t})}(\mathbb R^n,\gamma_{-1})$ into itself.

We now study the local operator $L_{\beta, {\rm loc}}^{\mathcal A}$. We consider the operator $L_\beta$ defined by
$$L_\beta(h)(t,x)=\int_0^\infty\int_{\mathbb R^n}h(s,z)K_\beta(t,x;s,z)dz\frac{ds}{s},\quad x\in\mathbb R^n,\;t>0,$$
where $K_\beta(t,x;s,z)=(st)^\beta\partial_u^{2\beta}W_u(x-z)_{|u=t+s}$, $x,z\in\mathbb R^n$ and $s,t>0$.

We have that, for every $\ell\in \mathbb{N}$,
 \begin{equation}\label{3.2.1}
\partial_u^\ell W_u(z)=\frac{1}{2^\ell\sqrt{\pi}}\frac{1}{(2u)^{\frac{1}{2}+\ell}}\widetilde{H}_{2\ell}\left(\frac{z}{\sqrt{2u}}\right)=\frac{1}{(4u)^\ell}H_{2\ell}\Big(\frac{z}{\sqrt{2u}}\Big)W_u(z),\quad z\in\mathbb R,\;u>0.
\end{equation}
We prove this equality by induction on $\ell$. Since $\widetilde{H}'_n(z)=-\widetilde{H}_{n+1}(z)$ and $H_{n+1}(z)=2zH_n(z)-2nH_{n-1}(z)$, $z\in \mathbb{R}$ and $n\in \mathbb{N}$, we have that
\begin{align*}
 \partial_u W_u(z)&=\partial _u\left(\frac{1}{\sqrt{2\pi u}}\widetilde{H}_0\Big(\frac{z}{\sqrt{2u}}\Big)\right)=\frac{1}{\sqrt{\pi}(2u)^{\frac{3}{2}}}\left(-\widetilde{H}_0\Big(\frac{z}{\sqrt{2u}}\Big)+\frac{z}{\sqrt{2u}}\widetilde{H}_1\Big(\frac{z}{\sqrt{2u}}\Big)\right)\\
 &=\frac{1}{2\sqrt{\pi}}\frac{1}{(2u)^{3/2}}\widetilde{H}_2\Big(\frac{z}{\sqrt{2u}}\Big),\quad z\in\mathbb R,\;u>0.
\end{align*}
Suppose that \eqref{3.2.1} is true for certain $\ell\in\mathbb N$. Then,
\begin{align*}
\partial_u^{\ell+1} W_u(z) & =\frac{1}{2^{2\ell +\frac{1}{2}}\sqrt{\pi}}\left(-\Big(\frac{1}{2}+\ell\Big)\frac{1}{u^{\frac{3}{2}+\ell}}\widetilde{H}_{2\ell}\Big(\frac{z}{\sqrt{2u}}\Big)+\frac{z}{2\sqrt{2}u^{2+\ell}}\widetilde{H}_{2\ell +1}\Big(\frac{z}{\sqrt{2u}}\Big)\right) \\
& =\frac{1}{2^\ell\sqrt{\pi}}\frac{1}{(2u)^{\frac{3}{2}+\ell}}\left(-(2\ell +1)\widetilde{H}_{2\ell}\Big(\frac{z}{\sqrt{2u}}\Big)+\frac{z}{\sqrt{2u}}\widetilde{H}_{2\ell +1}\Big(\frac{z}{\sqrt{2u}}\Big)\right)\\
&=\frac{1}{2^{\ell +1}\sqrt{\pi}}\frac{1}{(2u)^{\frac{3}{2}+\ell}}\widetilde{H}_{2\ell +2}\Big(\frac{z}{\sqrt{2u}}\Big),\quad z\in\mathbb R,\;u>0.
\end{align*}
From \eqref{3.2.1} it follows that, for every $m\in \mathbb{N}$ and $x=(x_1,...,x_n)\in \mathbb{R}^n$,
\begin{align*}
\partial_u^mW_u(x)&=\partial _u\Big(\prod_{i=1}^n W_u(x_i)\Big)=W_u(x)\sum_{|r|=m} {m\choose {r_1\ldots r_n}}\prod_{i=1}^n\frac{1}{(4u)^{r_i}}H_{2r_i}\left(\frac{x_i}{\sqrt{2u}}\right)\\
&=\frac{W_u(x)}{(4u)^m}\sum_{|r|=m} {m\choose {r_1\ldots r_n}}H_{2r}\left(\frac{x}{\sqrt{2u}}\right).
\end{align*}
We define $\mathcal M_{\beta,{\rm loc}}=L_{\beta,{\rm loc}}^{\mathcal A}-L_{\beta,{\rm loc}}$. Let us see that $\mathcal M_{\beta,{\rm loc}}$ is bounded from $L^r_{L^q_X((0,\infty),\frac{dt}{t})}(\mathbb R^n,\gamma_{-1})$ into itself, for every $1\leq r<\infty$.

By Minkowski integral and H\"older inequalities and \eqref{Jensen} we get, for every $x\in \mathbb{R}^n$,
\begin{align*}
  \|\mathcal{M}_{\beta, {\rm loc}}^{\mathcal A}(h)(\cdot,x)\|_{L^q_X((0,\infty),\frac{dt}{t})}&\\
  &\hspace{-3cm}\leq \int_{\mathbb R^n}\Big\|K_{\beta}^{\mathcal A}(\cdot,x;\cdot,z)-K_\beta(\cdot,x;\cdot,z)\Big\|_{L^{q'}_{L^q((0,\infty),\frac{dt}{t})}((0,\infty),\frac{ds}{s})}\|h(\cdot,z)\|_{L^q_X((0,\infty),\frac{ds}{s})}\mathcal{X}_{N_1}(x,z)dz\\
  &\hspace{-3cm}\leq C\int_{\mathbb R^n}\Big\|u^{2\beta}\partial_u^{2\beta}(T_u^{\mathcal A}(x,z)-W_u(x-z))\Big\|_{L^{q'}((0,\infty),\frac{du}{u})}\|h(\cdot,z)\|_{L^q_X((0,\infty),\frac{ds}{s})}\mathcal{X}_{N_1}(x,z)dz. 
\end{align*}

Our objective is to established that
\begin{equation}\label{objective}
\Big\|u^{2\beta}\partial_u^{2\beta}(T_u^{\mathcal A}(x,z)-W_u(x-z))\Big\|_{L^{q'}((0,\infty),\frac{du}{u})}\leq C\frac{\sqrt{1+|x|}}{|x-z|^{n-1/2}},\quad (x,z)\in N_1.
\end{equation}
Then, since $\mathcal M_{\beta,{\rm loc}}$ is a local operator, by virtue of \cite[Proposition 3.2.5]{Sa} we can conclude that  $\mathcal M_{\beta,{\rm loc}}$ is bounded from $L^r_{L^q_X((0,\infty),\frac{dt}{t})}(\mathbb R^n,\gamma_{-1})$ into itself, for every $1\leq r<\infty$. Recall that the operator
$$
\mathbb T(g)(x)=\int_{\mathbb R^n}\frac{\sqrt{1+|x|}}{|x-z|^{n-1/2}}\mathcal{X}_{N_1}(x,y)g(y)dy,\quad x\in \mathbb{R}^n,
$$
is bounded from $L^r(\mathbb R^n,dx)$ into itself, for every $1\leq r<\infty$ (see the proof of Theorem \ref{Th1.2}).

Let us show \eqref{objective}. By considering \eqref{cuentaderiv} and \eqref{3.2.1} we have that
\begin{align*}
    \partial_u^{2\beta}(T_u^{\mathcal A}(x,z)-W_u(x-z))&\\
    &\hspace{-2cm}=T_u^\mathcal{A}(x,z)\sum_{(j,r,s,\ell)\in \mathbb{J}} c_{j,r,s,\ell}\left(\frac{e^{-u}}{\sqrt{1-e^{-2u}}}\right)^{2|s|-|\ell|}H_{\ell}(z)H_{2s-\ell}\left(\frac{x-e^{-u}z}{\sqrt{1-e^{-2u}}}\right)\\
    &\hspace{-2cm}\quad -\frac{W_u(x-z)}{(4u)^{2\beta}}\sum_{|r|= 2\beta}{2\beta\choose {r_1\ldots r_n}}H_{2r}\left(\frac{x-z}{\sqrt{2u}}\right),\quad x,z\in \mathbb{R}^n,\;u>0.
\end{align*}
Here we understand that $(j,r,s,\ell)\in \mathbb{J}$ when $j=0,...,2\beta$, and $r,s,\ell\in \mathbb{N}^n$ satisfy $|r|=j$ and $\ell_i\leq s_i\leq r_i$, $i=1,...,n$. For every $(j,r,s,\ell)\in \mathbb{J}$ the constant $c_{j,r,s,\ell}$ is given by
$$
c_{j,r,s,\ell}=\frac{(-1)^{|s|+|\ell|}}{2^{|s|}}{2\beta \choose j}n^{2\beta-j}{j\choose r_1\ldots r_n}\prod_{i=1}^n\stirling{r_i}{s_i}{s_i \choose \ell_i}.
$$
Denoting by $\mathbb{K}=\mathbb{J}\setminus \{(2\beta,r,r,\ell): r,\ell \in \mathbb{N}^n, |r|=2\beta, \ell =(0,...,0)\}$, we can write 
\begin{align*}
    \partial_u^{2\beta}(T_u^{\mathcal A}(x,z)-W_u(x-z))&\\
    &\hspace{-4cm}=T_u^\mathcal{A}(x,z)\sum_{(j,r,s,\ell)\in \mathbb{K}} c_{j,r,s,\ell}\left(\frac{e^{-u}}{\sqrt{1-e^{-2u}}}\right)^{2|s|-|\ell|}H_{\ell}(z)H_{2s-\ell}\left(\frac{x-e^{-u}z}{\sqrt{1-e^{-2u}}}\right)\\
    &\hspace{-4cm}\quad +\left(T_u^\mathcal{A}(x,z)\frac{e^{-4\beta u}}{2^{2\beta}(1-e^{-2u})^{2\beta}}\sum_{|r|= 2\beta}{2\beta\choose {r_1\ldots r_n}}H_{2r}\left(\frac{x-e^{-u}z}{\sqrt{1-e^{-2u}}}\right)\right.\\
    &\hspace{-4cm}\quad 
    \left.
     -\frac{W_u(x-z)}{(4u)^{2\beta}}\sum_{|r|= 2\beta}{2\beta\choose {r_1\ldots r_n}}H_{2r}\left(\frac{x-z}{\sqrt{2u}}\right)\right)\\
     &\hspace{-4cm}=I_1(u,x,z)+I_2(u,x,z),\quad x,z\in \mathbb{R}^n,\;u>0.
\end{align*}
We have that
$$
|I_1(u,x,z)|\leq Ce^{-nu}\sum_{(j,r,s,\ell)\in \mathbb{K}}\frac{(1+|z|)^{|\ell|}e^{-c\frac{|x-e^{-u}z|^2}{1-e^{-2u}}}}{(1-e^{-2u})^{\frac{n}{2}+|s|-\frac{|\ell|}{2}}},\quad x,z\in \mathbb{R}^n,\;u>0.
$$
Then, by making the change of variables  $u=\log\frac{1+v}{1-v}$, $v\in (0,\infty)$, we get
\begin{align*}
\Big\|u^{2\beta}I_1(u,x,z)\Big\|_{L^{q'}((0,\infty),\frac{du}{u})}&\leq C\sum_{(j,r,s,\ell)\in \mathbb{K}}\left(\int_0^\infty \Big(\frac{u^{2\beta}e^{-nu}(1+|z|)^{|\ell|}e^{-c\frac{|x-e^{-u}z|^2}{1-e^{-2u}}}}{(1-e^{-2u})^{\frac{n}{2}+|s|-\frac{|\ell|}{2}}}\Big)^{q'}\frac{du}{u}\right)^{1/q'}\\
&\hspace{-4cm}\leq  Ce^{-c(|x|^2-|z|^2)}\sum_{(j,r,s,\ell)\in \mathbb{K}}\left(\int_0^1 \Big(\log\frac{1+v}{1-v}\Big)^{2\beta q'-1}\frac{(1-v)^{nq'-1}(1+|z|)^{q'|\ell|}e^{-c(v|x+z|^2+\frac{|x-z|^2}{v})}}{v^{(\frac{n}{2}+|s|-\frac{|\ell|}{2})q'}}dv\right)^{1/q'}\\
&\hspace{-4cm}\leq  C\sum_{(j,r,s,\ell)\in \mathbb{K}}\left(\int_0^1 \frac{(1+|z|)^{q'|\ell|}e^{-c(v|x+z|^2+\frac{|x-z|^2}{v})}}{v^{(\frac{n}{2}+|s|-\frac{|\ell|}{2}-2\beta)q'+1}}dv\right)^{1/q'},\quad (x,z)\in N_1.
\end{align*}
In the last inequality we have taken into account that $||x|^2-|z|^2|\leq C$, $(x,z)\in N_1$, that $\log\frac{1+v}{1-v}\leq Cv$, $v\in (0,\frac{1}{2})$, and that $(-\log(1-v))^{2\beta q'-1}(1-v)^{nq'-1}\leq C$, $v\in (\frac{1}{2},1)$.

Let $(j,r,s,\ell)\in \mathbb{K}$. In the case that $j\leq 2\beta -1$ or $j=2\beta$ and $s\not =r$ we have that $|s|\leq 2\beta -1$, and we can write
\begin{align*}
\frac{(1+|z|)^{q'|\ell|} e^{-c(v|x+z|^2+\frac{|x-z|^2}{v})}}{v^{(\frac{n}{2}+|s|-\frac{|\ell|}{2}-2\beta)q'+1}}&\leq C\frac{(1+|x+z|^{q'|\ell|}+|x-z|^{q'|\ell|}) e^{-c(v|x+z|^2+\frac{|x-z|^2}{v})}}{v^{(\frac{n}{2}-1-\frac{|\ell|}{2})q'+1}}\\
&\leq C \frac{e^{-c\frac{|x-z|^2}{v}}}{v^{(\frac{n}{2}-1)q'+1}}\leq C \frac{e^{-c\frac{|x-z|^2}{v}}}{v^{(\frac{n}{2}-\frac{1}{4})q'+1}},\quad x,z\in \mathbb{R}^n,\;v\in (0,1).
\end{align*}
If $j=2\beta$ and $s=r$, then $|s|=2\beta$ and $|\ell| \geq 1$. In this case we have that, for every $x,z\in \mathbb{R}^n$ and $v\in (0,1)$, 
$$
\frac{(1+|z|)^{q'|\ell|} e^{-c(v|x+z|^2+\frac{|x-z|^2}{v})}}{v^{(\frac{n}{2}+|s|-\frac{|\ell|}{2}-2\beta)q'+1}}\leq C\frac{(1+|z|)^{\frac{q'}{2}}e^{-c(v|x+z|^2+\frac{|x-z|^2}{v})}}{v^{(\frac{n}{2}-\frac{1}{4})q'+1}}\leq \frac{(1+|z|)^{\frac{q'}{2}}e^{-c\frac{|x-z|^2}{v}}}{v^{(\frac{n}{2}-\frac{1}{4})q'+1}}.
$$ 
Hence, we deduce that
$$
\Big\|u^{2\beta}I_1(u,x,z)\Big\|_{L^{q'}((0,\infty),\frac{du}{u})}\leq C\sqrt{1+|z|}\left(\int_0^1\frac{e^{-c\frac{|x-z|^2}{v}}}{v^{(\frac{n}{2}-\frac{1}{4})q'+1}}dv\right)^{1/q'}\leq C\frac{\sqrt{1+|z|}}{|x-z|^{n-\frac{1}{2}}},\quad (x,z)\in N_1.
$$
On the other hand, we can write, for each $x,z\in \mathbb{R}^n$ and $u>0$,
\begin{align*}
    |I_2(u,x,z)|&\\
    &\hspace{-1.5cm}=
\frac{1}{2^{2\beta}\pi ^{\frac{n}{2}}}\left|\frac{e^{-(4\beta +n)u}}{(1-e^{-2u})^{\frac{n}{2}+2\beta}}\sum_{|r|=2\beta}{2\beta\choose {r_1\ldots r_n}}\left[\widetilde{H}_{2r}\Big(\frac{x-e^{-u}z}{\sqrt{1-e^{-2u}}}\Big)-\widetilde{H}_{2r}\Big(\frac{x-z}{\sqrt{2u}}\Big)\right]\right.\\
    &\quad 
    \left.
     +\sum_{|r|=2\beta}{2\beta\choose {r_1\ldots r_n}}\widetilde{H}_{2r}\Big(\frac{x-z}{\sqrt{2u}}\Big)\left[\frac{e^{-(4\beta +n)u}}{(1-e^{-2u})^{\frac{n}{2}+2\beta}}-\frac{1}{(2u)^{\frac{n}{2}+2\beta}}\right]\right|\\
     &\hspace{-1.5cm}\leq C\left(\frac{1}{u^{\frac{n}{2}+2\beta}}\sum_{|r|=2\beta}\left|\widetilde{H}_{2r}\Big(\frac{x-e^{-u}z}{\sqrt{1-e^{-2u}}}\Big)-\widetilde{H}_{2r}\Big(\frac{x-z}{\sqrt{2u}}\Big)\right|+e^{-c\frac{|x-z|^2}{u}}\left|\frac{e^{-(4\beta +n)u}}{(1-e^{-2u})^{\frac{n}{2}+2\beta}}-\frac{1}{(2u)^{\frac{n}{2}+2\beta}}\right|\right).
\end{align*}
We observe that, since $\widetilde{H}_n'(s)=-\widetilde{H}_{n+1}(s)$, $s\in \mathbb{R}$, $n\in \mathbb{N}$, and $e^{-\frac{|x-e^{-u}z|^2}{1-e^{-u}}}\leq Ce^{-c\frac{|x-z|^2}{u}}$, $(x,z)\in N_1$, $u\in (0,1)$, we get for every $r=(r_1,...,r_n)\in \mathbb{N}^n$,
\begin{align*} 
    \left|\widetilde{H}_{2r}\Big(\frac{x-e^{-u}z}{\sqrt{1-e^{-2u}}}\Big)-\widetilde{H}_{2r}\Big(\frac{x-z}{\sqrt{2u}}\Big)\right|&\\
    &\hspace{-5cm}\leq \sum_{k=1}^n\prod_{i=1}^{k-1}\left|\widetilde{H}_{2r_i}\Big(\frac{x_i-z_i}{\sqrt{2u}}\Big)\right|\left|\widetilde{H}_{2r_k}\Big(\frac{x_k-e^{-u}z_k}{\sqrt{1-e^{-2u}}}\Big)-\widetilde{H}_{2r_k}\Big(\frac{x_k-z_k}{\sqrt{2u}}\Big)\right|\prod_{i=k+1}^{n}\left|\widetilde{H}_{2r_i}\Big(\frac{x_i-e^{-u}z_i}{\sqrt{1-e^{-2u}}}\Big)\right|\\
    &\hspace{-5cm}\leq Ce^{-c\frac{|x-z|^2}{u}}\sum_{k=1}^n\left|\frac{x_k-e^{-u}z_k}{\sqrt{1-e^{-2u}}}-\frac{x_k-z_k}{\sqrt{2u}}\right|\\
    &\hspace{-5cm}\leq Ce^{-c\frac{|x-z|^2}{u}}\sum_{k=1}^n\left(|z_k|\sqrt{1-e^{-2u}}+|x_k-z_k|\Big|\frac{1}{\sqrt{1-e^{-2u}}}-\frac{1}{\sqrt{2u}}\Big|\right)\\
    &\hspace{-5cm}\leq Ce^{-c\frac{|x-z|^2}{u}}\sqrt{u}(1+|z|),\quad (x,z)\in N_1,\;u\in (0,1).
\end{align*}
We also have that, when $u\in (0,1)$
$$
\left|\frac{e^{-(4\beta +n)u}}{(1-e^{-2u})^{\frac{n}{2}+2\beta}}-\frac{1}{(2u)^{\frac{n}{2}+2\beta}}\right|\leq \left|\frac{e^{-(4\beta +n)u}-1}{(1-e^{-2u})^{\frac{n}{2}+2\beta}}\right|+\left|\frac{1}{(1-e^{-2u})^{\frac{n}{2}+2\beta}}-\frac{1}{(2u)^{\frac{n}{2}+2\beta}}\right|\leq \frac{C}{u^{\frac{n}{2}+2\beta -1}}.
$$
Thus, we conclude that
$$
|I_2(u,x,z)|\leq C\frac{e^{-c\frac{|x-z|^2}{u}}}{u^{\frac{n}{2}+2\beta}}(\sqrt{u}(1+|z|)+u)\leq C(1+|x|) \frac{e^{-c\frac{|x-z|^2}{u}}}{u^{\frac{n}{2}+2\beta-\frac{1}{2}}},\quad (x,z)\in N_1,\;u\in (0,1).
$$
Also we have that
$$
|I_2(u,x,z)|\leq C\left(\frac{e^{-(4\beta +n)u}}{(1-e^{-2u})^{\frac{n}{2}+2\beta}}+\frac{1}{u^{\frac{n}{2}+2\beta}}\right)\leq \frac{C}{u^{\frac{n}{2}+2\beta}},\quad x,z\in \mathbb{R}^n,\;u\in (0,\infty ).
$$
By these estimations we get, since $\sqrt{m(x)}\sim (1+|x|)^{-1}$, $x\in \mathbb{R}^n$,
\begin{align*}
\Big\|u^{2\beta}I_2(u,x,z)\Big\|_{L^{q'}((0,\infty),\frac{du}{u})}^{q'}&\leq C\left((1+|x|)^{q'}\int_0^{m(x)}\frac{e^{-c\frac{|x-z|^2}{u}}}{u^{(\frac{n}{2}-\frac{1}{2})q'+1}}du+\int_{m(x)}^\infty \frac{du}{u^{\frac{n}{2}q'+1}}\right)\\
&\hspace{-2cm}\leq C\left((1+|x|)^{q'}m(x)^{\frac{q'}{4}}\int_0^{m(x)}\frac{e^{-c\frac{|x-z|^2}{u}}}{u^{(\frac{n}{2}-\frac{1}{4})q'+1}}du+\frac{1}{m(x)^{\frac{n}{2}q'}}\right)\\
&\hspace{-2cm}\leq C\left(\frac{\sqrt{1+|x|}}{|x-y|^{n-\frac{1}{2}}}\right)^{q'},\quad (x,z)\in N_1.
\end{align*}
Hence, \eqref{objective} is established.

We have seen that $L_\beta^{\mathcal A}$ is bounded from $L^r_{L^q_X((0,\infty),\frac{dt}{t})}(\mathbb R^n,\gamma_{-1})$ into itself, with $1< r<\infty$, provided that so is $L_{\beta,{\rm loc}}$.

We are going to prove that $L_{\beta,{\rm loc}}$ is bounded from $L^r_{L^q_X((0,\infty),\frac{dt}{t})}(\mathbb R^n,\gamma_{-1})$ into itself, for every  $1< r<\infty$. We use vector valued Calder\'on-Zygmund theory.

We have that
$$
|\partial_t^{2\beta}W_t(z)|\leq C\frac{e^{-c|z|^2/t}}{t^{n/2+2\beta}},\;\;\;\;z\in\mathbb R^n\;\mbox{and}\;t>0.
$$
We get
\begin{align*}
\int_{\mathbb R^n}\int_0^\infty|K_\beta(t,x;s,z)|\frac{ds}{s}dz& \leq C \int_{\mathbb R^n}\int_0^\infty(st)^\beta\frac{e^{-c\frac{|x-z|^2}{s+t}}}{(s+t)^{\frac{n}{2}+2\beta}}\frac{ds}{s}dz \\
& \leq C\int_{\mathbb R^n}\int_0^\infty\frac{(st)^\beta}{\left(1+\frac{|x-z|^2}{s+t}\right)^{\frac{n}{2}+2\beta}(s+t)^{\frac{n}{2}+2\beta}}\frac{ds}{s}dz \\
& \leq C\int_{\mathbb R^n}\int_0^\infty\frac{(st)^\beta}{\left(s+t+|x-z|^2\right)^{\frac{n}{2}+2\beta}}\frac{ds}{s}dz \\
& \leq C\int_0^\infty (st)^\beta\int_0^\infty\frac{\rho^{n-1}}{\left(s+t+\rho^2\right)^{\frac{n}{2}+2\beta}}d\rho\frac{ds}{s} \\
& \leq C\int_0^\infty\frac{(st)^\beta}{(s+t)^{2\beta}}\frac{ds}{s}=C\int_0^\infty\frac{u^{\beta -1}}{(1+u)^{2\beta}}du=C,\quad x\in\mathbb R^n\;\mbox{and}\;t>0.
\end{align*}
It follows that
\begin{align*}
 \|L_\beta(h)\|^q_{L^q_{L^q_X((0,\infty),\frac{dt}{t})}(\mathbb R^n,dx)}&=\int_{\mathbb R^n}\int_0^\infty\left(\int_{\mathbb R^n}\int_0^\infty|K_\beta(t,x;s,z)|\|h(s,z)\|_X\frac{ds}{s}dz\right)^q\frac{dt}{t}dx \\
&\hspace{-3cm}\leq \int_{\mathbb R^n}\int_0^\infty\int_{\mathbb R^n}\int_0^\infty|K_\beta(t,x;s,z)|\|h(s,z)\|_X^q\frac{ds}{s}dz\left(\int_{\mathbb R^n}\int_0^\infty |K_\beta(t,x;s,z)|\frac{ds}{s}dz\right)^{q/q'}\frac{dt}{t}dx  \\
&\hspace{-3cm}\leq C\|h\|^q_{L^q_{L^q_X((0,\infty),\frac{dt}{t})}(\mathbb R^n,dx)},\quad h\in L^q_{L^q_X((0,\infty),\frac{dt}{t})}(\mathbb R^n,dx).
\end{align*}
For every $x,z\in\mathbb R^n$, $x\neq z$, we define the operator $K_\beta(x,z)$ by
$$
K_\beta(x,z)(g)(t)=\int_0^\infty K_\beta(t,x;s,z)g(s)\frac{ds}{s},\;\;\;\;t\in(0,\infty).
$$
As before we have that
$$
\int_0^\infty|K_\beta(t,x;s,z)|\frac{ds}{s}\leq C\int_0^\infty \frac{(st)^\beta}{(s+t+|x-z|^2)^{\frac{n}{2}+\beta}}\frac{ds}{s}\leq\frac{C}{|x-z|^{n}},\quad x,z\in\mathbb R^n,\,x\neq z.
$$
Then, for every $g\in L^q_X((0,\infty),\frac{dt}{t})$,
\begin{align*}
\|K_\beta(x,z)(g)\|_{ L^q_X((0,\infty),\frac{dt}{t})}^q&=\int_0^\infty\left|\int_0^\infty K_\beta(t,x;s,z)g(s)\frac{ds}{s}\right|^q\frac{dt}{t} \\
&\hspace{-2cm} \leq \int_0^\infty\left(\int_0^\infty  K(t,x;s,z)\frac{ds}{s}\right)^{q/q'}\int_0^\infty K_\beta(t,x;s,z)\|g(s)\|^q\frac{ds}{s}dt \\
&\hspace{-2cm} \leq \frac{C}{|x-z|^{(\frac{q}{q'}+1)n}}\|g\|^q_{L^q_X((0,\infty),\frac{dt}{t})}\leq \frac{C}{|x-z|^{qn}}\|g\|^q_{L^q_X((0,\infty),\frac{dt}{t})},\;\;\;\;x,z\in\mathbb R^n,\,x\neq z.
\end{align*}
Hence,
$$
\|K_\beta(x,z)\|_{ L^q_X((0,\infty),\frac{dt}{t})\rightarrow L^q_X((0,\infty),\frac{dt}{t})}\leq\frac{C}{|x-z|^n},\quad x,z\in\mathbb R^n,\,x\neq z.
$$
In a similar way we can see that, for every $i=1,...,n$,
$$
\|\partial_{x_i}K_\beta(x,z)\|_{L^q_X((0,\infty),\frac{dt}{t})\rightarrow L^q_X((0,\infty),\frac{dt}{t})}\leq\frac{C}{|x-z|^{n+1}},\quad x,z\in\mathbb R^n,\,x\neq z,
$$
and
$$
\|\partial_{z_i}K_\beta(x,z)\|_{ L^q_X((0,\infty),\frac{dt}{t})\rightarrow L^q_X((0,\infty),\frac{dt}{t})}\leq\frac{C}{|x-z|^{n+1}},\quad x,z\in\mathbb R^n,\, x\neq z,
$$
where $\partial_{x_i}K_\beta(x,z)$ and $\partial_{z_i}K_\beta(x,z)$ are understood as the integral operators defined by the kernels $\partial_{x_i}K_\beta(t,x;s,z)$ and $\partial_{z_i}K_\beta(t,x;s,z)$, respectively.

We also have that, for every $h\in C_c(\mathbb R^n)\otimes C_c(0,\infty)\otimes X$,
$$
K_\beta(h)(t,x)=\left(\int_{\mathbb R^n}K_\beta(x,z)[h(\cdot,x)]dz\right)(t),
$$
for almost all $(t,x)\notin \supp h$. Here the integral is understood in the $L^2((0,\infty),\frac{dt}{t})$-Bochner sense.
According to vector valued Calder\'on-Zygmund theory , $K_\beta$ defines a bounded operator from $L^r_{L^q_X((0,\infty),\frac{dt}{t})}(\mathbb R^n,dx)$ into itself, for every $1<r<\infty$.

By \cite[Propositions 3.2.5 and 3.2.7]{Sa}  $K_{\beta,{\rm loc}}$ is bounded from $L^r_{L^q_X((0,\infty),\frac{dt}{t})}(\mathbb R^n,\gamma_{-1})$ into itself, for every $1<r<\infty$.

The proof of (ii) $\Rightarrow$ (i) is thus finished.

In order to prove that (iii) $\Rightarrow$ (i) we can proceed as in the proof of (ii) $\Rightarrow$ (i) with the semigroup $\{T_t^{\mathcal A}\}_{t>0}$ replaced by $\{P_t^{\mathcal A}\}_{t>0}$. We can also argue as follows. Suppose that (iii) is true. According to (\ref{gPoissongHeat}) and by taking into account that $g_{\beta_1,0,\{T_t^\mathcal{A}\}_{t>0}}^{q,X}(f)\leq g_{\beta_2,0,\{T_t^\mathcal{A}\}_{t>0}}^{q,X}(f)$ when $0<\beta_1\leq \beta_2$, we have that
$$
\|f\|_{L^p_X(\mathbb R^n,\gamma_{-1})}\leq C \|g_{m,0,\{T_t^{\mathcal A}\}_{t>0}}^{q,X}(f)\|_{L^p(\mathbb R^n,\gamma_{-1})},\;\;\;\;f\in L^p(\mathbb R^n,\gamma_{-1}),
$$
where $m\in\mathbb N$ and $m-1\leq\beta <m$. Then by using the property we have just proved, (i) holds.

\subsection{Proof of Corollary \ref{cor1.4}}
The equivalence stated in Corollary \ref{cor1.4} follows by using Theorems \ref{Th1.2} and \ref{Th1.3} and by taking into account that $X$ is isomorphic to a Hilbert space if and only if $X$ has 2-martingale type and 2-martingale cotype.


\section{ Proofs of Theorems \ref{Th1.5} and \ref{Th1.6}}

\subsection{Proof of Theorem \ref{Th1.5}} Since $\{P_t^\mathcal{A}\}_{t>0}$ is a subordinated Poisson semigroup of the symmetric diffusion semigroup $\{T_t^\mathcal{A}\}_{t>0}$ the property $(i)\Rightarrow (ii)$ can be deduced from \cite[Theorem 1.6]{Hy}. We recall that $X^*$ is UMD provided that $X$ is UMD.

Assume now that $(ii)$ is true. We are going to see that $X$ is UMD. Our first objective is to see that there exists $C>0$ such that
$$
\|g_{\{P_t\}_{t>0}}^{1,X}(f)\|_{L^p_X(\mathbb{R}^n,\gamma_{-1})}\leq C\|f\|_{L^p_X(\mathbb{R}^n,\gamma_{-1})},\quad f\in L^p_X(\mathbb{R}^n,\gamma_{-1}).
$$
We define the global and local operators in the usual way and consider the operator
\begin{align*}
G(f)(t,x,w) &=t\big[\partial _tP_{t,{\rm loc}}^{\mathcal A}(f(\cdot ,w))(x)-\partial _tP_{t,{\rm loc}}(f(\cdot, w))(x)\big]\\
&=\int_{\mathbb{R}^n}\mathcal{X}_{N_1}(x,y)t\partial _t[P_t^\mathcal{A}(x,y)-P_t(x-y)]f(y,w)dy,\quad x\in \mathbb{R}^n,\;t>0, \;w\in \Omega.
\end{align*}
According to \eqref{difPoisson} we have that
$$
H(x,y)=\left\|t\partial _t[P_t^\mathcal{A}(x,y)-P_t(x-y)]\right\|_{L^2((0,\infty ),\frac{dt}{t})}\leq C\frac{\sqrt{1+|x|}}{|x-y|^{n-\frac{1}{2}}},\quad (x,y)\in N_1.
$$
Then, the operator $\mathfrak{H}$ defined by
$$
\mathfrak{H}(f)(x,w)=\int_{\mathbb{R}^n}H(x,y)f(y,w)\mathcal{X}_{N_1}(x,y)dy,\quad x\in \mathbb{R}^n,\;w\in \Omega,
$$
is bounded from $L^r_X(\mathbb{R}^n,dx)$ into itself, for every $1\leq r\leq \infty$.

By using Minkowski inequality we deduce that
$$
\Big\|\|G(f)(\cdot , x,w)\|_{L^2((0,\infty ),\frac{dt}{t})}\Big\|_{L^r_X(\mathbb{R}^n,dx)}\leq C\|f\|_{L^r_X(\mathbb{R}^n,dx)}, $$
for every $f\in L^r_X(\mathbb{R}^n,dx)$, $1\leq r<\infty$. Then, since the operators are local (see \cite[Proposition 3.2.5]{Sa})
$$
\big\|g_{\{P_t^\mathcal{A}\}_{t>0};\,{\rm loc}}^{1,X}(f)-g_{\{P_t\}_{t>0};\,{\rm loc}}^{1,X}(f)\big\|_{L^r_X(\mathbb{R}^n,\gamma_{-1})}\leq C\|f\|_{L^r_X(\mathbb{R}^n,\gamma_{-1})},
$$
for every $f\in L^r_X(\mathbb{R}^n,\gamma_{-1})$, $1\leq r<\infty$.
We now consider the operator $\mathcal{M}_{\rm glob}$ defined by
$$
\mathcal{M}_{\rm glob}(f)(t,x,w)=\int_{\mathbb{R}^n}t\partial _tP_t^\mathcal{A}(x,y)f(y,w)\mathcal{X}_{N_1^c}(x,y)dy,\quad x\in \mathbb{R}^n,\;t>0,\;w\in \Omega.
$$

As in \eqref{gPoissongHeat} we obtain
$$
\big\|t\partial _tP_t^\mathcal{A}(x,y)\big\|_{L^2((0,\infty ),\frac{dt}{t})}\leq C\big\|t\partial _tT_t^\mathcal{A}(x,y)\big\|_{L^2((0,\infty ),\frac{dt}{t})},\quad x,y\in \mathbb{R}^n.
$$

According to (\ref{A2}) when $m=1$, $k=0$ and $q=2$ it follows that, for every $r\in (1,\infty)$ there exists a positive function $F_r$ defined in $\mathbb{R}^n\times \mathbb{R}^n$ such that
$$
\|t\partial_tT_t^\mathcal{A}(x,y)\|_{L^2((0,\infty ),\frac{dt}{t})}\leq F_r(x,y),\quad (x,y)\in N_1^c,
$$
being 
$$
\sup_{x\in \mathbb{R}^n}\int_{\mathbb{R}^n}e^{\frac{|x|^2-|y|^2}{r}}F_r(x,y)\mathcal{X}_{N_1^c}(x,y)dy<\infty,
$$
and
$$
\sup_{y\in \mathbb{R}^n}\int_{\mathbb{R}^n}e^{\frac{|x|^2-|y|^2}{r}}F_r(x,y)\mathcal{X}_{N_1^c}(x,y)dx<\infty.
$$
Then,
$$
\Big\|\|\mathcal{M}_{\rm glob}(f)(\cdot ,x,w)\|_{L^2((0,\infty ),\frac{dt}{t})}\Big\|_{L^r_X(\mathbb{R}^n,\gamma_{-1})}\leq C\|f\|_{L^r_X(\mathbb{R}^n,\gamma_{-1})},
$$
for every $f\in L^r(\mathbb{R}^n,\gamma_{-1})$ and $r\in (1,\infty )$.

It follows that 
$$
\|g_{\{P_t^\mathcal{A}\}_{t>0};\,{\rm glob}}^{1,X}(f)\|_{L^r_X(\mathbb{R}^n,\gamma_{-1})}\leq C\|f\|_{L^r_X(\mathbb{R}^n,\gamma_{-1})},
$$
for each $f\in L^r(\mathbb{R}^n,\gamma_{-1})$ and $r\in (1,\infty )$.

Since $(ii)$ holds we have that
$$
\|g_{\{P_t^\mathcal{A}\}_{t>0};\,{\rm loc}}^{1,X}(f)\|_{L^p_X(\mathbb{R}^n,\gamma_{-1})}\leq C\|f\|_{L^p_X(\mathbb{R}^n,\gamma_{-1})},\quad f\in L^p_X(\mathbb{R}^n,\gamma_{-1}).
$$
We can write 
\begin{align*}
\|g_{\{P_t\}_{t>0};\,{\rm loc}}^{1,X}(f)\|_{L^p_X(\mathbb{R}^n,\gamma_{-1})}&\leq \|g_{\{P_t^\mathcal{A}\}_{t>0};\,{\rm loc}}^{1,X}(f)-g_{\{P_t\}_{t>0};\,{\rm loc}}^{1,X}(f)\|_{L^p_X(\mathbb{R}^n,\gamma_{-1})}\\
&\quad +\|g_{\{P_t^\mathcal{A}\}_{t>0};\,{\rm loc}}^{1,X}(f)\|_{L^p_X(\mathbb{R}^n,\gamma_{-1})}.
\end{align*}
Hence, $g_{\{P_t\}_{t>0};\,{\rm loc}}^{1,X}$ is bounded from $L^p_X(\mathbb{R}^n,\gamma_{-1})$ into itself, and also from $L^p_X(\mathbb{R}^n,dx)$ into itself (see \cite[Proposition 3.2.5]{Sa}).

By using now the corresponding properties of invariance (see \cite[p. 21]{HTV}) we deduce that
$$
\|g_{\{P_t\}_{t>0}}^{1,X}(f)\|_{L^p_X(\mathbb{R}^n,dx)}\leq C\|f\|_{L^p_X(\mathbb{R}^n,dx)},\quad f\in L^p_X(\mathbb{R}^n,dx).
$$
By changing $X$ by $X^*$ and $p$ by $p'$ we also infer from $(ii)$ that
$$
\|g_{\{P_t\}_{t>0}}^{1,X^*}(f)\|_{L^{p'}_{X^*}(\mathbb{R}^n,dx)}\leq C\|f\|_{L^{p'}_{X^*}(\mathbb{R}^n,dx)},\quad f\in L^{p'}_{X^*}(\mathbb{R}^n,dx).
$$
Let now $f\in L^2(\mathbb{R}^n,\gamma_{-1})\otimes X$ and $h\in L^2(\mathbb{R}^n,\gamma_{-1})\otimes X^*$. Since $X$ is order continuous the K\"othe and topological dual coincides. We can write the following polarization equality
\begin{align*}
    \int_{\mathbb{R}^n}\langle f(x,\cdot),h(x,\cdot)\rangle _{X,X^*}dx&=4\int_{\mathbb{R}^n}\int_0^\infty \langle t\partial _tP_t(f(\cdot ,w))(x),t\partial _tP_t(h(\cdot ,w))(x)\rangle_{X,X^*}\frac{dt}{t}dx\\
    &=4\int_{\mathbb{R}^n}\int_0^\infty \int_{\Omega}t\partial _tP_t(f(\cdot ,w))(x)t\partial _tP_t(h(\cdot ,w))(x)d\mu (w)\frac{dt}{t}dx.
\end{align*}
By using H\"older inequality it follows that
\begin{align*}
\left|\int_{\mathbb{R}^n}\langle f(x,\cdot),h(x,\cdot )\rangle_{X,X^*}dx\right|&\\
&\hspace{-4cm}\leq C\int_{\mathbb{R}^n}\int_{\Omega}\big\|t\partial_tP_t(f(\cdot ,w))(x)\big\|_{L^2((0,\infty ),\frac{dt}{t})}\big\|t\partial_tP_t(h(\cdot ,w))(x)\big\|_{L^2((0,\infty ),\frac{dt}{t})}d\mu (w)dx\\
&\hspace{-4cm}\leq C\int_{\mathbb{R}^n}\big\|\|t\partial_tP_t(f(\cdot ,w))(x)\|_{L^2((0,\infty ),\frac{dt}{t})}\|_X\big\|\|t\partial_tP_t(h(\cdot ,w))(x)\|_{L^2((0,\infty ),\frac{dt}{t})}\big\|_{X^*}dx\\
&\hspace{-4cm}\leq C\Big\|g_{\{P_t\}_{t>0}}^{1,X}(f)\Big\|_{L^p(\mathbb{R}^n,dx)}\Big\|g_{\{P_t\}_{t>0}}^{1,X^*}(h)\Big\|_{L^{p'}_{X^*}(\mathbb{R}^n,dx)}\leq C\Big\|g_{\{P_t\}_{t>0}}^{1,X}(f)\Big\|_{L^p(\mathbb{R}^n,dx)}\|h\|_{L^{p'}_{X^*}(\mathbb{R}^n,dx)}.
\end{align*}
Since $L^{p'}_{X^*}(\mathbb{R}^n,dx)$ is norming in $L^p_X(\mathbb{R}^n,dx)$ it follows that
$$
\|f\|_{L^p_X(\mathbb{R}^n,dx)}\leq C\big\|g_{\{P_t\}_{t>0}}^{1,X}(f)\big\|_{L^p_X(\mathbb{R}^n,dx)}.
$$
From \cite[Theorem 1.6]{Hy} we deduce that $X$ is UMD and the proof is complete.

\noindent {\bf Remark}. The property established in \cite[Theorem 1.6]{Hy} suggests to ask whether in Theorem \ref{Th1.5} $(ii)$ we can replace the condition involving the dual space $X^*$ by this other one
$$
\|f\|_{L^p_X(\mathbb{R}^n,\gamma_{-1})}\leq C\big\|g_{\{P_t^\mathcal{A}\}_{t>0}}^{1,X}(f)\big\|_{L^p_X(\mathbb{R}^n,\gamma_{-1})},\quad f\in L^p_X(\mathbb{R}^n,\gamma_{-1}),
$$
and then the UMD property for $X$ can be again deduced  from $(ii)$. At this moment we do not know what is the answer for this question.

\subsection{Proof of Theorem \ref{Th1.6}}
In \cite[Theorems 1.4 and 1.5]{BR} the Banach spaces with the UMD-property were characterized by using $L^p(\mathbb{R}^n,\gamma_{-1})$-boundedness properties of the Riesz transforms associate with the operator $\mathcal A$. For every $i=1,...,n$, the Riesz transform $\mathcal{R}_i^{\mathcal{A}}$ in the $\mathcal{A}$-setting is defined by
$$
\mathcal{R}_i^\mathcal{A}(f)(x)=\lim_{\varepsilon \rightarrow 0^+}\int_{|x-y|>\varepsilon}\mathcal{R}_i^\mathcal{A}(x,y)f(y)dy,\quad \mbox{ for almost all }x\in \mathbb{R}^n,
$$
for every $f\in L^p(\mathbb{R}^n, \gamma_{-1})$, $1\leq p<\infty$, where
$$
\mathcal{R}_i^\mathcal{A}(x,y)=\frac{1}{\sqrt{\pi}}\int_0^\infty \partial_{x_i}T_t^\mathcal{A}(x,y)\frac{dt}{\sqrt{t}},\quad x,y\in \mathbb{R}^n.
$$
For every $i=1,...,n$, $\mathcal{R}_i^\mathcal{A}$ is bounded from $L^p(\mathbb{R}^n,\gamma_{-1})$ into itself, when $1<p<\infty$, and from $L^1(\mathbb{R}^n,\gamma_{-1})$ into $L^{1,\infty}(\mathbb{R}^n,\gamma_{-1})$; moreover, $\mathcal{R}_i^\mathcal{A}$ can be extended, in the usual way, to $L^p(\mathbb{R}^n,\gamma_{-1})\otimes Y$, $1\leq p<\infty$, where $Y$ is a Banach space. In \cite[Theorem 1.4]{BR} it was established that $Y$ has the UMD property if and only if  $\mathcal{R}_i^\mathcal{A}$, $i=1,...,n$, can be extended to $L^p_Y(\mathbb{R}^n,\gamma_{-1})$ for some (equivalently, for any) $1<p<\infty$ as a bounded operator from $L^p_Y(\mathbb{R}^n,\gamma_{-1})$ into itself.

Suppose that $(ii)$ holds. Let $i=1,...,n$. We have that, for every $k\in \mathbb{N}^n$, $\partial_{x_i}\widetilde{H}_k=-\widetilde{H}_{k+e_i}$, where $e_i=(e_{ij})_{j=1}^n$, being $e_{ij}=1$, $i=j$, and $e_{ij}=0$, if $i\not =j$. Then, since 
$$
\mathcal{A}^{-1/2}f=\sum_{\ell\in \mathbb{N}^n}\frac{1}{\sqrt{n+|\ell|}}c_\ell(f)\widetilde{H}_\ell,\quad f\in L^2(\mathbb{R}^n, \gamma_{-1}),
$$
we have that
$$
\mathcal{R}_i^\mathcal{A}\widetilde{H}_k=-\frac{1}{\sqrt{n+|k|}}\widetilde{H}_{k+e_i},\quad k\in \mathbb{N}^n.
$$
On the other hand, we can write
$$
P_t^\mathcal{A-I}(f)=\sum_{\ell\in \mathbb{N}^n}e^{-t\sqrt{n+|\ell|-1}}c_\ell(f)\widetilde{H}_\ell,\quad f\in L^2(\mathbb{R}^n, \gamma_{-1}).
$$
Then we get
$$
\partial _tP_t^\mathcal{A-I}(\mathcal{R}_i^\mathcal{A}\widetilde{H}_k)=-\frac{1}{\sqrt{n+|k|}}\partial _tP_t^\mathcal{A-I}(\widetilde{H}_{k+e_i})=e^{-t\sqrt{n+|k|}}\widetilde{H}_{k+e_i},\quad k\in \mathbb{N}^n.
$$
and also, by considering that 
$$
P_t^\mathcal{A}(f)=\sum_{\ell\in \mathbb{N}^n}e^{-t\sqrt{n+|\ell|}}c_\ell(f)\widetilde{H}_\ell,\quad f\in L^2(\mathbb{R}^n, \gamma_{-1}),
$$
it follows that
$$
\partial _{x_i}P_t^\mathcal{A}(\widetilde{H}_k)=-e^{-t\sqrt{n+|k|}}\widetilde{H}_{k+e_i},\quad k\in \mathbb{N}^n.
$$
Let us denote by $F={\rm span }\, \{\widetilde{H}_k\}_{k\in \mathbb{N}^n}$ the linear space generated by $\{\widetilde{H}_k\}_{k\in \mathbb{N}^n}$. For every $f\in F$ we have that
\begin{equation}\label{4.2.1}
\partial _{x_i}P_t^\mathcal{A}(f)=-\partial _tP_t^\mathcal{A-I}(\mathcal{R}_i^\mathcal{A}(f)).
\end{equation}
We recall that $F$ is a dense subspace of $L^p(\mathbb{R}^n,\gamma_{-1})$, for every $1\leq p<\infty$. It is clear that $\mathcal{R}_i^\mathcal{A}f\in F$, for every $f\in F$ and $i=1,...,n$. 

\noindent Since $X$ is order continuous, and by taking into account that, when $n\geq 2$, 0 is not an eigenvalue of $\mathcal{A}-I$, it follows that, for every $f\in L^2(\mathbb{R}^n,\gamma_{-1})\otimes X$ and $g\in L^2(\mathbb{R}^n,\gamma_{-1})\otimes X^*$, 
$$
\int_{\mathbb{R}^n}\langle f(x,\cdot),g(x,\cdot)\rangle_{X,X^*}d\gamma_{-1}(x)=4\int_{\mathbb{R}^n}\int_0^\infty \langle t\partial _tP_t^\mathcal{A-I}(f(\cdot,w))(x),t\partial _t P_t^\mathcal{A-I}(g(\cdot,w))(x)\rangle\frac{dt}{t}d\gamma_{-1}(x).
$$
As above, by taking into account that
$$
\big\|g_{\{P_t^\mathcal{A-I}\}_{t>0}}^{1,X^*}(h)\big\|_{L^{p'}_{X^*}(\mathbb{R}^n,\gamma_{-1})}\leq C\|h\|_{L^{p'}_{X^*}(\mathbb{R}^n,\gamma_{-1})},\quad h\in L^{p'}_{X^*}(\mathbb{R}^n,\gamma_{-1}),
$$
we obtain that, for every $f\in (L^2(\mathbb{R}^n,\gamma_{-1})\cap L^p(\mathbb{R}^n,\gamma_{-1}))\otimes X$,
\begin{equation}\label{4.2.2}
    \|f\|_{L^p_X(\mathbb{R}^n,\gamma_{-1})}\leq C\big\|g_{\{P_t^{\mathcal{A}-I}\}_{t>0}}^{1,X}(f)\big\|_{L^p_X(\mathbb{R}^n,\gamma_{-1})}.
\end{equation}
According to \eqref{4.2.1} and \eqref{4.2.2} we get
\begin{align*}
\big\|\mathcal{R}_i^\mathcal{A}(f)\big\|_{L^p_X(\mathbb{R}^n,\gamma_{-1})}&\leq C\big\|g_{\{P_t^\mathcal{A-I}\}_{t>0}}^{1,X}(\mathcal{R}_i^\mathcal{A}(f))\big\|_{L^p_X(\mathbb{R}^n,\gamma_{-1})}= C\big\|g_{i,\{P_t^\mathcal{A}\}_{t>0}}^X(f)\big\|_{L^p_X(\mathbb{R}^n,\gamma_{-1})}\\
&\leq C\|f\|_{L^p_X(\mathbb{R}^n,\gamma_{-1})},\quad f\in F\otimes X.
\end{align*}
By \cite[Theorem 1.4]{BR} we conclude that the space $X$ has the UMD property.

We observe that when $n=1$, we have that
\begin{align*}
\int_{\mathbb{R}}\langle (f-E_0(f))(x,\cdot),(g-E_0(g))(x,\cdot)\rangle_{X,X^*}d\gamma_{-1}(x)&\\
&\hspace{-4cm}=4\int_{\mathbb{R}}\int_0^\infty \langle t\partial _tP_t^\mathcal{A-I}(f(\cdot,w))(x),t\partial _t P_t^\mathcal{A-I}(g(\cdot,w))(x)\rangle\frac{dt}{t}d\gamma_{-1}(x),
\end{align*}
for every $f\in L^2(\mathbb{R},\gamma_{-1})\otimes X$ and $g\in L^2(\mathbb{R},\gamma_{-1})\otimes X^*$, where $E_0(f)$ is the projection of $f$ into the subspace generated by $\widetilde{H}_0$. Property \eqref{4.2.2} can be written in the following way
$$
\|f-E_0(f)\|_{L^p_X(\mathbb{R},\gamma_{-1})}\leq C\big\|g_{\{P_t^{\mathcal{A}-I}\}_{t>0}}^{1,X}(f)\big\|_{L^p_X(\mathbb{R},\gamma_{-1})},
$$
 for every $f\in (L^2(\mathbb{R},\gamma_{-1})\cap L^p(\mathbb{R},\gamma_{-1}))\otimes X$. Since $E_0(\mathcal{R}_1^\mathcal{A}(f))=0$, for $f\in F\otimes X$, with the same argument as before we can deduce that $X$ is UMD when $n=1$.
 
Suppose now that the Banach space $X$ has the UMD property. Let $i=1,...,n$. Our objective is to see that, for certain $C>0$, 
\begin{equation}\label{O1}
\big\|g_{\{P_t^{\mathcal{A}-I}\}_{t>0}}^{1,X^*}(f)\big\|_{L^{p'}_{X^*}(\mathbb{R}^n,\gamma_{-1})}\leq C\|f\|_{L^{p'}_{X^*}(\mathbb{R}^n,\gamma_{-1})},\quad f\in L^{p'}_{X^*}(\mathbb{R}^n,\gamma_{-1}),
\end{equation}
and
\begin{equation}\label{O2}
\big\|g_{i,\{P_t\}_{t>0}}^X(f)\big\|_{L^p_X(\mathbb{R}^n,\gamma_{-1})}\leq C\|f\|_{L^p_X(\mathbb{R}^n,\gamma_{-1})}, \quad f\in L^p_X(\mathbb{R}^n,\gamma_{-1}).
\end{equation}
Let us consider first the estimation for $g_{\{P_t^{\mathcal{A}-I}\}_{t>0}}^{1,X^*}$. We define the local and global operators in the usual way and write
$$
    g_{\{P_t^\mathcal{A-I}\}_{t>0}}^{1,X^*}(f)\leq g_{\{P_t^\mathcal{A-I}\}_{t>0};\,{\rm glob}}^{1,X^*}(f)+g_{\{P_t^\mathcal{A-I}\}_{t>0};\,{\rm loc}}^{1,X^*}(f)-g_{\{P_t\}_{t>0};\,{\rm loc}}^{1,X^*}(f)+g_{\{P_t\}_{t>0};\,{\rm loc}}^{1,X^*}(f),
$$
for every $f\in L_X^p(\mathbb{R}^n,\gamma_{-1})$.

By using Minkowski inequality we get, for each $x\in \mathbb{R}^n$ and $w\in \Omega$,
$$
g_{\{P_t^\mathcal{A-I}\}_{t>0};\,{\rm glob}}^{1,X^{*}}(f(\cdot, w))(x)\leq C\int_{\mathbb{R}^n}\mathcal{X}_{N_1^c}(x,y)|f(y,w)|\big\|t\partial _tP_t^\mathcal{A-I}(x,y)\big\|_{L^2((0,\infty )\frac{dt}{t})}dy,
$$
and 
\begin{align*}
    \big|g_{\{P_t^\mathcal{A-I}\}_{t>0};\,{\rm loc}}^{1,X^*}(f(\cdot, w))(x)-g_{\{P_t\}_{t>0};\,{\rm loc}}^{1,X^*}(f(\cdot , w))(x)\big|&\\
    &\hspace{-6cm}\leq C \int_{\mathbb{R}^n}\mathcal{X}_{N_1}(x,y)|f(y,w)|\big\|t\partial _t\big[P_t^\mathcal{A-I}(x,y)-P_t(x-y)\big]\big\|_{L^2((0,\infty )\frac{dt}{t})}dy.
\end{align*}

We recall that 
$$
P_t^\mathcal{A-I}(x,y)=\frac{1}{\sqrt{\pi}}\int_0^\infty e^{-u}T_{t^2/(4u)}^\mathcal{A-I}(x,y)\frac{du}{\sqrt{u}},\quad x,y\in \mathbb{R}^n, \;t>0.
$$
We obtain, for every $x,y\in \mathbb{R}^n$ and $t>0$,
$$
t\partial _tP_t^\mathcal{A-I}(x,y)=\frac{t^2}{2\sqrt{\pi}}\int_0^\infty e^{-u}\partial_vT_v^{\mathcal{A}-I}(x,y)_{|v=\frac{t^2}{4u}}\frac{du}{u^{\frac{3}{2}}}=\frac{t}{\sqrt{\pi}}\int_0^\infty e^{-\frac{t^2}{4v}}\partial_vT_v^{\mathcal{A}-I}(x,y)\frac{dv}{\sqrt{v}}.
$$
Then, Minkowski inequality leads to
\begin{align*}
\big\|t\partial _tP_t^\mathcal{A-I}(x,y)\big\|_{L^2((0,\infty ),\frac{dt}{t})}&\leq C\int_0^\infty \big\|te^{-\frac{t^2}{4v}}\big\|_{L^2((0,\infty ),\frac{dt}{t})}|\partial _vT_v^\mathcal{A-I}(x,y)|\frac{dv}{\sqrt{v}}\\
&\hspace{-1cm}\leq C \int_0^\infty |\partial _vT_v^\mathcal{A-I}(x,y)|dv=\big\|v\partial_vT_v^{\mathcal{A}-I}(x,y)\big\|_{L^1((0,\infty ),\frac{dv}{v})},\quad x,y\in \mathbb{R}^n.
\end{align*}
We now observe that, since $n\geq 2$, the estimation \eqref{A2} is also valid when $m=1$, $k=0$, $q=1$ and $\mathcal{A}$ is replaced by $\mathcal{A}-I$. Then, for $0<\delta<\eta <1$, we have that
$$
\big\|t\partial _tP_t^\mathcal{A-I}(x,y)\big\|_{L^2((0,\infty ),\frac{dt}{t})}\leq C \left\{\begin{array}{ll}
e^{-(\eta+\delta)\frac{|x|^2}{2} +(\eta-\delta)\frac{|y|^2}{2}},&\langle x,y \rangle \leq  0, \\[0.2cm]
|x+y|^ne^{\frac{\eta}{2}(|y|^2-|x|^2)-\frac{\delta}{2}|x+y||x-y|},&\langle x,y \rangle > 0. \end{array}\right.,\quad (x,y)\in N_1^c.
$$
By choosing $\eta-\delta<\frac{2}{p'}<\eta+\delta$ we deduce that
$$
\sup_{x\in \mathbb{R}^n}\int_{\mathbb{R}^n}e^{\frac{|x|^2-|y|^2}{p'}}\mathcal{X}_{N_1^c}(x,y)\big\|t\partial _tP_t^\mathcal{A-I}(x,y)\big\|_{L^2((0,\infty ),\frac{dt}{t})}dy<\infty,
$$
and
$$
\sup_{y\in \mathbb{R}^n}\int_{\mathbb{R}^n}e^{\frac{|x|^2-|y|^2}{p'}}\mathcal{X}_{N_1^c}(x,y)\big\|t\partial _tP_t^\mathcal{A-I}(x,y)\big\|_{L^2((0,\infty ),\frac{dt}{t})}dx<\infty.
$$
It follows that the operator $g_{\{P_t^\mathcal{A-I}\}_{t>0};\,{\rm glob}}^{1,X^*}$ is bounded from $L^{p'}_{X^*}(\mathbb{R}^n,\gamma_{-1})$ into itself.

On the other hand, by proceeding as above and by taking into account that estimation \eqref{diferencia} holds for $q=1$ and $\mathcal{A}-I$ (because $n\geq 2$) we obtain
\begin{align*}
\big\|t\partial _t\big[P_t^\mathcal{A-I}(x,y)-P_t(x-y)\big]\big\|_{L^2((0,\infty )\frac{dt}{t})}&\leq C \Big\|v\partial_v[T_v^{\mathcal{A}-I}(x,y)-W_v(x-y)]\Big\|_{L^1((0,\infty ),\frac{dv}{v})}\\
&\leq C\frac{\sqrt{1+|x|}}{|x-y|^{n-\frac{1}{2}}},\quad (x,y)\in N_1.
\end{align*}
We deduce that 
$$
 \big\|g_{\{P_t^\mathcal{A-I}\}_{t>0};\,{\rm loc}}^{1,X^*}(f(\cdot, w))(x)-g_{\{P_t\}_{t>0};\,{\rm loc}}^{1,X^*}(f(\cdot , w))(x)\big\|_{L^{p'}_{X^*}(\mathbb{R}^n,\gamma_{-1})}\leq C\|f\|_{L^{p'}_{X^*}(\mathbb{R}^n,\gamma_{-1})},
 $$
for every $f\in L^{p'}_{X^*}(\mathbb{R}^n,\gamma_{-1})$.

Finally, since $X^*$ has the UMD property, \cite[Theorem 1.6]{Hy} leads to $g_{\{P_t\}_{t>0}}^{1,X^*}$ is a bounded operator from $L^{p'}_{X^*}(\mathbb{R}^n,\gamma_{-1})$ into itself. 
We can also write
\begin{align*}
    \big\|t\partial _tP_t(x-y)\big\|_{L^2((0,\infty )\frac{dt}{t})}&\leq C\int_0^\infty |\partial _vW_v(x-y)|dv\leq C\int_0^\infty \frac{e^{-c\frac{|x-y|^2}{v}}}{v^{\frac{n}{2}+1}}dv\\
    &\leq \frac{C}{|x-y|^{n}}\leq C\frac{\sqrt{1+|x|}}{|x-y|^{n-\frac{1}{2}}},\quad (x,y)\in N_1^c.
\end{align*}
Thus, $g_{\{P_t\}_{t>0};\,{\rm glob}}^{1,X^*}$, and as a consequence $g_{\{P_t\}_{t>0};\,{\rm loc}}^{1,X^*}$, are  bounded operators from $L^{p'}_{X^*}(\mathbb{R}^n,\gamma_{-1})$ into itself. Thus \eqref{O1} is established.

In order to prove \eqref{O2} we proceed in the same way. In this case we write
\begin{align*}
\big\|t\partial_{x_i}P_t^\mathcal{A}(x,y)\big\|_{L^2((0,\infty ),\frac{dt}{t})}&= \Big\|\frac{t^2}{2\sqrt{\pi}}\int_0^\infty \frac{e^{-\frac{t^2}{4u}}}{u^\frac{3}{2}}\partial _{x_i}T_u^\mathcal{A}(x,y)du\Big\|_{L^2((0,\infty ),\frac{dt}{t})}\\
&\leq C\int_0^\infty \big\|t^2e^{-\frac{t^2}{4u}}\big\|_{L^2((0,\infty ),\frac{dt}{t})}|\partial _{x_i}T_u^\mathcal{A}(x,y)|\frac{du}{u^\frac{3}{2}}\\
&\leq C\int_0^\infty |\partial _{x_i}T_u^\mathcal{A}(x,y)|\frac{du}{\sqrt{u}},\quad x,y\in \mathbb{R}^n.
\end{align*}
Since
$$
\partial _{x_i}T_u^\mathcal{A}(x,y)=\frac{2e^{-nu}(x_i-e^{-u}y_i)}{(1-e^{-2u})^{\frac{n}{2}+1}}e^{-\frac{|x-e^{-u}y|^2}{1-e^{-2u}}},\quad x,y\in \mathbb{R}^n,\,u>0,
$$
it follows that
\begin{align*}
    \big\|t\partial _{x_i}P_t^\mathcal{A}(x,y)\big\|_{L^2((0,\infty ),\frac{dt}{t})}&\leq C\int_0^\infty e^{-nu}\frac{e^{-\eta\frac{|x-e^{-u}y|^2}{1-e^{-2u}}}}{(1-e^{-2u})^{\frac{n+1}{2}}}\frac{du}{\sqrt{u}}\\
    &\leq Ce^{-\eta(|x|^2-|y|^2)}\int_0^\infty e^{-nu}\frac{e^{-\eta\frac{|y-e^{-u}x|^2}{1-e^{-2u}}}}{(1-e^{-2u})^{\frac{n+1}{2}}}\frac{du}{\sqrt{u}},\quad x,y\in \mathbb{R}^n,
\end{align*}
with $0<\eta <1$.

By proceeding as in the estimation of $R_{m,k}^\mathcal{A}$ in \eqref{A2} we get, for very $(x,y)\notin N_1$,
$$
\big\|t\partial _{x_i}P_t^\mathcal{A}(x,y)\big\|_{L^2((0,\infty ),\frac{dt}{t})}\leq C\left\{
\begin{array}{ll}
e^{-\eta |x|^2},&\mbox{ if }\langle x,y\rangle\leq 0, \\[0.2cm]
\displaystyle |x+y|^n\exp\Big(-\eta\big(\frac{|x+y||x-y|}{2}-\frac{|y|^2-|x|^2}{2}\big)\Big),&\mbox{ if }\langle x,y\rangle\geq 0.
\end{array}
\right.
$$
By choosing $\frac{1}{p}<\eta<1$ we can deduce, as before, that $g_{i,\{P_t^\mathcal{A}\}_{t>0};\,{\rm glob}}^{X}$ is bounded from $L^p_X(\mathbb{R}^n,\gamma_{-1})$ into itself.

On the other hand, as above, we have that
$$
\big\|t\partial_{x_i}[P_t^\mathcal{A}(x,y)-P_t(x-y)]\big\|_{L^2((0,\infty ),\frac{dt}{t})}\leq C\int_0^\infty |\partial _{x_i}[T_u^\mathcal{A}(x,y)-W_u(x-y)]|\frac{du}{\sqrt{u}},\quad x,y\in \mathbb{R}^n.
$$ 
We can write
\begin{align*}
    \partial _{x_i}(T_u^\mathcal{A}(x,y)-W_u(x-y)) & =-\frac{2(x_i-y_ie^{-u})}{1-e^{-u}}T_u^{\mathcal A}(x,y)+\frac{x_i-y_i}{u}W_u(x-y)  \\
    & \hspace{-4cm} = \frac{x_i-y_i}{u}(W_u(x-y)-T_u^{\mathcal A}(x,y))+\left(\frac{x_i-y_i}{u}-\frac{2(x_i-y_ie^{-u})}{1-e^{-2u}}\right)T_u^{\mathcal A}(x,y),\quad x,y\in \mathbb{R}^n,\;u>0.
\end{align*}
Then, we get
$$
|\partial _{x_i}(T_u^\mathcal{A}(x,y)-W_u(x-y))|\leq C\left(\frac{e^{-nu}}{(1-e^{-2u})^{\frac{n+1}{2}}}+\frac{1}{u^{\frac{n+1}{2}}}\right)\leq \frac{C}{u^{\frac{n+1}{2}}},\quad x,y\in \mathbb{R}^n,\;u>0.
$$
Also, by using that
\begin{align*}
\left|\frac{x_i-y_i}{u}-\frac{2(x_i-y_ie^{-u})}{1-e^{-2u}}\right|&\leq C\left( |x-y|\left|\frac{1}{2u}-\frac{1}{1-e^{-2u}}\right|+\frac{1}{1-e^{-2u}}|x_i-y_i-(x_i-e^{-u}y_i)|\right)\\
&\leq C(|x-y|+|y|),\quad x,y\in \mathbb{R}^n,\;u>0,
\end{align*}
the estimation
$$
e^{-\frac{|x-e^{-u}|^2}{1-e^{-2u}}}\leq Ce^{-c\frac{|x-y|^2}{u}},\quad (x,y)\in N_1,\,u\in (0,1),
$$
and \eqref{3.2.2}, we obtain that, when $(x,y)\in N_1$ and $u\in (0,1)$, 
\begin{align*}
  |\partial _{x_i}(T_u^\mathcal{A}(x,y)-W_u(x-y)) |&\leq C e^{-c\frac{|x-y|^2}{u}}\left(\frac{|x-y|}{u^{\frac{n}{2}}}\Big(1+|y|^2+\frac{|y|}{\sqrt{u}}\Big)+\frac{e^{-nu}(|x-y|+|y|)}{(1-e^{-2u})^{\frac{n}{2}}}\right)\\
  &\hspace{-2cm} \leq C\frac{e^{-c\frac{|x-y|^2}{u}}}{u^{\frac{n-1}{2}}}\Big(1+|y|^2+\frac{|y|}{\sqrt{u}}\Big).
\end{align*}
Taking into account \eqref{K12infty}, \eqref{3.2.3} and \eqref{3.2.4} for $q=1$, we get 
\begin{align*}
    \big\|t\partial_{x_i}[P_t^\mathcal{A}(x,y)-P_t(x-y)]\big\|_{L^2((0,\infty ),\frac{dt}{t})}&\\
    &\hspace{-5cm}\leq C\left(\int_0^{m(x)}\frac{e^{-c\frac{|x-y|^2}{u}}}{u^{\frac{n}{2}}}\Big(1+|y|^2+\frac{|y|}{\sqrt{u}}\Big)du+\int_{m(x)}^\infty\frac{du}{u^{\frac{n}{2}+1}} \right)\leq C\frac{\sqrt{1+|x|}}{|x-y|^{n-\frac{1}{2}}},\quad (x,y)\in N_1,
\end{align*}
and deduce that
$$
\big\|g_{i,\{P_t^\mathcal{A}\}_{t>0};\,{\rm loc}}^{X}(f(\cdot,w))(x)-g_{i,\{P_t\}_{t>0};\,{\rm loc}}^{X}(f(\cdot,w))(x)\big\|_{L^p_X(\mathbb{R}^n,\gamma_{-1})}\leq C\|f\|_{L^p_X(\mathbb{R}^n,\gamma_{-1})},
$$
for every $f\in L^p_X(\mathbb{R}^n\gamma_{-1})$.

We now study $g_{i,\{P_t\}_{t>0};\,{\rm loc}}^{X}$ by using the vector-valued Calder\'on-Zygmund theory. We consider the operator $G_i$ defined by
$$
G_i(f)(t,x,w)=t\partial_{x_i}P_t(f(\cdot ,w))(x),\quad x\in \mathbb{R}^n,\;t>0,\;w\in \Omega .
$$
This is a convolution operator defined by the kernel
$$
K_t^i(x)=c_n\frac{x_it^2}{(t^2+2|x|^2)^{\frac{n+3}{2}}},\quad x\in \mathbb{R}^n,\;t>0,
$$
where $c_n=-2^{\frac{n}{2}+1}(n+1)\pi^{-\frac{n+1}{2}}\Gamma(\frac{n+1}{2})$.
The Fourier transforms $\widehat{K_t^i}$ of $K_t^i$ is given by
$$
\widehat{K_t^i}(y)=Ciy_ie^{-ct|y|},\quad y\in \mathbb{R}^n,\;t>0.
$$
Then,
$$
\big\|\widehat{K_t^i}(y)\big\|_{L^2((0,\infty ),\frac{dt}{t})}
=C\left(\int_0^\infty y_i^2t^2e^{-2ct|y|}\frac{dt}{t}\right)^{\frac{1}{2}}\leq C,\quad y\in \mathbb{R}^n.
$$
We also have that
\begin{align*}
\big\|K_t^i(x)\big\|_{L^2((0,\infty ),\frac{dt}{t})}
&=c_n\left(\int_0^\infty \frac{x_i^2t^3}{(t^2+2|x|^2)^{n+3}}dt\right)^{\frac{1}{2}}\\
&\leq C\left(\int_0^\infty\frac{dt}{(t+|x|)^{2n+1}}\right)^{\frac{1}{2}}\leq \frac{C}{|x|^n},\quad x\in \mathbb{R}^n\setminus\{0\}.
\end{align*}
In a similar way we can see that, for $j=1,...,n$,
$$
\big\|\partial_{x_j}K_t^i(x)\big\|_{L^2((0,\infty ),\frac{dt}{t})}\leq \frac{C}{|x|^{n+1}},\quad x\in \mathbb{R}^n\setminus\{0\}.
$$
By using Bourgain's results (see \cite{Bou} and \cite{Bou1}) we conclude that the operator $G_i$ is bounded from $L^p_X(\mathbb{R}^n,dx)$ into $L^p_{X(L^2((0,\infty ),\frac{dt}{t}))}(\mathbb{R}^n,dx)$.

Here $X(L^2((0,\infty ),\frac{dt}{t}))$ represents the K\"othe Bochner function space consisting of all the measurable functions $g:\Omega\longrightarrow L^2((0,\infty ),\frac{dt}{t})$ such that the function 
$
h(w)=\big\|g(w)\big\|_{L^2((0,\infty ),\frac{dt}{t})}$, $w\in \Omega,
$
belongs to $X$. $X(L^2((0,\infty ),\frac{dt}{t}))$ is endowed with the natural norm $\|\cdot\|_{X(L^2((0,\infty ),\frac{dt}{t}))}$ defined by
$$
\|g\|_{X(L^2((0,\infty ),\frac{dt}{t}))}=\big\|\big\|g(\cdot,w)\big\|_{L^2((0,\infty ),\frac{dt}{t})}\big\|_X,\quad g\in X(L^2((0,\infty ),\frac{dt}{t})).
$$
Hence, the operator $g_{i,\{P_t\}_{t>0}}^{X}$ is bounded from $L^p_X(\mathbb{R}^n, dx)$ into itself.

According to \cite[Propositions 3.2.5 and 3.2.7]{Sa}, $g_{i,\{P_t\}_{t>0};\,{\rm loc}}^{X}$ is bounded from $L^p_X(\mathbb{R}^n,\gamma_{-1})$ into itself.

We conclude that $g_{i,\{P_t^\mathcal{A}\}_{t>0}}^{X}$ is bounded from $L^p_X(\mathbb{R}^n,\gamma_{-1})$ into itself. Thus the property $(ii)$ is established.

\bibliographystyle{acm}

\end{document}